\documentclass[a4paper,reqno]{amsart}
\usepackage[all]{xy}           
\usepackage{amssymb}           
\usepackage{hyperref}
\usepackage{eucal}
\usepackage{graphicx}
\usepackage{epsfig}
\usepackage{wrapfig}

\numberwithin{equation}{section}

\newtheorem{definition}{Definition}[section]
\newtheorem{lemma}[definition]{Lemma}
\newtheorem{theorem}[definition]{Theorem}
\newtheorem{proposition}[definition]{Proposition}
\newtheorem{corollary}[definition]{Corollary}
\newtheorem{remarkth}[definition]{Remark}
\newtheorem{example}[definition]{Example}
\newenvironment{remark}{\begin{remarkth}\upshape}{\hfill$\diamond$\end{remarkth}}

\renewcommand{\emph}[1]{{\bfseries\itshape{#1}}}

\newcommand{\R}{\mathbb{R}}      
\newcommand{\N}{\mathbb{N}}      
\newcommand{\F}{\mathbb{F}}
\newcommand{\J}{\mathbb{J}}
\newcommand{\T}{\mathbb{T}}

\newcommand{\lcf}{\lbrack\! \lbrack}
\newcommand{\rcf}{\rbrack\! \rbrack}
\newcommand\map[3]{#1\ \colon\ #2\longrightarrow#3}
\newcommand{\lvec}[1]{\overleftarrow{#1}}
\newcommand{\rvec}[1]{\overrightarrow{#1}}

\newcommand{\I}{I\mkern-7muI}
\newcommand{\tr}{\operatorname{Tr}}
\newcommand{\e}{\mathrm{e}}

\makeatletter
\newcommand\prol{\@ifstar{\@proldf}{\@prolpf}}  
\def\@prolpf{\@ifnextchar[{\@prolpf@wrt}{\@prolpf@}}
\def\@prolpf@wrt[#1]#2{\@ifnextchar[{\@prolpf@wrt@at{#1}{#2}}{\@prolpf@wrt@{#1}{#2}}}
\def\@prolpf@wrt@at#1#2[#3]{\prolsymbol^{#1}_{#3}#2}
\def\@prolpf@wrt@#1#2{\prolsymbol^{#1}#2}
\def\@prolpf@#1{\@ifnextchar[{\@prolpf@at{#1}}{\@prolpf@@{#1}}}
\def\@prolpf@at#1[#2]{\prolsymbol_{#2}#1}
\def\@prolpf@@#1{\prolsymbol#1}
\def\@proldf{\@ifnextchar[{\@proldf@wrt}{\@proldf@}}
\def\@proldf@wrt[#1]#2{\@ifnextchar[{\@proldf@wrt@at{#1}{#2}}{\@proldf@wrt@{#1}{#2}}}
\def\@proldf@wrt@at#1#2[#3]{\prolsymbol^{*#1}_{#3}#2}
\def\@proldf@wrt@#1#2{\prolsymbol^{*#1}#2}
\def\@proldf@#1{\@ifnextchar[{\@proldf@at{#1}}{\@proldf@@{#1}}}
\def\@proldf@at#1[#2]{\prolsymbol^*_{#2}#1}
\def\@proldf@@#1{\prolsymbol^*#1}
\def\prolsymbol{\mathcal{T}}
\makeatother

\def\lcf{\lbrack\! \lbrack}
\def\rcf{\rbrack\! \rbrack}
\newcommand{\pai}[2]{\langle #1, #2\rangle}

\newcommand{\set}[2]{\left\{\,#1\left.\vphantom{#1#2}\,\right\vert\,#2\,
                \right\}}

\newcommand{\Sec}[1]{\operatorname{Sec}(#1)}

\begin{document}

\title[Discrete Nonholonomic Mechanics]{Discrete Nonholonomic Lagrangian Systems on Lie Groupoids}

\author[D.\ Iglesias]{David Iglesias}
\address{D.\ Iglesias:
Instituto de Matem\'aticas y F{\'\i}sica Fundamental, Consejo
Superior de Investigaciones Cient\'{\i}ficas, Serrano 123, 28006
Madrid, Spain} \email{iglesias@imaff.cfmac.csic.es}

\author[J.\ C.\ Marrero]{Juan C.\ Marrero}
\address{Juan C.\ Marrero:
Departamento de Matem\'atica Fundamental, Facultad de
Ma\-te\-m\'a\-ti\-cas, Universidad de la Laguna, La Laguna,
Tenerife, Canary Islands, Spain} \email{jcmarrer@ull.es}

\author[D.\ Mart\'{\i}n de Diego]{David Mart\'{\i}n de Diego}
\address{D.\ Mart\'{\i}n de Diego:
Instituto de Matem\'aticas y F{\'\i}sica Fundamental, Consejo
Superior de Investigaciones Cient\'{\i}ficas, Serrano 123, 28006
Madrid, Spain} \email{d.martin@imaff.cfmac.csic.es}

\author[E.\ Mart\'{\i}nez]{Eduardo Mart\'{\i}nez}
\address{Eduardo Mart\'{\i}nez:
Departamento de Matem\'atica Aplicada, Facultad de Ciencias,
Universidad de Zaragoza, 50009 Zaragoza, Spain}
\email{emf@unizar.es}

\thanks{This work has been partially supported by MICYT (Spain)
Grants MTM 2006-03322, MTM 2004-7832, MTM 2006-10531 and
S-0505/ESP/0158 of the CAM. D. Iglesias thanks MEC for a
``Juan de la Cierva" research contract.}

\keywords{}

\subjclass[2000]{}

\begin{abstract}
This paper studies the construction of geometric integrators for
nonholonomic systems. We derive the nonholonomic discrete
Euler-Lagrange equations in a setting which permits to deduce
geometric integrators  for continuous nonholonomic systems
(reduced or not). The formalism is given  in terms of Lie
groupoids, specifying a discrete Lagrangian  and a constraint
submanifold on it. Additionally, it is necessary to fix a vector
subbundle of the Lie algebroid associated to the Lie groupoid. We
also discuss the existence of nonholonomic evolution operators in
terms of the discrete nonholonomic Legendre transformations and in
terms of adequate decompositions of the prolongation of the Lie
groupoid. The characterization of the reversibility of the
evolution operator and the discrete nonholonomic momentum equation
are also considered. Finally, we illustrate with  several
classical examples the wide range of application of the theory
(the discrete nonholonomic constrained particle, the Suslov
system, the Chaplygin sleigh, the Veselova system, the rolling
ball on a rotating table and the two wheeled planar mobile robot).
\end{abstract}

\maketitle

\tableofcontents

\section{Introduction}

In the paper of  Moser and Veselov \cite{Mose} dedicated to the
complete integrability of certain dynamical systems, the authors
proposed a discretization of the tangent bundle $TQ$ of a
configuration space $Q$ replacing it by the product $Q\times Q$,
approximating a tangent vector on $Q$ by a pair of `close' points
$(q_0, q_1)$. In this sense, the continuous Lagrangian function
$L: TQ\longrightarrow \R$ is replaced by a discretization $L_d:
Q\times Q\longrightarrow \R$. Then,  applying a suitable
variational principle, it is possible to derive the discrete
equations of motion. In the regular case, one obtains an evolution
operator, a map which assigns to each pair $(q_{k-1}, q_k)$ a pair
$(q_k, q_{k+1})$, sharing many properties with the continuous
system, in particular, symplecticity, momentum conservation and a
good energy behavior. We refer to \cite{mawest} for an excellent
review in discrete Mechanics (on $Q \times Q$) and its numerical
implementation.

On the other hand, in \cite{Mose,VeVe}, the authors also
considered discrete Lagrangians defined on a Lie group $G$ where
the evolution operator is given by a diffeomorphism of $G$.

All the above examples led to A.\ Weinstein~\cite{weinstein} to
study discrete mechanics on Lie groupoids. A Lie groupoid is a
geometric structure that includes as particular examples the case
of cartesian products $Q\times Q$ as well as Lie groups and other
examples as Atiyah or action Lie groupoids \cite{Mac}. In a recent
paper \cite{groupoid}, we studied discrete Lagrangian and
Hamiltonian Mechanics on Lie groupoids, deriving from a
variational principle the discrete Euler-Lagrange equations. We
also introduced a symplectic 2-section (which is preserved by the
Lagrange evolution operator) and defined the Hamiltonian evolution
operator,  in terms of the discrete Legendre transformations,
which is a symplectic map with respect to the canonical symplectic
2-section on the prolongation of the dual of the Lie algebroid of
the given groupoid. These techniques include as particular cases
the classical discrete Euler-Lagrange equations, the discrete
Euler-Poincar\'e equations (see \cite{Bobe,Bobe1,Mars3,MaPeSh})
and the discrete Lagrange-Poincar\'e equations. In fact, the
results in \cite{groupoid} may be applied in the construction of
geometric integrators for continuous Lagrangian systems which are
invariant under the action of a symmetry Lie group (see also
\cite{JaLeMaWe} for the particular case when the symmetry Lie
group is abelian).

From the perspective of geometric integration, there are a great
interest in introducing new geometric techniques for developing
numerical integrators since standard methods often introduce some
spurious effects like dissipation in conservative systems
\cite{Hair,Sanz}. The case of dynamical systems subjected to
constraints is also of considerable interest. In particular, the
case of holonomic constraints is well established in the
literature of geometric integration, for instance, in simulation
of molecular dynamics
 where the constraints may be molecular bond lengths or
angles and also in multibody dynamics (see \cite{Hair,LeRe} and
references therein).

By contrast, the construction of geometric integrators  for the
case of  nonholonomic constraints is less well understood. This
type of constraints appears, for instance, in mechanical models of
convex rigid bodies rolling without sliding on a surface
\cite{NF}. The study of systems with nonholonomic constraints goes
back to the XIX century. The equations of motion were obtained
applying either D'Alembert's principle of virtual work or Gauss
principle of least constraint. Recently, many authors have shown a
new interest in that theory and also in its relation  to the new
developments in control theory and  robotics using geometric
techniques (see, for instance,
\cite{BS,Bl,BlKrMaMu,Cort,Koiller,LeMa2,LD}).

Geometrically, nonholonomic constraints are globally described by
a submanifold ${\mathcal M}$ of the velocity phase space $TQ$. If
${\mathcal M}$ is a vector subbundle of $TQ$, we are dealing with
the case of linear constraints and, in the case ${\mathcal M}$ is
an affine subbundle, we are in the case of affine constraints.
Lagrange-D'Alembert's  or Chetaev's principles allow us to
determine the set of possible values of the constraint forces only
from the set of admissible kinematic states, that is, from the
constraint manifold ${\mathcal M}$ determined by the vanishing of
the nonholonomic constraints $\phi^a$. Therefore,  assuming that
the dynamical properties of the system are mathematically
described by a Lagrangian function $L: TQ\longrightarrow \R$ and
by a constraint submanifold ${\mathcal M}$,  the equations of
motion, following Chetaev's principle, are
\[\left[ \frac{d}{dt}\left( \frac{\partial L}{\partial \dot
q^i}\right) - \frac{\partial L}{\partial q^i} \right] \delta
q^i=0\; ,
\]
where $\delta q^i$ denotes the virtual displacements verifying
$\displaystyle{ \frac{\partial \phi^a}{\partial \dot{q}^i}\delta
q^i =0 }$.
 By using the Lagrange multiplier rule,
we obtain that
\begin{equation}\label{aqq}
\frac{d}{dt}\left( \frac{\partial L}{\partial \dot
q^i}\right)-\frac{\partial L}{\partial
q^i}=\bar{\lambda}_a\frac{\partial \phi^i}{\partial \dot{q}^i} \;
,
\end{equation}
with the condition $\dot{q}(t)\in {\mathcal M}$,
$\bar{\lambda}_{a}$ being the Lagrange multipliers to be
determined. Recently, J. Cort\'es {\em et al} \cite{CoLeMaMa} (see
also \cite{CoMa,Mest,MeLa}) proposed a unified framework for
nonholonomic systems in the Lie algebroid setting that we will use
along this paper generalizing some previous work for free
Lagrangian mechanics on Lie algebroids (see, for instance,
\cite{LeMaMa,M,Ma2,Ma4}).

The  construction of geometric integrators for Equations
(\ref{aqq}) is very recent. In fact, in  \cite{MS} appears as an
open problem:
\begin{quote}
...The problem for the more general class of non-holonomic
constraints is still open, as is the question of the correct
analogue of symplectic integration for non-holonomically
constrained Lagrangian systems...
\end{quote}
Numerical integrators derived from discrete variational principles
have proved their adaptability to  many situations: collisions,
classical field theory, external forces...\cite{Mars1,mawest} and
it also seems very adequate for nonholonomic systems, since
nonholonomic equations of motion come  from H\"{o}lder's variational
principle which  is not a standard variational principle
\cite{Arno}, but admits an adequate discretization. This is the
procedure introduced   by J. Cort\'es and S. Mart{\'\i}nez
\cite{Cort,CoSMa} and followed by other authors
\cite{F,FZ,FeZe,McLPerl} extending, moreover, the results to
nonholonomic systems defined on Lie groups (see also \cite{LeMaSa}
for a different approach using generating functions).

In this paper, we tackle the problem from the unifying point of
view of Lie groupoids (see \cite{CoLeMaMa} for the continuous
case). This technique permits to recover all the previous methods
in the literature \cite{CoSMa,FZ,McLPerl} and consider new cases
of great importance in nonholonomic dynamics. For instance, using
action Lie groupoids, we may discretize LR-nonholonomic systems
such as the Veselova system or using Atiyah Lie groupoids we find
discrete versions for the reduced equations of  nonholonomic
systems with symmetry.

The paper is structured as follows. In section 2 we review some
basic results on Lie algebroids and Lie groupoids. In particular,
we describe the prolongation of a Lie groupoid \cite{Sau}, which
has a double structure of Lie groupoid and Lie algebroid. Then, we
briefly expose the geometric structure of discrete unconstrained
mechanics on Lie groupoids: Poincar\'e-Cartan sections, Legendre
transformations... The main results of the paper appear in section
3, where the geometric structure of discrete nonholonomic systems
on Lie groupoids is considered. In particular, given a discrete
Lagrangian $L_d:\Gamma\to \R$ on a Lie groupoid $\Gamma$, a
constraint distribution ${\mathcal D}_c$ in the Lie algebroid
$E_\Gamma$ of $\Gamma$ and a discrete constraint submanifold
${\mathcal M}_c$ in $\Gamma$, we obtain the nonholonomic discrete
Euler-Lagrange equations from a discrete Generalized H\"older's
principle (see section 3.1). In addition, we characterize the
regularity of the nonholonomic system in terms of the nonholonomic
Legendre transformations and decompositions of the prolongation of
the Lie groupoid. In the case when the system is regular, we can
define the nonholonomic evolution operator. An interesting
situation, studied in in Section 3.4, is that of reversible
discrete nonholonomic Lagrangian systems, where the Lagrangian and
the discrete constraint submanifold are invariants with respect to
the inversion of the Lie groupoid. The particular example of
reversible systems in the pair groupoid $Q\times Q$ was first
studied in \cite{McLPerl}. We also define the discrete
nonholonomic momentum map. In order to give an idea of the breadth
and flexibility of the proposed formalism, several examples are
discussed, including their regularity and their reversibility:
\begin{itemize}
\item[-] Discrete holonomic Lagrangian systems on a Lie groupoid,
which are a generalization of the Shake algorithm for holonomic
systems \cite{Hair,LeRe,mawest}; \item[-] Discrete nonholonomic
systems on the pair groupoid, recovering the equations first
considered in \cite{CoSMa}. An explicit example of this situation
is the discrete nonholonomic constrained particle. \item[-]
Discrete nonholonomic systems on Lie groups, where the equations
that are obtained are the so-called discrete
Euler-Poincar\'e-Suslov equations (see \cite{FZ}). We remark that,
although our equations coincide with those in \cite{FZ}, the
technique developed in this paper is different to the one in that
paper. Two explicit examples which we describe here are the Suslov
system and the Chaplygin sleigh. \item[-] Discrete nonholonomic
Lagrangian systems on an action Lie groupoid. This example is
quite interesting since it allows us to discretize a well-known
nonholonomic LR-system: the Veselova system (see \cite{VeVe}; see
also \cite{FeJo}). For this example, we obtain a discrete system
that is not reversible and we show that the system is regular in a
neighborhood  around the manifold of units. \item[-] Discrete
nonholonomic Lagrangian systems on an Atiyah Lie groupoid. With
this example, we are able to discretize reduced systems, in
particular, we concentrate on the example of the discretization of
the equations of motion of a rolling ball without sliding on a
rotating table with constant angular velocity. \item[-] Discrete
Chaplygin systems, which are regular systems $(L_d, {\mathcal
M}_c,{\mathcal D}_c)$ on the Lie groupoid $\Gamma
\rightrightarrows M$, for which $(\alpha ,\beta )\circ
i_{{\mathcal M}_c}:{\mathcal M}_c\to M\times M$ is a
diffeomorphism and $\rho \circ i_{{\mathcal D}_c}:{\mathcal D}_c
\to TM$ is an isomorphism of vector bundles, $(\alpha, \beta)$
being the source and target of the Lie groupoid $\Gamma$ and
$\rho$ being the anchor map of the Lie algebroid $E_\Gamma$. This
example includes a discretization of the two wheeled planar mobile
robot.
\end{itemize}
We conclude our paper with future lines of work.

\section{Discrete Unconstrained Lagrangian Systems on Lie Groupoids}
\subsection{Lie algebroids}
A \emph{Lie algebroid} $E$ over a manifold $M$ is a real vector
bundle $\tau: E \to M$ together with a Lie bracket $\lcf\cdot,
\cdot\rcf$ on the space $\Sec{\tau}$ of the global cross sections
of $\tau: E \to M$ and a bundle map $\rho:E \to TM$, called
\emph{the anchor map}, such that if we also denote by
$\rho:\Sec{\tau}\to {\frak X}(M)$ the homomorphism of
$C^\infty(M)$-modules induced by the anchor map then
\begin{equation}\label{Leib}
 \lcf X,fY\rcf=f\lcf X,Y\rcf + \rho(X)(f)Y,
\end{equation}
for $X,Y\in \Sec{\tau}$ and $f\in C^\infty(M)$ (see \cite{Mac}).

 If $(E,\lcf\cdot,\cdot\rcf,\rho)$ is a Lie algebroid over $M$
 then the anchor map  $\rho:\Sec{\tau}\to {\frak X}(M)$ is a
 homomorphism between the Lie algebras
 $(\Sec{\tau},\lcf\cdot,\cdot\rcf)$ and $({\frak
 X}(M),[\cdot,\cdot])$. Moreover, one may define the differential
 $d$ of $E$ as follows:
 \begin{equation}\label{dif}
\begin{aligned}
d \mu(X_0,\dots, X_k)&=\sum_{i=0}^{k}
(-1)^i\rho(X_i)(\mu(X_0,\dots,
\widehat{X_i},\dots, X_k)) \\
&+ \sum_{i<j}(-1)^{i+j}\mu(\lcf X_i,X_j\rcf,X_0,\dots,
\widehat{X_i},\dots,\widehat{X_j},\dots ,X_k),
\end{aligned}
\end{equation}
for $\mu\in \Sec{\wedge^k \tau^*}$ and $X_0,\dots ,X_k\in
\Sec{\tau}.$ $d$ is a cohomology operator, that is, $d^2=0$. In
particular, if $f: M\longrightarrow \R$ is a real smooth function
then $df(X)=\rho(X)f,$ for $X\in \Sec{\tau}$.

Trivial examples of Lie algebroids are a real Lie algebra of
finite dimension (in this case, the base space is a single point)
and the tangent bundle of a manifold $M.$

On the other hand, let $(E,\lcf\cdot,\cdot\rcf,\rho)$ be a Lie
algebroid of rank $n$  over a manifold $M$ of dimension $m$ and
$\pi:P\to M$ be a fibration. We consider the subset of $E\times
TP$
\[
{\mathcal T}^EP=\set{(a,v)\in E\times TP}{(T\pi)(v)=\rho(a)},
\]
where $T\pi: TP \to TM$ is the tangent map to $\pi$. Denote by
$\tau^{\pi}: {\mathcal T}^EP \to P$ the map given by
$\tau^{\pi}(a, v) = \tau_{P}(v)$, $\tau_{P}: TP \to P$ being the
canonical projection. If $dim P = p$, one may prove that
${\mathcal T}^EP$ is a vector bundle over $P$ of rank $n + p -m$
with vector bundle projection $\tau^{\pi}: {\mathcal T}^EP \to P$.

A section $\tilde{X}$ of $\tau^{\pi}: {\mathcal T}^EP \to P$ is
said to be \emph{projectable} if there exists a section $X$ of
$\tau: E \to M$ and a vector field $U \in {\frak X}(P)$ which is
$\pi$-projectable to the vector field $\rho(X)$ and such that
$\tilde{X}(p) = (X(\pi(p)), U(p))$, for all $p \in P$. For such a
projectable section $\tilde{X}$, we will use the following
notation $\tilde{X} \equiv (X, U)$. It is easy to prove that one
may choose a local basis of projectable sections of the space
$\Sec{\tau^{\pi}}$.

The vector bundle $\tau^{\pi}: {\mathcal T}^EP \to P$ admits a Lie
algebroid structure $(\lcf \cdot , \cdot \rcf^{\pi}, \rho^{\pi})$.
Indeed, if $(X, U)$ and $(Y, V)$ are projectable sections then
\[
\lcf (X, U), (Y, V) \rcf^{\pi} = (\lcf X, Y\rcf, [U, V]),
\makebox[.4cm]{}  \rho^{\pi}(X, U) = U.
\]

$({\mathcal T}^EP,\lcf\cdot,\cdot\rcf^\pi,\rho^\pi)$ is the
\emph{$E$-tangent bundle to $P$ or the prolongation of $E$ over
the fibration $\pi:P\to M$} (for more details, see \cite{LeMaMa}).

Now, let $(E,\lcf\cdot,\cdot \rcf,\rho)$ (resp.,
$(E',\lcf\cdot,\cdot\rcf', \rho')$) be a Lie algebroid over a
manifold $M$ (resp., $M'$) and suppose that $\Psi: E \to E'$ is a
vector bundle morphism over the map $\Psi_{0}: M \to M'$. Then,
the pair $(\Psi, \Psi_{0})$ is said to be a \emph{Lie algebroid
morphism} if
\begin{equation}\label{dd'}
d ((\Psi, \Psi_{0})^*\phi')= (\Psi, \Psi_{0})^*(d'\phi'), \;\;\;
\text{ for all }\phi'\in \Sec{\wedge^k(\tau')^*} \text{ and for all
}k,
\end{equation}
where $d$ (resp., $d'$) is the differential of the Lie algebroid
$E$ (resp., $E'$) (see \cite{LeMaMa}). In the particular case when
$M = M'$ and $\Psi_{0} = Id$ then (\ref{dd'}) holds if and only if
\[
\lcf \Psi \circ X, \Psi \circ Y \rcf' = \Psi \lcf X, Y \rcf,
\makebox[.3cm]{} \rho'(\Psi X) = \rho(X), \makebox[.3cm]{} \mbox{
for } X, Y \in \Sec{\tau}.
\]

\subsection{Lie groupoids}\label{secLg}

A Lie groupoid over a differentiable manifold $M$ is a
differentiable manifold $\Gamma$ together with the following
structural maps:
\begin{itemize}
\item A pair of submersions $\alpha:\Gamma\to M$, \emph{the
source}, and $\beta:\Gamma\to M,$ \emph{the target}. The maps
$\alpha$ and $\beta$ define the set of \emph{composable pairs }
\[
\Gamma_2=\set{(g,h)\in G\times G}{\beta(g)=\alpha(h)}.
\]
\item A \emph{multiplication} $m: \Gamma_{2} \to \Gamma$, to be
denoted simply by $m(g,h)=gh$, such that
\begin{itemize}
\item $\alpha(gh)=\alpha(g)$ and $\beta(gh)=\beta(h)$. \item
$g(hk)=(gh)k$.
\end{itemize}

\item An \emph{identity section} $\epsilon: M \to \Gamma$ such
that
\begin{itemize}
\item $\epsilon(\alpha(g))g=g$ and $g\epsilon(\beta(g))=g$.
\end{itemize}
\item An \emph{inversion map} $i: \Gamma \to \Gamma$, to be
simply denoted by $i(g)=g^{-1}$, such that
\begin{itemize}
\item $g^{-1}g=\epsilon(\beta(g))$ and
$gg^{-1}=\epsilon(\alpha(g))$.
\end{itemize}
\end{itemize}
A Lie groupoid $\Gamma$ over a set $M$ will be simply denoted by
the symbol $\Gamma \rightrightarrows M$.

On the other hand, if $g \in \Gamma$ then the
\emph{left-translation} by $g$ and the
\emph{right-translation} by $g$ are the diffeomorphisms
$$
\begin{array}{lll}
l_{g}: \alpha^{-1}(\beta(g)) \longrightarrow
\alpha^{-1}(\alpha(g))&; \; \;& h \longrightarrow
l_{g}(h) = gh, \\
r_{g}: \beta^{-1}(\alpha(g)) \longrightarrow
\beta^{-1}(\beta(g))&; \; \;& h \longrightarrow r_{g}(h) = hg.
\end{array}
$$
Note that $l_{g}^{-1} = l_{g^{-1}}$ and $r_{g}^{-1} = r_{g^{-1}}$.

A vector field $\tilde{X}$ on $\Gamma$ is said to be
\emph{left-invariant} (resp., \emph{right-invariant}) if it is
tangent to the fibers of $\alpha$ (resp., $\beta$) and
$\tilde{X}(gh) = (T_{h}l_{g})(\tilde{X}_{h})$ (resp.,
$\tilde{X}(gh) = (T_{g}r_{h})(\tilde{X}(g)))$, for $(g,h) \in
\Gamma_{2}$.

Now, we will recall the definition of the \emph{Lie algebroid
associated with $\Gamma$}.

We consider the vector bundle $\tau: E_\Gamma \to M$, whose fiber
at a point $x \in M$ is $(E_{\Gamma})_{x} = V_{\epsilon(x)}\alpha
= Ker (T_{\epsilon(x)}\alpha)$. It is easy to prove that there
exists a bijection between the space $\Sec{\tau}$ and the set of
left-invariant (resp., right-invariant) vector fields on $\Gamma$.
If $X$ is a section of $\tau: E_{\Gamma} \to M$, the corresponding
left-invariant (resp., right-invariant) vector field on $\Gamma$
will be denoted $\lvec{X}$ (resp., $\rvec{X}$), where
\begin{equation}\label{linv}
\lvec{X}(g) = (T_{\epsilon(\beta(g))}l_{g})(X(\beta(g))),
\end{equation}
\begin{equation}\label{rinv}
\rvec{X}(g) = -(T_{\epsilon(\alpha(g))}r_{g})((T_{\epsilon
(\alpha(g))}i)( X(\alpha(g)))),
\end{equation}
for $g \in \Gamma$. Using the above facts, we may introduce a Lie
algebroid structure $(\lcf\cdot , \cdot\rcf, \rho)$ on
$E_{\Gamma}$, which is defined by
\begin{equation}\label{LA}
\lvec{\lcf X, Y\rcf} = [\lvec{X}, \lvec{Y}], \makebox[.3cm]{}
\rho(X)(x) = (T_{\epsilon(x)}\beta)(X(x)),
\end{equation}
for $X, Y \in \Sec{\tau}$ and $x \in M$. Note that
\begin{equation}\label{RL}
\rvec{\lcf X, Y\rcf} = -[\rvec{X}, \rvec{Y}], \makebox[.3cm]{}
[\rvec{X}, \lvec{Y}] = 0,
\end{equation}
(for more details, see \cite{CDW,Mac}).

Given two Lie groupoids $\Gamma \rightrightarrows M$ and $\Gamma'
\rightrightarrows M'$, a \emph{morphism of Lie groupoids} is a
smooth map $\Phi: \Gamma \to \Gamma'$ such that
\[
(g, h) \in \Gamma_{2} \Longrightarrow (\Phi(g), \Phi(h)) \in
(\Gamma')_{2}
\]
and
\[
\Phi(gh) = \Phi(g)\Phi(h).
\]
A morphism of Lie groupoids $\Phi: \Gamma \to \Gamma'$ induces a
smooth map $\Phi_{0}: M \to M'$ in such a way that
\[
\alpha' \circ \Phi = \Phi_{0} \circ \alpha, \makebox[.3cm]{}
\beta' \circ \Phi = \Phi_{0} \circ \beta, \makebox[.3cm]{} \Phi
\circ \epsilon = \epsilon' \circ \Phi_{0},
\]
$\alpha$, $\beta$ and $\epsilon$ (resp., $\alpha'$, $\beta'$ and
$\epsilon'$) being the source, the target and the identity section
of $\Gamma$ (resp., $\Gamma'$).

Suppose that $(\Phi, \Phi_{0})$ is a morphism between the Lie
groupoids $\Gamma \rightrightarrows M$ and $\Gamma'
\rightrightarrows M'$ and that $\tau: E_\Gamma \to M$ (resp.,
$\tau': E_{\Gamma'} \to M'$) is the Lie algebroid of $\Gamma$
(resp., $\Gamma'$). Then, if $x \in M$ we may consider the linear
map $E_{x}(\Phi): (E_\Gamma)_{x} \to (E_{\Gamma'})_{\Phi_{0}(x)}$
defined by
\begin{equation}\label{Amor}
E_{x}(\Phi)(v_{\epsilon(x)}) =
(T_{\epsilon(x)}\Phi)(v_{\epsilon(x)}), \; \; \mbox{ for }
v_{\epsilon(x)} \in (E_\Gamma)_{x}.
\end{equation}
In fact, we have that the pair $(E(\Phi), \Phi_{0})$ is a morphism
between the Lie algebroids $\tau: E_\Gamma \to M$ and $\tau':
E_{\Gamma'} \to M'$ (see \cite{Mac}).

Trivial examples of Lie groupoids are Lie groups and the pair or
banal groupoid $M\times M$, $M$ being an arbitrary smooth
manifold. The Lie algebroid of a Lie group $\Gamma$ is just the
Lie algebra ${\frak g}$ of $\Gamma$. On the other hand, the Lie
algebroid of the pair (or banal) groupoid $M\times M$ is the
tangent bundle $TM$ to $M$.

Apart from the Lie algebroid $E_\Gamma$ associated with a Lie
groupoid $\Gamma\rightrightarrows M$, other interesting Lie
algebroids associated with $\Gamma$ are the following ones:
\begin{itemize}
\item {\bf The $E_\Gamma$- tangent bundle to $E_\Gamma^*$:}
\end{itemize}
Let ${\mathcal T}^{E_\Gamma}E_{\Gamma}^{*}$ be the
$E_\Gamma$-tangent bundle to $E^*_\Gamma$, that is,
\[
{\mathcal
T}_{\Upsilon_x}^{E_\Gamma}E_\Gamma^*=\set{(v_x,X_{\Upsilon_x})\in
(E_\Gamma)_x\times T_{\Upsilon_x}
E_\Gamma^*}{(T_{\Upsilon_x}\tau^*)(X_{\Upsilon_x})=(T_{\epsilon(x)}\beta)(v_x)}
\]
for $\Upsilon_{x}\in (E_\Gamma^*)_x,$ with $x\in M.$ As we know,
${\mathcal T}^{E_\Gamma}E_\Gamma^*$ is a  Lie algebroid over
$E_\Gamma^*$.

We may  introduce the canonical section $\Theta$ of the vector
bundle $({\mathcal T}^{E_{\Gamma}}E_{\Gamma}^{*})^* \to
E_{\Gamma}^{*}$ as follows:
\[
\Theta(\Upsilon_{x})(a_{x}, X_{\Upsilon_{x}}) =
\Upsilon_{x}(a_{x}),
\]
for $\Upsilon_{x} \in (E_{\Gamma}^{*})_{x}$ and $(a_{x},
X_{\Upsilon_{x}}) \in {\mathcal
T}_{\Upsilon_x}^{E_{\Gamma}}E_{\Gamma}^{*}$. $\Theta$ is called
the \emph{Liouville section associated with $E_\Gamma$}. Moreover,
we define \emph{the canonical symplectic section} $\Omega$
associated with $E_{\Gamma}$ by $\Omega = -d\Theta$, where $d$ is
the differential on the Lie algebroid ${\mathcal
T}^{E_{\Gamma}}E_{\Gamma}^{*} \to E_{\Gamma}^{*}$. It is easy to
prove that $\Omega$ is nondegenerate and closed, that is, it is a
symplectic section of ${\mathcal T}^{E_{\Gamma}}E_{\Gamma}^{*}$
(see \cite{LeMaMa}).

Now, if $Z$ is a section of $\tau:{E_\Gamma}\to M$ then there is a
unique vector field $Z^{*c}$ on $E_\Gamma^*$, \emph{the complete
lift of $X$ to $E_\Gamma^*$}, satisfying the two following
conditions:
\begin{enumerate}
\item $Z^{*c}$ is $\tau^*$-projectable on $\rho(Z)$ and \item
$Z^{*c}(\widehat{X})=\widehat{\lcf Z,X \rcf}$ \end{enumerate} for
$X\in \Sec{\tau}$ (see \cite{LeMaMa}). Here, if $X$ is a section of
$\tau:E_\Gamma\to M$ then $\widehat{X}$ is the linear function
$\widehat{X}\in C^\infty(E^*)$ defined by
\[
\widehat{X}(a^*)=a^*(X(\tau^*(a^*))), \mbox{ for all } a^*\in E^*.
\]
Using the vector field $Z^{*c}$, one may introduce \emph{the
complete lift} $Z^{*\bf{c}}$ of $Z$ as the section of
$\tau^{\tau^*}:{\mathcal T}^{E_\Gamma}E_\Gamma^*\to E^*_{\Gamma}$
defined by
\begin{equation}\label{cl*1}
Z^{*\bf{c}}(a^*)=(Z(\tau^*(a^*)),Z^{*c}(a^*)),\mbox{ for } a^*\in
E^*.
\end{equation}

$Z^{*\bf{c}}$ is just the Hamiltonian section of $\widehat{Z}$
with respect to the canonical symplectic section $\Omega$
associated with $E_\Gamma$. In other words,
\begin{equation}\label{cl*}
i_{Z^{*\bf{c}}}\Omega=d\widehat{Z},
\end{equation}
where $d$ is the differential of the Lie algebroid
$\tau^{\tau^*}:{\mathcal T}^{E_\Gamma}E_\Gamma^*\to E_\Gamma^*$
(for more details, see \cite{LeMaMa}).

\begin{itemize}\item {\bf The Lie algebroid $\widetilde{\tau_\Gamma}:{\mathcal
T}^\Gamma\Gamma\to \Gamma$ }:\end{itemize}

Let ${\mathcal T}^\Gamma\Gamma$ be the Whitney sum
$V\beta\oplus_\Gamma V\alpha$ of the vector bundles $V\beta\to
\Gamma$ and $V\alpha\to \Gamma$, where $V\beta$ (respectively,
$V\alpha$) is the vertical bundle of $\beta$ (respectively,
$\alpha$). Then, the vector bundle
$\widetilde{\tau}_\Gamma:{\mathcal T}^\Gamma\Gamma\equiv
V\beta\oplus_\Gamma V\alpha\to \Gamma$ admits a Lie algebroid
structure $(\lcf\cdot,\cdot\rcf^{{\mathcal
T}^\Gamma\Gamma},\rho^{{\mathcal T}^\Gamma\Gamma})$. The anchor
map $\rho^{{\mathcal T}^\Gamma\Gamma}$ is given by

\[
(\rho^{{\mathcal T}^\Gamma\Gamma})(X_g,Y_g)=X_g+Y_g
\]
and the Lie bracket bracket $\lcf\cdot,\cdot\rcf^{{\mathcal
T}^\Gamma\Gamma}$ on the space $\Sec{\widetilde{\tau}_\Gamma}$  is
characterized for the following relation
\[
\lcf(\rvec{X},\lvec{Y}),(\rvec{X'},\lvec{Y'})\rcf^{{\mathcal
T}^\Gamma\Gamma}=(-\rvec{\lcf X,X'\rcf},\lvec{\lcf Y,Y'\rcf}),
\]
for $X,Y,X',Y'\in \Sec{\tau}$ (for more details, see
\cite{groupoid}).

On other hand, if $X$ is a section of $\tau:E_\Gamma\to M$, one
may define the sections $X^{(1,0)},X^{(0,1)}$ (the $\beta$ and
$\alpha$-lifts) and $X^{(1,1)}$ (the complete lift) of $X$ to
$\widetilde{\tau}_\Gamma:{\mathcal T}^\Gamma\Gamma\to \Gamma$ as
follows:
\[
X^{(1,0)}(g)=(\rvec{X}(g),0_g),\;\;\;
X^{(0,1)}(g)=(0_g,\lvec{X}(g)),\;\;\mbox{ and }
X^{(1,1)}(g)=(-\rvec{X}(g),\lvec{X}(g)).
\]
We have that
\[
\lcf X^{(1,0)},Y^{(1,0)}\rcf^{{\mathcal T}^\Gamma\Gamma}=-\lcf
X,Y\rcf^{(1,0)}\;\;\lcf X^{(0,1)},Y^{(1,0)}\rcf^{{\mathcal
T}^\Gamma\Gamma}=0,\]\[\;\;\;\lcf
X^{(0,1)},Y^{(0,1)}\rcf^{{\mathcal T}^\Gamma\Gamma}=\lcf
X,Y\rcf^{(0,1)},
\]
and, as a consequence,
\[
\lcf X^{(1,1)},Y^{(1,0)}\rcf^{{\mathcal T}^\Gamma\Gamma}=\lcf
X,Y\rcf^{(1,0)},\;\;\;\lcf X^{(1,1)},Y^{(0,1)}\rcf^{{\mathcal
T}^\Gamma\Gamma}=\lcf X,Y\rcf^{(0,1)},\]\[\;\;\;\lcf
X^{(1,1)},Y^{(1,1)}\rcf^{{\mathcal T}^\Gamma\Gamma}=\lcf
X,Y\rcf^{(1,1)}.
\]
 Now, if $g,h\in \Gamma$ one may introduce
the linear monomorphisms $\;_h^{(1,0)}:
(E_\Gamma)_{\alpha(h)}^*\to ({\mathcal T}_{h}^{\Gamma}\Gamma)^{*}
\equiv V_{h}^{*}\beta \oplus V_{h}^{*}\alpha$ and $\;_g^{(0,1)}:
(E_\Gamma)_{\beta(g)}^*\to ({\mathcal T}_{g}^{\Gamma}\Gamma)^{*}
\equiv V_{g}^{*}\beta \oplus V_{g}^{*}\alpha$ given by
\begin{equation}\label{(1,0)h}
\gamma_h^{(1,0)}(X_h,Y_h)=\gamma(T_h(i\circ r_{h^{-1}})(X_h)),
\end{equation}
\begin{equation}\label{(0,1)g}
\gamma_g^{(0,1)}(X_g,Y_g)=\gamma((T_gl_{g^{-1}})(Y_g)),
\end{equation}
for $(X_g,Y_g)\in {\mathcal T}_g^\Gamma\Gamma$ and $(X_h,Y_h)\in
{\mathcal T}_h^\Gamma\Gamma.$

Thus, if $\mu$ is a section of $\tau^*:E_\Gamma^*\to M$, one may
define the corresponding lifts $\mu^{(1,0)}$ and $\mu^{(0,1)}$ as
the sections of $\widetilde{\tau_\Gamma}^*: ({\mathcal
T}^\Gamma\Gamma)^{*} \to \Gamma$ given by
\[
\begin{array}{rr}
\mu^{(1,0)}(h)=\mu_h^{(1,0)},& \mbox{for } h\in \Gamma,\\
\mu^{(0,1)}(g)=\mu_g^{(0,1)},& \mbox{for } g\in \Gamma.
\end{array}
\]

Note that if $g\in \Gamma$ and $\{X_A\}$ (respectively, $\{Y_B\}$)
is a local basis of $\Sec{\tau}$ on an open subset $U$
(respectively, $V$) of $M$ such that $\alpha(g)\in U$
(respectively, $\beta(g)\in V$) then $\{X_A^{(1,0)},Y_B^{(0,1)}\}$
is a local basis of $\Sec{\widetilde{\tau}_\Gamma}$ on the open
subset $\alpha^{-1}(U)\cap \beta^{-1}(V).$ In addition, if
$\{X^A\}$ (respectively, $\{Y^B\}$) is the dual basis of $\{X_A\}$
(respectively, $\{Y_B\}$) then $\{(X^A)^{(1,0)},(Y^B)^{(0,1)}\}$
is the dual basis of $\{X_A^{(1,0)},Y_B^{(0,1)}\}$.

\subsection{Discrete Unconstrained Lagrangian Systems}
(See \cite{groupoid} for details)
 A \emph{discrete unconstrained Lagrangian system on a Lie groupoid}
consists of a Lie groupoid ${\Gamma}\rightrightarrows M$ (the
\emph{discrete
  space}) and a \emph{discrete Lagrangian} $L_d: {\Gamma} \to \R$.

\subsubsection{Discrete unconstrained Euler-Lagrange equations.}
 An \emph{admissible sequence of order $N$ } on the Lie groupoid
 $\Gamma$ is an element $(g_1,\dots ,g_N)$ of $\Gamma^N\equiv
 \Gamma\times  \dots \times\Gamma$ such that
 $(g_k,g_{k+1})\in \Gamma_2$, for $k=1,\dots, N-1$.

 An admissible sequence $(g_1,\dots ,g_N)$ of order $N$ is
 a solution of the \emph{ discrete unconstrained Euler-Lagrange equations } for
 $L_d$ if
 \[
 d^o[L_d\circ l_{g_k} + L_d\circ r_{g_{k+1}}\circ
 i](\epsilon(x_k))_{|(E_\Gamma)_{x_k}}=0
 \]
 where $\beta(g_k)=\alpha(g_{k+1})=x_k$ and $d^o$ is the standard
 differential on $\Gamma$, that is, the differential of the Lie
 algebroid $\tau_\Gamma:T\Gamma\to \Gamma$ (see \cite{groupoid}).

 The \emph{ discrete unconstrained Euler-Lagrange operator}
 $D_{DEL}L_d:\Gamma_2\to E_\Gamma^*$ is given by
 \[
(D_{DEL}L_d)(g,h)=d^o[L_d\circ l_g + L_d\circ r_h\circ
i](\epsilon(x))_{|(E_\Gamma)_x}=0,\] for $(g,h)\in \Gamma_2$, with
$\beta(g)=\alpha(h)=x\in M$ (see \cite{groupoid}).

Thus, an admissible sequence $(g_1,\dots ,g_N)$ of order $N$ is a
solution of the discrete unconstrained Euler-Lagrange equations if
and only if
\[
(D_{DEL}L_d)(g_k,g_{k+1})=0,\;\;\; \mbox{ for }k=1,\dots ,N-1.
\]

\subsubsection{ Discrete Poincar\'e-Cartan sections} Consider the
Lie algebroid $\widetilde{\tau_\Gamma}:{\mathcal
T}^\Gamma\Gamma\equiv V\beta\oplus_{\Gamma} V\alpha \to {\Gamma}$,
and define the \emph{Poincar\'e-Cartan 1-sections }
$\Theta_{L_d}^-, \Theta_{L_d}^+\in \Sec{
(\widetilde{\tau_\Gamma})^*}$ as follows
\begin{equation}\label{5.16'}
\Theta_{L_d}^-(g)(X_g, Y_g)= -X_g(L_d),\;\;\;\;\;
\Theta_{L_d}^+(g)(X_g, Y_g)= Y_g(L_d),
\end{equation}
 for each $g\in {\Gamma}$ and $(X_g, Y_g)\in {\mathcal T}_g^\Gamma\Gamma\equiv V_g\beta\oplus
V_g\alpha$.

Since $dL_d=\Theta_{L_d}^+ - \Theta_{L_d}^- $ and so, using
$d^2=0,$ it follows that $d\Theta_{L_d}^+=d\Theta_{L_d}^-$. This
means that there exists a unique 2-section
$\Omega_{L_d}=-d\Theta_{L_d}^+=-d\Theta_{L_d}^-$, which will be
called the \emph{Poincar\'e-Cartan} 2-section. This 2-section will
be important to study the symplectic character of the discrete
unconstrained Euler-Lagrange equations.

If $g$ is an element of $\Gamma$ such that $\alpha(g)=x$ and
$\beta(g)=y$ and $\{X_A\}$ (respectively, $\{Y_B\}$) is a local
basis of $\Sec{\tau}$ on the open subset $U$ (respectively,
$V$) of $M$, with $x\in U$ (respectively, $y\in V$), then on
$\alpha^{-1}(U)\cap \beta^{-1}(V)$ we have that
\begin{equation}\label{OmegaLd}
\begin{array}{rcl}\Theta_{L_d}^-&=&-\rvec{X_A}(L)(X^A)^{(1,0)},\;\;\;
\Theta_{L_d}^+=\lvec{Y_B}(L)(Y^B)^{(0,1)},\\
\Omega_{L_d}&=&-\rvec{X_A}(\lvec{Y_B}(L_d))(X^A)^{(1,0)}\wedge
(Y^B)^{(0,1)}
\end{array}
\end{equation}
where $\{X^A\}$ (respectively, $\{Y^B\}$) is the dual basis of
$\{X_A\}$ (respectively, $\{Y_B\}$) (for more details, see
\cite{groupoid}).

\subsubsection{ Discrete unconstrained
Lagrangian evolution operator} Let
$\Upsilon: {\Gamma}\to {\Gamma}$ be a smooth map such that:
\begin{enumerate}
\item[-] $\hbox{graph}(\Upsilon)\subseteq {\Gamma}_2$, that is, $(g, \Upsilon(g))\in
  {\Gamma}_2$, for all $g\in {\Gamma}$ ($\Upsilon$ is a \emph{second order
  operator}) and
\item[-] $(g, \Upsilon(g))$ is a solution of the discrete
unconstrained Euler-Lagrange
  equations, for all $g\in {\Gamma}$, that is,
  $(D_{{DEL}}L_d)(g,\Upsilon(g))=0,$ for all $g\in {\Gamma}.$
\end{enumerate}
In such a case
\begin{equation}\label{5.22}
  \lvec{X}(g)(L_d)-\rvec{X}(\Upsilon(g))(L_d)=0,
\end{equation}
for every section $X$ of $\tau:E_{\Gamma}\to M$ and every $g\in
{\Gamma}.$ The map $\Upsilon: {\Gamma}\to {\Gamma}$ is called a
\emph{discrete flow} or a \emph{discrete unconstrained Lagrangian
  evolution operator for $L_d$}.

Now, let $\Upsilon:{\Gamma}\to {\Gamma}$ be a second order
operator. Then, the prolongation $\prol[]{\Upsilon}:{\mathcal
T}^\Gamma\Gamma\equiv V\beta\oplus_{\Gamma} V\alpha\to {\mathcal
T}^\Gamma\Gamma \equiv V\beta\oplus_{\Gamma} V\alpha$ of
$\Upsilon$ is the Lie algebroid morphism over $\Upsilon:
{\Gamma}\to {\Gamma}$ defined as follows (see \cite{groupoid}):
\begin{equation}\label{poi}
\begin{array}{rcl}
  \prol[]{\Upsilon}[g](X_g, Y_g)&=&((T_g(r_{g\Upsilon(g)}\circ i))(Y_g),
  (T_g\Upsilon)(X_g)\\&&+(T_g\Upsilon)(Y_g)-T_g(r_{g\Upsilon(g)}\circ i)(Y_g)) ,
  \end{array}
\end{equation}
for all $(X_g, Y_g)\in {\mathcal T}_g^\Gamma\Gamma \equiv
V_g\beta\oplus V_g\alpha$. Moreover, from (\ref{linv}),
(\ref{rinv}) and (\ref{poi}), we obtain that
\begin{equation}\label{zaa}
  \prol[]{\Upsilon}[g](\rvec{X}(g), \lvec{Y}(g))=(-\rvec{Y}(\Upsilon(g)),
  (T_g\Upsilon)(\rvec{X}(g)+\lvec{Y}(g)) +\rvec{Y}(\Upsilon(g))) ,
\end{equation}
for all $X, Y\in \Sec{\tau}$.

Using (\ref{poi}), one may prove that (see \cite{groupoid}):
\begin{enumerate}
\item The map $\Upsilon$ is a discrete unconstrained Lagrangian
evolution operator for
  $L_d$ if and only if
  $(\prol[]{\Upsilon},\Upsilon)^*\Theta_{L_d}^-=\Theta_{L_d}^+$.
\item The map $\Upsilon$ is a discrete  unconstrained Lagrangian
evolution operator for
  $L_d$ if and only if
  $(\prol[]{\Upsilon},\Upsilon)^*\Theta_{L_d}^--\Theta_{L_d}^-=dL_d$.
\item If $\Upsilon$ is discrete unconstrained Lagrangian evolution
operator then
  $$(\prol[]{\Upsilon},\Upsilon)^*\Omega_{L_d}=\Omega_{L_d}.$$
\end{enumerate}

\subsubsection{ Discrete unconstrained Legendre transformations} Given a
Lagrangian $L_d: {\Gamma}\to \R$ we define the \emph{discrete
unconstrained Legendre transformations} $\F^{-}L_d: {\Gamma}\to
E^*_{\Gamma}$ and $\F^{+}L_d: {\Gamma}\to E^*_{\Gamma}$ by (see
\cite{groupoid})
\begin{align*}
  (\F^{-}L_d)(h)(v_{\epsilon(\alpha(h))})&= -v_{\epsilon(\alpha(h))}(L_d\circ
  r_h\circ i),\;\;\; \mbox{ for }
  v_{\epsilon(\alpha(h))}\in (E_{\Gamma})_{\alpha(h)},\\
  (\F^{+}L_d)(g)(v_{\epsilon(\beta(g))})&= v_{\epsilon(\beta(g))}(L_d\circ
  l_g), \mbox{ for } v_{\epsilon(\beta(g))}\in
  (E_{\Gamma})_{\beta(g)}.
\end{align*}
Now, we introduce the prolongations $\prol[\Gamma]{\F^{-}L_d}:
{\mathcal T}^\Gamma\Gamma\equiv V\beta\oplus_{\Gamma} V\alpha\to
\prol[E_{\Gamma}]{E^*_{\Gamma}}$ and $\prol[\Gamma]{\F^{+}L_d}:
{\mathcal T}^\Gamma\Gamma\equiv V\beta\oplus_{\Gamma} V\alpha\to
\prol[E_{\Gamma}]{E^*_{\Gamma}}$ by
\begin{eqnarray}
  \label{TF-L}\prol[\Gamma]{\F^{-}L_d}[h](X_h, Y_h) &= (T_h(i\circ
  r_{h^{-1}})(X_h),
  (T_h\F^{-}L_d)(X_h)+(T_h\F^{-}L_d)(Y_h)),\\\label{TF+L}
  \prol[\Gamma]{\F^{+}L_d}[g](X_g, Y_g) &= ((T_g l_{g^{-1}})(Y_g),
  (T_g\F^{+}L_d)(X_g)+(T_g\F^{+}L_d)(Y_g)),
\end{eqnarray}
for all $h,g\in {\Gamma}$ and $(X_h, Y_h)\in {\mathcal
T}_h^\Gamma\Gamma \equiv V_h\beta\oplus V_h\alpha$ and
$(X_g,Y_g)\in {\mathcal T}^\Gamma_g\Gamma\equiv V_g\beta \oplus
V_g\alpha$ (see \cite{groupoid}).  We observe that the discrete
Poincar\'e-Cartan 1-sections and 2-section are related to the
canonical Liouville section of
$(\prol[E_{\Gamma}]{E^*_{\Gamma}})^*\to E^*_{\Gamma}$ and the
canonical symplectic section of
$\wedge^2(\prol[E_{\Gamma}]{E^*_{\Gamma}})^*\to E^*_{\Gamma}$ by
pull-back under the discrete unconstrained Legendre
transformations, that is (see \cite{groupoid}),
\begin{eqnarray}
 \label{OmOml1} (\prol[\Gamma]{\F^{-}L_d}, \F^-L_d)^*\Theta= \Theta^-_{L_d},&&
  (\prol[\Gamma]{\F^{+}L_d}, \F^+L_d)^*\Theta= \Theta^+_{L_d},\\
 \label{OmOml} (\prol[\Gamma]{\F^{-}L_d}, \F^-L_d)^*\Omega= \Omega_{L_d}, &&
  (\prol[\Gamma]{\F^{+}L_d}, \F^+L_d)^*\Omega= \Omega_{L_d} .
\end{eqnarray}

\subsubsection{  Discrete regular Lagrangians}\label{subsection1.3.5} A discrete
Lagrangian $L_d:{\Gamma}\to \R$ is said to be \emph{regular} if
the Poincar{\'e}-Cartan $2$-section $\Omega_{L_d}$ is
nondegenerate on the Lie algebroid
$\widetilde{\tau}_\Gamma:{\mathcal T}^\Gamma\Gamma\equiv
V\beta\oplus_\Gamma V\alpha\to \Gamma$ (see \cite{groupoid}). In
\cite{groupoid}, we obtained some necessary and sufficient
conditions for a discrete Lagrangian on a Lie groupoid ${\Gamma}$
to be regular that we summarize as follows:
\begin{eqnarray*}
  L_d \hbox{ is regular} &\kern-2pt \Longleftrightarrow &\kern-2pt  \hbox{The Legendre
    transformation } \F^+L_d  \hbox{ is a local
    diffeomorphism}\\
  &\kern-2pt\Longleftrightarrow&\kern-2pt  \hbox{The Legendre transformation } \F^-L_d
  \hbox{ is a local
    diffeomorphism}
\end{eqnarray*}
Locally, we deduce that $L_d$ is regular if and only if for every
$g\in \Gamma$ and every local  basis $\{X_A\}$ (respectively,
$\{Y_B\}$) of $\Sec{\tau}$ on an open subset $U$
(respectively, $V$) of $M$ such that $\alpha(g)\in U$
(respectively, $\beta(g)\in V$) we have that the matrix
$(\rvec{X}_A(\lvec{Y}_B(L_d)))$ is regular on $\alpha^{-1}(U)\cap
\beta^{-1}(V).$

Now, let $L_d:\Gamma\to \R$ be a discrete Lagrangian and $g$ be a
point of $\Gamma$. We define the $\R$-bilinear map
$G_g^{L_d}:(E_\Gamma)_{\alpha(g)}\oplus (E_\Gamma)_{\beta(g)}\to
\R$ given by
\begin{equation}\label{1.21'}
G_g^{L_d}(a,b)=\Omega_{L_d}(g)((-T_{\epsilon(\alpha(g))}(r_g\circ
i)(a),0),(0,(T_{\epsilon(\beta(g))}l_g)(b))). \end{equation}

Then, using (\ref{OmegaLd}), we have that

\begin{proposition}\label{prop1.1}
The discrete Lagrangian $L_d:\Gamma\to \R$ is regular if and only
if $G_g^{L_d}$ is nondegenerate, for all $g\in \Gamma$, that is,
\[
G_g^{L_d}(a,b)=0,\mbox{ for all }b\in
(E_\Gamma)_{\beta(g)}\Rightarrow a=0
\]
(respectively, $G_g^{L_d}(a,b)=0,$ for all $a\in
(E_\Gamma)_{\alpha(g)}\Rightarrow b=0$).
\end{proposition}

 On the other hand, if $L_d:\Gamma\to \R$ is a
discrete Lagrangian on a Lie groupoid $\Gamma$ then we have that
\[
\tau^*\circ \F^-L_d=\alpha,\;\;\; \tau^*\circ \F^+L_d=\beta,
\]
where $\tau^*:E_\Gamma^*\to M$ is the vector bundle projection.
Using these facts, (\ref{TF-L}) and (\ref{TF+L}), we deduce the
following result.

\begin{proposition}\label{prop1.2}
Let $L_d:\Gamma\to \R$ be a discrete Lagrangian function. Then,
the following conditions are equivalent:
\begin{enumerate}
\item $L_d$ is regular. \item The linear map ${\mathcal
T}_{h}^\Gamma\F^-L_d:V_h\beta\oplus V_h\alpha\to {\mathcal
T}_{\F^-L_d(h)}^{E_\Gamma}E_\Gamma^*$ is a linear isomorphism, for
all $h\in \Gamma.$ \item The linear map ${\mathcal
T}_g^\Gamma\F^+L_d:V_g\beta \oplus V_g\alpha\to {\mathcal
T}_{\F^+L_d(g)}^{E_\Gamma}{E_\Gamma}^*$ is a linear isomorphism,
for all $g\in \Gamma$.
\end{enumerate}
\end{proposition}

Finally, let $L_d:\Gamma\to \R$ be a regular discrete Lagrangian
function and $(g_0,h_0)\in \Gamma\times \Gamma$ be a solution of
the discrete Euler-Lagrange equations for $L_d$. Then, one may
prove (see \cite{groupoid}) that there exist two open subsets
$U_{0}$ and $V_{0}$ of $\Gamma$, with $g_{0} \in U_{0}$ and $h_{0}
\in V_{0}$, and there exists a (local) discrete unconstrained
Lagrangian evolution operator $\Upsilon_{L_{d}}: U_{0} \to V_{0}$
such that:
\begin{enumerate}
\item
$\Upsilon_{L_{d}}(g_{0}) = h_{0}$,

\item
$\Upsilon_{L_{d}}$ is a diffeomorphism and

\item $\Upsilon_{L_{d}}$ is unique, that is, if $U'_{0}$ is an
open subset of $\Gamma$, with $g_{0} \in U_{0}',$ and
$\Upsilon'_{L_{d}}: U'_{0} \to \Gamma$ is a (local) discrete
Lagrangian evolution operator then
\[
\Upsilon_{L_d |U_{0}\cap U_{0}'} = \Upsilon'_{L_d |U_{0}\cap
U_{0}'}.
\]

\end{enumerate}

\section{Discrete Nonholonomic (or constrained) Lagrangian systems on Lie groupoids}

\subsection{Discrete Generalized H\"{o}lder's principle}
\label{DMLg}
 Let ${\Gamma}$ be  a Lie groupoid with structural
maps
\[
\alpha, \beta: {\Gamma} \to M, \; \; \epsilon: M \to {\Gamma}, \;
\; i: {\Gamma} \to {\Gamma}, \; \; m: {\Gamma}_{2} \to {\Gamma}.
\]
Denote by $\tau:E_{\Gamma}\to M$ the Lie algebroid associated to
${\Gamma}$. Suppose that the rank of $E_{\Gamma}$ is $n$ and that
the dimension of $M$ is $m$.

A generalized discrete nonholonomic (or constrained) Lagrangian
system on $\Gamma$ is determined by:
\begin{itemize}
\item[-]  a \emph{regular discrete Lagrangian} $L_d: {\Gamma}
\longrightarrow \R$, \item[-] a  \emph{constraint distribution},
${\mathcal D}_c$, which is a vector subbundle  of the bundle
$E_{\Gamma}\to M$ of admissible directions. We will denote by
$\tau_{{\mathcal D}_c}:{\mathcal D}_c\to M$ the vector bundle
projection and by $i_{{\mathcal D}_c}:{\mathcal D}_c\to E_\Gamma$
the canonical inclusion.

 \item[-] a \emph{discrete constraint embedded
submanifold} ${\mathcal M}_c$ of ${\Gamma}$, such that $dim
{\mathcal M}_c = dim {\mathcal D}_c = m + r$, with $r \leq n.$ We
will denote by $i_{{\mathcal M}_c}:{\mathcal M}_c\to \Gamma$ the
canonical inclusion.
\end{itemize}
\begin{remark}
Let $L_d:\Gamma \to \R$ be a regular discrete Lagrangian on a Lie
groupoid $\Gamma$ and ${\mathcal M}_c$ be a submanifold of $\Gamma$
such that $\epsilon ({\mathcal M})\subseteq {\mathcal M}_c$. Then,
$dim {\mathcal M}_c=m+r$, with $0\leq r\leq m$. Moreover, for every
$x\in M$, we may introduce the subspace ${\mathcal D}_c(x)$ of
$E_\Gamma (x)$ given by
\[
{\mathcal D}_c(x)=T_{\epsilon (x)}{\mathcal M}_c\cap E_\Gamma (x).
\]
Since the linear map $T_{\epsilon (x)}\alpha : T_{\epsilon
(x)}{\mathcal M}_c \to T_xM$ is an epimorphism, we deduce that
$dim {\mathcal D}_c(x)=r$. In fact, ${\mathcal D}_c= \bigcup
_{x\in M} {\mathcal D}_c(x)$ is a vector subbundle of $E_\Gamma$
(over $M$) of rank $r$. Thus, we may consider the discrete
nonholonomic system $(L_d,{\mathcal M}_c,{\mathcal D}_c)$ on the
Lie groupoid $\Gamma$.
\end{remark}

For $g\in {\Gamma}$ fixed, we consider the following set of
\emph{admissible sequences} of order $N$:
\[
{\mathcal C}^N_{g}\kern-1.1pt=\set{(g_1, \ldots, g_N)\in
{\Gamma}^N}{(g_k, g_{k+1})\in {\Gamma}_2, \hbox{ for } k=1,.., N-1
\hbox{ and } g_1 \ldots g_N=g }.
\]
Given a  tangent vector  at $(g_1, \ldots, g_N)$ to the manifold
${\mathcal C}^N_{g}$, we may write it as the tangent vector at
$t=0$ of a curve in ${\mathcal C}^N_{g}$, $t\in (-\varepsilon,
\varepsilon)\subseteq \R\longrightarrow c(t)$ which  passes
through $(g_1, \ldots, g_N)$ at $t=0$. This type of curves is of
the form
\[
c(t)=(g_1h_1(t), h_1^{-1}(t)g_2 h_2(t), \ldots,
h_{N-2}^{-1}(t)g_{N-1} h_{N-1}(t), h_{N-1}^{-1}(t)g_N)
\]
where $h_k(t)\in \alpha^{-1}(\beta(g_k)),$ for all $t,$ and
$h_k(0)=\epsilon(\beta(g_k))$ for $k=1,\ldots, N-1$.

Therefore, we may  identify the   tangent space to ${\mathcal
C}^N_g$ at $(g_1,\dots ,g_N)$ with
\[
T_{(g_1, g_2, .., g_N)}{\mathcal C}^N_g\equiv\set{(v_1, v_2,
\ldots, v_{N-1})}{v_k\in (E_\Gamma)_{x_k} \hbox{ and }
x_k=\beta(g_k), 1\leq k\leq N-1}.
\]
Observe that  each $v_k$ is the tangent vector to the curve $h_k$
at $t=0$.

The curve $c$ is called a \emph{variation} of $(g_1, \ldots, g_N)$
and $(v_1, v_2, \ldots, v_{N-1})$ is called an \emph{infinitesimal
variation} of $(g_1, \ldots, g_N)$.

Now, we define the  \emph{discrete action sum} associated to the
discrete Lagrangian $L_d: {\Gamma}\longrightarrow \R$ as
\[
\begin{array}{rcrcl}
{\mathcal S} L_d&:& {\mathcal C}^N_g& \longrightarrow &\R\\
          & &(g_1, \ldots, g_N)&\longmapsto&\displaystyle{ \sum_{k=1}^{N}
          L_d(g_k)}.
\end{array}
\]
We define the \emph{variation}  $\delta {\mathcal S} L_d: T_{(g_1,
\ldots, g_N)}{\mathcal C}^N_g\to \R$  as
\begin{eqnarray*}
\delta {\mathcal S} L_d(v_1, \ldots,
v_{N-1})&\kern-5pt=&\kern-5pt\frac{d}{dt}\Big|_{t=0}{\mathcal S} L_d(c(t))\\
&\kern-10pt =& \kern-10pt \frac{d}{dt}\Big|_{t=0} \left\{
L_d(g_1h_1(t))+L_d(h_1^{-1}(t)g_2 h_2(t))\right. \\
&\kern-10pt & \kern-10pt\left. + \ldots + L_d(h_{N-2}^{-1}(t)g_{N-1} h_{N-1}(t))+
L_d(h^{-1}_{N-1}(t)g_N) \right\}\\
&\kern-10pt =&\kern-10pt \sum_{k=1}^{N-1}\left(d^o (L_d\circ l_{g_k})(\epsilon
(x_k))(v_k)+d^o (L_d\circ r_{g_{k+1}}\circ
i)(\epsilon(x_{k}))(v_k)\right)
\end{eqnarray*}
where $d^o$ is the standard differential on ${\Gamma}$, i.e.,
$d^o$ is the differential of the Lie algebroid
$\tau_{\Gamma}:T{\Gamma}\to {\Gamma}.$ It is obvious from the last
expression that the definition of variation $\delta {\mathcal S}
L_d$ does not depend on the choice of variations $c$ of the
sequence $g$ whose infinitesimal variation is $(v_1, \ldots,
v_{N-1})$.

Next, we will introduce the subset $({\mathcal V}_c)_g$ of
$T_{(g_1, \ldots, g_N)}{\mathcal C}^N_g$ defined by
\[
({\mathcal V}_c)_g=\set{(v_1, \ldots, v_{N-1})\in T_{(g_1, \ldots,
g_N)}{\mathcal C}^N_g}{\forall k\in \{1,\ldots, N-1\},\quad v_k\in
{\mathcal D}_c}.
\]

Then, we will say that a sequence in ${\mathcal C}^N_g$ satisfying
the constraints determined by ${\mathcal M}_c$ is a
\emph{H\"older-critical point } of the discrete action sum
${\mathcal S} L_d$ if the restriction of $\delta {\mathcal S} L_d$
to $({\mathcal V}_c)_g$ vanishes, i.e.
\[
\delta {\mathcal S} L_d\Big|_{({\mathcal V}_c)_g}=0.
\]

\begin{definition}[Discrete H\"older's principle]
Given $g\in {\Gamma}$, a sequence $(g_1, \ldots, g_N)$ $\in {\mathcal
C}^N_g$ such that $g_k\in {\mathcal M}_c$, $1\leq k\leq N$, is a
solution of the discrete nonholonomic  Lagrangian system
determined by $(L_d, {\mathcal M}_c, {\mathcal D}_c)$ if and only
if  $(g_1, \ldots, g_N)$ is a H\"older-critical point of
${\mathcal S}L_d$.
\end{definition}

If $(g_1,\dots,g_N)\in {\mathcal C}_g^N\cap ({\mathcal M}_c\times
\dots \times {\mathcal M}_c)$  then $(g_1,\dots ,g_N)$ is a
solution of the nonholonomic discrete Lagrangian system if and
only if
\[
\sum_{k=1}^{N-1}(d^o(L_d\circ l_{g_k}) + d^o(L_d\circ
r_{g_{k+1}}\circ i))(\epsilon(x_k))_{|({\mathcal D}_c)_{x_k}}=0,
\]
where $\beta(g_k)=\alpha(g_{k+1})=x_k.$ For $N=2$, we obtain that
$(g,h)\in \Gamma_2\cap ({\mathcal M}_c\times {\mathcal M}_c)$
(with $\beta(g)=\alpha(h)=x$) is a solution if
\[
d^o(L_d\circ l_g + L_d\circ r_h\circ i)(\epsilon (x))_{|({\mathcal
D}_c)_x}=0.
\]

These equations will be called the \emph{ discrete nonholonomic
Euler-Lagrange equations for the system $(L_d,{\mathcal
M}_c,{\mathcal D}_c).$}

 Let $(g_1, \dots , g_N)$ be an element of ${\mathcal
C}^{N}_{g}$. Suppose that $\beta(g_{k}) = \alpha(g_{k+1}) =
x_{k}$, $1 \leq k \leq N -1$, and that $\{X_{Ak}\} = \{X_{ak},
X_{\alpha k}\}$ is a local adapted basis of $\Sec{\tau}$ on an
open subset $U_{k}$ of $M$, with $x_{k} \in U_{k}$. Here,
$\{X_{ak}\}_{1\leq a \leq r}$ is a local basis of
$\Sec{\tau_{{\mathcal D}_{c}}}$ and, thus, $\{X^{\alpha k}\}_{r+1
\leq \alpha \leq n}$ is a local basis of the space of sections of
the vector subbundle $\tau_{{\mathcal D}_{c}^0}: {\mathcal
D}_{c}^0 \to M$, where ${\mathcal D}_{c}^0$ is the annihilator of
${\mathcal D}_{c}$ and $\{X^{ak}, X^{\alpha k}\}$ is the dual
basis of $\{X_{ak}, X_{\alpha k}\}$. Then, the sequence $(g_{1},
\dots , g_{N})$ is a solution of the discrete nonholonomic
equations if $(g_1,\dots,g_N)\in {\mathcal M}_c\times\dots \times
{\mathcal M}_c$ and it satisfies the following closed system of
difference equations
\begin{eqnarray*}
0&=&\sum_{k=1}^{N-1}\left[\lvec{X_{ak}}\big({g_k})(L_d)-\rvec{X_{ak}}\big({g_{k+1}})(L_d)
\right],
\\
& = &\sum_{k=1}^{N-1}
\left[\pai{dL_d}{(X_{ak})^{(0,1)}}(g_k)-\pai{dL_d}{(X_{ak})^{(1,0)}}(g_{k+1})\right],
\end{eqnarray*}
for $1\leq a \leq r$, $d$ being the differential of the Lie
algebroid $\pi^\tau:{\mathcal T}^{\Gamma}{\Gamma} \equiv
V\beta\oplus_{\Gamma} V\alpha\longrightarrow {\Gamma}$. For $N=2$
we obtain that $(g, h)\in {\Gamma}_2\cap ({\mathcal M}_c\times
{\mathcal M}_c)$ (with $\beta(g)=\alpha(h)=x$) is a solution if
\[
\lvec{X_a}(g)(L_d)-\rvec{X_a}(h)(L_d)=0
\]
where $\{X_a\}$ is a local basis of $\Sec{\tau_{{\mathcal
D}_c}}$ on an open subset $U$ of $M$ such that $x \in U$.

Next, we describe an alternative version of these difference
equations. First observe that using the Lagrange multipliers the
discrete nonholonomic equations are rewritten as
\[
d^o\left[L_d\circ l_{g}+L_d\circ r_{h}\circ
i\right](\epsilon(x))(v)=\lambda_{\alpha} X^{\alpha}(x)(v),
\]
for $v\in (E_{\Gamma})_{x}$, with $(g, h) \in \Gamma_{2} \cap
({\mathcal M}_c \times {\mathcal M}_c)$ and $\beta(g) = \alpha(h)
= x$. Here, $\{X^\alpha\}$ is a local basis of sections of the
annihilator ${\mathcal D}^0_c.$

Thus, the discrete nonholonomic equations are:
\[
\lvec{Y}(g)(L_d)-\rvec{Y}(h)(L_d)=\lambda_{\alpha}
(X^{\alpha})(Y)|_{\beta(g)},\quad  (g, h) \in \Gamma_{2} \cap
({\mathcal M}_c \times {\mathcal M}_c),
\]
for all $Y\in \Sec{\tau }$ or, alternatively,
\[
\pai{dL_d-\lambda_{\alpha}
(X^{\alpha})^{(0,1)}}{Y^{(0,1)}}(g)-\pai{dL_d}{Y^{(1,0)}}(h)=0,
\quad (g, h) \in \Gamma_{2} \cap ({\mathcal M}_c \times {\mathcal
M}_c),
\]
for all $Y \in \Sec{\tau}$.

On the other hand, we may define the \emph{discrete nonholonomic
Euler-Lagrange operator} $D_{DEL}(L_d,{\mathcal M}_c,{\mathcal
D}_c):\Gamma_2\cap ({\mathcal M}_c\times {\mathcal M}_c)\to
{\mathcal D}_c^*$ as follows
\[
D_{DEL}(L_d, {\mathcal M}_{c}, {\mathcal D}_{c})(g, h) =
d^o\left[L_d\circ l_{g}+L_d\circ r_{h}\circ i\right](\epsilon
(x))_{| ({\mathcal D}_c)_x},
\]
for $(g, h) \in \Gamma_{2} \cap ({\mathcal M}_c \times {\mathcal
M}_c)$, with $\beta (g) = \alpha (h) = x \in M$.

Then, we may characterize the solutions of the discrete
nonholonomic equations as the sequences $(g_{1}, \dots , g_{N})$,
with $(g_{k}, g_{k+1}) \in \Gamma_{2} \cap ({\mathcal M}_c \times
{\mathcal M}_{c})$, for each $k \in \{1, \dots , N-1\}$, and
\[
D_{DEL}(L_{d}, {\mathcal M}_{c}, {\mathcal D}_{c})(g_{k}, g_{k+1})
= 0.
\]

\begin{remark}\label{dE-Lopo}{\rm
\begin{enumerate}
\item The set $\Gamma_{2} \cap ({\mathcal M}_{c} \times {\mathcal
M}_{c})$ is not, in general, a submanifold of ${\mathcal M}_{c}
\times {\mathcal M}_{c}$.

\item Suppose that $\alpha_{{\mathcal M}_{c}}: {\mathcal M}_{c}
\to M$ and $\beta_{{\mathcal M}_{c}}: {\mathcal M}_{c} \to M$ are
the restrictions to ${\mathcal M}_{c}$ of $\alpha: \Gamma \to M$
and $\beta: \Gamma \to M$, respectively. If $\alpha_{{\mathcal
M}_{c}}$ and $\beta_{{\mathcal M}_{c}}$ are submersions then
$\Gamma_{2} \cap ({\mathcal M}_{c} \times {\mathcal M}_{c})$ is a
submanifold of ${\mathcal M}_{c} \times {\mathcal M}_{c}$ of
dimension $m + 2r$.
\end{enumerate}
}
\end{remark}

\subsection{Discrete Nonholonomic Legendre transformations}

Let $(L_{d}, {\mathcal M}_{c}, {\mathcal D}_{c})$ be a discrete
nonholonomic Lagrangian system. We define the \emph{discrete
nonholonomic Legendre transformations}
\[
\F^-(L_{d}, {\mathcal M}_{c}, {\mathcal D}_{c}): {\mathcal M}_{c}
\to {\mathcal D}_{c}^{*} \; \; \mbox{ and } \; \; \F^+(L_{d},
{\mathcal M}_{c}, {\mathcal D}_{c}): {\mathcal M}_{c} \to
{\mathcal D}_{c}^{*}
\]
as follows:
\begin{eqnarray}\label{F-LMD}
\F^{-}(L_{d}, {\mathcal M}_{c}, {\mathcal
D}_{c})(h)(v_{\epsilon(\alpha(h))})=-v_{\epsilon(\alpha(h))}(L_{d}
\circ r_{h} \circ i), \; \; \mbox{ for } v_{\epsilon(\alpha(h))}
\in {\mathcal D}_{c}(\alpha(h)),&&
\\
\label{F+LMD} \F^{+}(L_{d}, {\mathcal M}_{c}, {\mathcal
D}_{c})(g)(v_{\epsilon(\beta(g))})=v_{\epsilon(\beta(g))}(L_{d}
\circ l_{g}), \; \; \mbox{ for } v_{\epsilon(\beta(g))} \in
{\mathcal D}_{c}(\beta(g)).&&
\end{eqnarray}
If $\F^{-}L_{d}: \Gamma \to E_{\Gamma}^{*}$ and $\F^{+}L_{d}:
\Gamma \to E_{\Gamma}^{*}$ are the standard Legendre
transformations associated with the Lagrangian function $L_{d}$
and $i_{{\mathcal D}_{c}}^{*}: E_{\Gamma}^{*} \to {\mathcal
D}_{c}^{*}$ is the dual map of the canonical inclusion
$i_{{\mathcal D}_{c}}: {\mathcal D}_{c} \to E_{\Gamma}$ then
\begin{equation}\label{Rel}
\F^{-}(L_{d}, {\mathcal M}_{c}, {\mathcal D}_{c}) = i_{{\mathcal
D}_{c}}^{*} \circ \F^{-}L_{d} \circ i_{{\mathcal M}_{c}}, \; \;
\F^{+}(L_{d}, {\mathcal M}_{c}, {\mathcal D}_{c}) = i_{{\mathcal
D}_{c}}^{*} \circ \F^{+}L_{d} \circ i_{{\mathcal M}_{c}}.
\end{equation}

\begin{remark}\label{1'}{\rm
\begin{enumerate}
\item
Note that
\begin{equation}\label{CDLT}
\tau_{{\mathcal D}_{c}}^{*} \circ \F^{-}(L_{d}, {\mathcal M}_{c},
{\mathcal D}_{c}) = \alpha_{{\mathcal M}_{c}}, \; \; \;
\tau_{{\mathcal D}_{c}}^{*} \circ \F^{+}(L_{d}, {\mathcal M}_{c},
{\mathcal D}_{c}) = \beta_{{\mathcal M}_{c}}.
\end{equation}

\item If $D_{DEL}(L_{d}, {\mathcal M}_{c}, {\mathcal D}_{c})$ is
the discrete nonholonomic Euler-Lagrange operator then
\[
D_{DEL}(L_{d}, {\mathcal M}_{c}, {\mathcal D}_{c})(g, h) =
\F^{+}(L_{d}, {\mathcal M}_{c}, {\mathcal D}_{c})(g) -
\F^{-}(L_{d}, {\mathcal M}_{c}, {\mathcal D}_{c})(h),
\]
for $(g, h) \in \Gamma_{2} \cap ({\mathcal M}_{c} \times {\mathcal
M}_{c})$.
\end{enumerate}
}
\end{remark}

On the other hand, since by assumption $L_{d}: \Gamma \to \R$ is a
regular discrete Lagrangian function, we have that the discrete
Poincar\'{e}-Cartan $2$-section $\Omega_{L_{d}}$ is symplectic on
the Lie algebroid $\tilde{\tau_{\Gamma}}: {\mathcal
T}^{\Gamma}{\Gamma} \to \Gamma$. Moreover, the regularity of $L$
is equivalent to the fact that the Legendre transformations
$\F^{-}L_{d}$ and $\F^{+}L_{d}$ to be local diffeomorphisms (see
Subsection \ref{subsection1.3.5}).

Next, we will obtain necessary and sufficient conditions for the
discrete nonholonomic Legendre transformations associated with the
system $(L_d,{\mathcal M}_c,{\mathcal D}_c)$ to be local
diffeomorphisms.

Let $F$ be the vector subbundle (over $\Gamma$) of
$\tilde{\tau_{\Gamma}}: {\mathcal T}^{\Gamma}\Gamma \to \Gamma$
whose fiber at the point $h \in \Gamma$ is
\[
F_{h} = \set{\gamma^{(1,0)}_{h}}{\gamma \in {\mathcal
D}_{c}(\alpha(h))^{0}}^{0} \subseteq {\mathcal
T}_{h}^{\Gamma}\Gamma.
\]
In other words,
\[
F_{h}^{0} = \set{\gamma^{(1,0)}_{h}}{\gamma \in {\mathcal
D}_{c}(\alpha(h))^{0}}.
\]
Note that the rank of $F$ is $n+r$.

We also consider the vector subbundle $\bar{F}$ (over $\Gamma$) of
$\tilde{\tau_{\Gamma}}: {\mathcal T}^{\Gamma}\Gamma \to \Gamma$ of
rank $n + r$ whose fiber at the point $g \in \Gamma$ is
\[
\bar{F}_{g} = \set{\gamma^{(0,1)}_{g}}{\gamma \in {\mathcal
D}_{c}(\beta(g))^{0}}^{0} \subseteq {\mathcal
T}_{g}^{\Gamma}\Gamma.
\]
\begin{lemma}\label{l2}
$F$ (respectively, $\bar{F}$) is a coisotropic vector subbundle of
the symplectic vector bundle $({\mathcal T}^{\Gamma}\Gamma,
\Omega_{L_{d}})$, that is,
\[
F_{h}^{\perp} \subseteq F_{h}, \; \; \mbox{ for every } h \in
\Gamma
\]
(respectively, $\bar{F}_{g}^{\perp} \subseteq \bar{F}_{g}$, for
every $g \in \Gamma$), where $F_{h}^{\perp}$ (respectively,
$\bar{F}_{g}^{\perp}$) is the symplectic orthogonal of $F_{h}$
(respectively, $\bar{F}_{g}$) in the symplectic vector space
$({\mathcal T}_{h}^{\Gamma}\Gamma, $ $\Omega_{L_{d}}(h))$
(respectively, $({\mathcal T}_{g}^{\Gamma}\Gamma,
\Omega_{L_{d}}(g))$).
\end{lemma}

\begin{proof}
If $h \in \Gamma$ we have that
\[
F_{h}^{\perp} = \flat_{\Omega_{L_{d}}(h)}^{-1}(F_{h}^{0}),
\]
$\flat_{\Omega_{L_{d}}(h)}: {\mathcal T}_{h}^{\Gamma}\Gamma \to
({\mathcal T}_{h}^{\Gamma}\Gamma)^{*}$ being the canonical
isomorphism induced by the symplectic form $\Omega_{L_{d}}(h)$.
Thus, using (\ref{OmegaLd}), we deduce that
\[
F_{h}^{\perp} = \set{\flat_{\Omega_{L_{d}}(h)}^{-1}(\gamma^{(1,
0)}_{h})}{\gamma \in {\mathcal D}_{c}(\alpha(h))^{0}} \subseteq
\{0\} \oplus V_{h}\alpha \subseteq F_{h}.
\]
The coisotropic character of $\bar{F}_{g}$ is proved in a similar
way.
\end{proof}

We also have the following result
\begin{lemma}\label{l3}
Let ${\mathcal T}^{\Gamma}\F^{-}L_{d}: {\mathcal T}^{\Gamma}\Gamma
\to {\mathcal T}^{E_{\Gamma}}E_{\Gamma}^{*}$ (respectively,
${\mathcal T}^{\Gamma}\F^{+}L_{d}: {\mathcal T}^{\Gamma}\Gamma \to
{\mathcal T}^{E_{\Gamma}}E_{\Gamma}^{*}$) be the prolongation of
the Legendre transformation $\F^{-}L_{d}: \Gamma \to
E_{\Gamma}^{*}$ (respectively, $\F^{+}L_{d}: \Gamma \to
E_{\Gamma}^{*}$). Then,
\[
({\mathcal T}^{\Gamma}_{h} \F^{-}L_{d})(F_{h}) = {\mathcal
T}^{{\mathcal D}_{c}}_{\F^{-}L_{d}(h)} E_{\Gamma}^{*} =\set{
(v_{\alpha(h)}, X_{\F^{-}L_{d}(h)}) \in {\mathcal
T}^{E_{\Gamma}}_{\F^{-}L_{d}(h)} E_{\Gamma}^{*}}{v_{\alpha(h)} \in
{\mathcal D}_{c}(\alpha(h))},
\]
for $h \in {\mathcal M}_{c}$ (respectively, \[({\mathcal
T}^{\Gamma}_{g} \F^{+}L_{d})(\bar{F}_{g}) = {\mathcal
T}^{{\mathcal D}_{c}}_{\F^{+}L_{d}(g)} E_{\Gamma}^{*} = \set{
(v_{\beta(g)}, X_{\F^{+}L_{d}(g)}) \in {\mathcal
T}^{E_{\Gamma}}_{\F^{+}L_{d}(g)} E_{\Gamma}^{*}}{v_{\beta(g)} \in
{\mathcal D}_{c}(\beta(g))},\] for $g \in {\mathcal M}_{c}$).
\end{lemma}

\begin{proof}
It follows using (\ref{(1,0)h}), (\ref{TF-L}) (respectively,
(\ref{(0,1)g}), (\ref{TF+L})) and Proposition \ref{prop1.2}.
\end{proof}

Now, we may prove the following theorem.

\begin{theorem}\label{t4}
Let $(L_{d}, {\mathcal M}_{c}, {\mathcal D}_{c})$ be a discrete
nonholonomic Lagrangian system. Then, the following conditions are
equivalent:
\begin{enumerate}
\item The discrete nonholonomic Legendre transformation
$\F^{-}(L_{d}, {\mathcal M}_{c}, {\mathcal D}_{c})$ (respectively,
$\F^{+}(L_{d}, {\mathcal M}_{c}, {\mathcal D}_{c})$) is a local
diffeomorphism.

\item
For every $h \in {\mathcal M}_{c}$ (respectively, $g \in {\mathcal
M}_{c}$)
\begin{equation}\label{Reg1}
(\rho^{{\mathcal
T}^{\Gamma}\Gamma})^{-1}(T_{h}{\mathcal M}_{c}) \cap F_{h}^{\perp}
= \{0 \}
\end{equation}
(respectively, $(\rho^{{\mathcal
T}^{\Gamma}\Gamma})^{-1}(T_{g}{\mathcal M}_{c}) \cap
\bar{F}_{g}^{\perp} = \{0 \}$).

\end{enumerate}
\end{theorem}

\begin{proof}

(i) $\Rightarrow$ (ii) If $h \in {\mathcal M}_{c}$ and $(X_{h},
Y_{h}) \in (\rho^{{\mathcal
T}^{\Gamma}\Gamma})^{-1}(T_{h}{\mathcal M}_{c}) \cap
F_{h}^{\perp}$ then, using the fact that $F_{h}^{\perp} \subseteq
\{0\} \oplus V_{h}\alpha$ (see the proof of Lemma \ref{l2}), we
have that $X_{h} = 0$. Therefore,
\begin{equation}\label{Yh}
Y_{h} \in V_{h}\alpha \cap T_{h}{\mathcal M}_{c}.
\end{equation}
Next, we will see that
\begin{equation}\label{TF-LY}
(T_{h}\F^{-}(L_{d}, {\mathcal M}_{c}, {\mathcal D}_{c}))(Y_{h}) =
0.
\end{equation}

From (\ref{CDLT}) and (\ref{Yh}), it follows that
$(T_{h}\F^{-}(L_{d}, {\mathcal M}_{c}, {\mathcal D}_{c}))(Y_{h})$
is vertical with respect to the projection $\tau_{{\mathcal
D}_{c}}^{*}: {\mathcal D}_{c}^{*} \to M$.

Thus, it is sufficient to prove that
\[
((T_{h}\F^{-}(L_{d}, {\mathcal M}_{c}, {\mathcal
D}_{c}))(Y_{h}))(\hat{Z}) = 0, \; \; \mbox{ for all } Z \in
\Sec{\tau_{{\mathcal D}_{c}}}.
\]
Here, $\hat{Z}: {\mathcal D}_{c}^{*} \to \R$ is the linear
function on ${\mathcal D}_{c}^{*}$ induced by the section $Z$.

Now, using (\ref{Rel}), we deduce that
\[
((T_{h}\F^{-}(L_{d}, {\mathcal M}_{c}, {\mathcal
D}_{c}))(Y_{h}))(\hat{Z}) = d (\hat{Z} \circ i_{{\mathcal
D}_{c}}^{*})((\F^-L_d)(h))(0, (T_{h}\F^{-}L_{d})(Y_h)),
\]
where $d$ is the differential of the Lie algebroid $\tau^{\tau^*}:
{\mathcal T}^{E_\Gamma}E_\Gamma^* \to E_\Gamma^*$.

Consequently, if $Z^{*{\bf c}}: E_{\Gamma}^{*} \to {\mathcal
T}^{E_{\Gamma}}E_{\Gamma}^{*}$ is the complete lift of $Z \in
\Sec{\tau}$, we have that (see (\ref{cl*})),
\begin{equation}\label{Zhat1}
\begin{array}{rcl}
((T_{h}\F^{-}(L_{d}, {\mathcal M}_{c}, {\mathcal
D}_{c}))(Y_{h}))(\hat{Z}) &=& \Omega(\F^{-}L_{d}(h))(Z^{*{\bf
c}}(\F^{-}L_{d}(h)),\\&&(0, (T_{h}\F^{-}L_{d})(Y_h)),
\end{array}
\end{equation}
$\Omega$ being the canonical symplectic section associated with
the Lie algebroid $E_{\Gamma}$.

On the other hand, since $Z \in \Sec{\tau_{{\mathcal
D}_{c}}}$, it follows that $Z^{*{\bf c}}(\F^{-}L_{d}(h))$ is in
${\mathcal T}^{{\mathcal D}_{c}}_{\F^{-}L_{d}(h)}E_{\Gamma}^{*}$
and, from Lemma \ref{l3}, we conclude that there exists $(X'_{h},
Y'_{h}) \in F_{h}$ such that
\begin{equation}\label{Zhat2}
({\mathcal T}_h^{\Gamma}\F^{-}L_{d})(X'_{h}, Y'_{h}) = Z^{*{\bf
c}}((\F^{-}L_{d})(h)).
\end{equation}

Moreover, using (\ref{TF-L}), we obtain that
\begin{equation}\label{Zhat3}
({\mathcal T}_h^{\Gamma}\F^{-}L_{d})(0, Y_{h}) = (0,
(T_{h}\F^{-}L_{d})(Y_h)).
\end{equation}

Thus, from (\ref{OmOml}), (\ref{Zhat1}), (\ref{Zhat2}) and
(\ref{Zhat3}), we deduce that
\[
((T_h\F^-(L_d,{\mathcal M},{\mathcal
D}_c))(Y_h))(\widehat{Z})=-\Omega_{L_d}(h)((0,Y_h),(X'_h,Y_h')).
\]
Therefore, since $(0,Y_h)\in F_h^\perp,$ it follows that
(\ref{TF-LY}) holds, which implies that $Y_h=0.$

This proves that $(\rho^{\mathcal
T^\Gamma\Gamma})^{-1}(T_h{\mathcal M}_c)\cap {F}_h^\perp =\{0\}.$

If $\F^+(L_d,{\mathcal M}_c,{\mathcal D}_c)$ is a local
diffeomorphism then, proceeding as above, we have that
$(\rho^{{\mathcal T}^\Gamma\Gamma})^{-1}(T_g{\mathcal M}_c)\cap
\bar{F}_g^\perp=\{0\},$ for all $g\in {\mathcal M}_c.$

(ii) $\Rightarrow$ (i) Suppose that $h\in {\mathcal M}_c$ and that
$Y_h$ is a tangent vector to ${\mathcal M}_c$ at $h$ such that
\begin{equation}\label{TF-LYO}
(T_h\F^{-}(L_d,{\mathcal M}_c,{\mathcal D}_c))(Y_h)=0.
\end{equation}
We have that $(T_h\alpha)(Y_h)=0$ and, thus,
\[
(0,Y_h)\in (\rho^{{\mathcal T^\Gamma}\Gamma})^{-1}(T_h{\mathcal
M}_c).
\]
We will see that $(0,Y_h)\in F_h^\perp,$ that is,

\begin{equation}\label{Fhort}
\Omega_{L_d}(h)((0,Y_h),(X'_h,Y'_h))=0,\;\;\; \mbox{ for }
(X_h',Y_h')\in F_h.
\end{equation}

Now, using (\ref{TF-L}) and (\ref{OmOml}), we deduce that
\[
\Omega_{L_d}(h)((0,Y_h),(X_h',Y_h'))=\Omega(\F^-
L_d(h))((0,(T_h\F^- L_d)(Y_h)),({\mathcal
T}^{\Gamma}_h\F^-L_d)(X_h',Y_h')).
\]
Therefore, from Lemma \ref{l3}, we obtain that
\[
\Omega_{L_d}(h)((0,Y_h),(X_h',Y_h'))=\Omega(\F^-L_d(h))(0,
(T_h\F^-L_d)(Y_h)), (v_{\alpha(h)},Y_{\F^-L_d(h)}))
\]
with $(v_{\alpha(h)}, Y_{\F^{-}L_{d}(h)}) \in {\mathcal
T}^{{\mathcal D}_{c}}_{\F^{-}L_{d}(h)} E_{\Gamma}^{*}$.

Next, we take a section $Z \in \Sec{\tau_{{\mathcal
D}_{c}}}$ such that $Z(\alpha(h)) = v_{\alpha(h)}$. Then (see
(\ref{cl*1})),
\[
(v_{\alpha(h)}, Y_{\F^{-}L_{d}(h)}) = Z^{*{\bf c}}(\F^-L_d(h)) +
(0, Y'_{\F^{-}L_{d}(h)}),
\]
where $Y'_{\F^{-}L_{d}(h)} \in T_{\F^{-}L_{d}(h)}E_{\Gamma}^{*}$
and $Y'_{\F^{-}L_{d}(h)}$ is vertical with respect to the
projection $\tau^*: E_{\Gamma}^{*} \to M$.

Thus, since (see Eq. (3.7) in \cite{LeMaMa})

\[
\Omega(\F^{-}L_{d}(h))((0, (T_{h}\F^{-}L_{d})(Y_{h})),
(0, Y'_{\F^{-}L_{d}(h)})) = 0,
\]
we have that
\[
\begin{array}{rcl}
\Omega_{L_{d}}(h) ((0, Y_{h}), (X'_{h}, Y'_{h}))& =&
-\Omega(\F^{-}L_{d}(h))(Z^{*{\bf c}}(\F^{-}L_{d}(h)), (0,
(T_{h}\F^{-}L_{d})(Y_{h}))) \\&=& -d (\hat{Z} \circ i_{{\mathcal
D}_{c}}^{*})(\F^{-}L_{d}(h))(0, (T_{h}\F^{-}L_{d})(Y_{h}))
\end{array}
\]
and, from (\ref{TF-LYO}), we deduce that (\ref{Fhort}) holds.

This proves that $Y_{h} \in  (\rho^{{\mathcal
T}^{\Gamma}\Gamma})^{-1}(T_{h}{\mathcal M}_{c}) \cap
F_{h}^{\perp}$ which implies that $Y_{h} = 0.$

Therefore, $\F^{-}(L_{d}, {\mathcal M}_{c}, {\mathcal D}_{c})$ is
a local diffeomorphism.

If  $(\rho^{{\mathcal T}^{\Gamma}\Gamma})^{-1}(T_{g}{\mathcal
M}_{c}) \cap \bar{F}_{g}^{\perp} = \{ 0 \}$ for all $g \in
{\mathcal M}_{c}$ then, proceeding as above, we obtain that
$\F^{+}(L_{d}, {\mathcal M}_{c}, {\mathcal D}_{c})$ is a local
diffeomorphism.
\end{proof}
Now, let $\rho^{{\mathcal T}^\Gamma\Gamma}:{\mathcal
T}^\Gamma\Gamma\to T\Gamma$ be the anchor map of the Lie algebroid
$\pi^\tau:{\mathcal T}^\Gamma\Gamma\to \Gamma.$ Then, we will
denote by ${\mathcal H}_h$ the subspace of ${\mathcal
T}_h^\Gamma\Gamma$ given by
\[{\mathcal H}_h=(\rho^{{\mathcal
T}^\Gamma\Gamma})^{-1}(T_h{\mathcal M}_c)\cap F_h,\;\;\; \mbox{
for } h\in {\mathcal M}_c.
\]
In a similar way, for every $g\in {\mathcal M}_c$ we will
introduce the subspace $\bar{\mathcal H}_g$ of ${\mathcal
T}^\Gamma_g\Gamma$ defined by
\[
\bar{\mathcal H}_g=(\rho^{{\mathcal
T}^\Gamma\Gamma})^{-1}(T_g{\mathcal M}_c)\cap \bar{F}_g. \]

On the other hand, let $h$ be a point of ${\mathcal M}_c$ and
$G_h^{L_d}:(E_\Gamma)_{\alpha(h)}\oplus (E_\Gamma)_{\beta(h)}\to
\R$ be the $\R$-bilinear map given by (\ref{1.21'}). We will
denote by $(\lvec{E}_\Gamma)^{{\mathcal M}_c}_{h}$ the subspace of
$(E_\Gamma)_{\beta(h)}$ defined by
\[
(\lvec{E}_\Gamma)_{h}^{{\mathcal M}_c}= \set{b\in
(E_\Gamma)_{\beta(h)}}{(T_{\epsilon(\beta(h))}l_h)(b)\in
T_h{\mathcal M}_c}
\]
and by $G_h^{L_dc}:({\mathcal D}_{c})_{\alpha(h)}\times
(\lvec{E}_\Gamma)_{h}^{{\mathcal M}_c}\to \R$ the restriction to
$({\mathcal D}_{c})_{\alpha(h)}\times (\lvec{E}_\Gamma)^{{\mathcal
M}_c}_{h}$ of the $\R$-bilinear map $G_h^{L_d}.$

In a similar way, if $g$ is a point of $\Gamma$ we will consider
the subspace $(\rvec{E}_\Gamma)^{{\mathcal M}_c}_{g}$ of
$(E_\Gamma)_{\alpha(g)}$ defined by
\[
(\rvec{E}_\Gamma)_{g}^{{\mathcal M}_c}=\set{a\in
(E_\Gamma)_{\alpha(g)}}{(T_{\epsilon(\alpha(g))}(r_g\circ
i))(a)\in T_g{\mathcal M}_c}
\] and the restriction
$\bar{G}_g^{L_dc}:(\rvec{E}_\Gamma)_{g}^{{\mathcal M}_c}\times
({\mathcal D}_c)_{\beta(g)}\to \R$ of $G_g^{L_d}$ to the space
$(\rvec{E}_\Gamma)_{g}^{{\mathcal M}_c}\times ({\mathcal
D}_c)_{\beta(g)}.$

\begin{proposition}\label{proposition2.7}
Let $(L_d,{\mathcal M}_c,{\mathcal D}_c)$ be a discrete
nonholonomic Lagrangian system. Then, the following conditions are
equivalent: \begin{enumerate} \item For every $h\in {\mathcal
M}_c$ (respectively, $g\in {\mathcal M}_c$)
\[
(\rho^{{\mathcal T}^\Gamma\Gamma})^{-1}(T_h{\mathcal M}_c)\cap
F_h^\perp=\{0\}
\]
(respectively, $(\rho^{{\mathcal
T}^\Gamma\Gamma})^{-1}(T_g{\mathcal M}_c)\cap
\bar{F}_g^\perp=\{0\}$). \item For every $h \in {\mathcal M}_{c}$
(respectively, $g \in {\mathcal M}_{c}$) the dimension of the
vector subspace ${\mathcal H}_{h}$ (respectively, $\bar{\mathcal
H}_{g}$) is $2r$ and the restriction to the vector subbundle
${\mathcal H}$ (respectively, $\bar{\mathcal H}$) of the
Poincar\'{e}-Cartan $2$-section $\Omega_{L_{d}}$ is nondegenerate.

\item For every $h\in {\mathcal M}_c$ (respectively, $g\in
{\mathcal M}_c$)
\[
\set{b\in (\lvec{E}_\Gamma)^{{\mathcal M}_c}_{h}}{
G_h^{L_dc}(a,b)=0, \forall a\in ({\mathcal
D}_c)_{\alpha(h)}}=\{0\}
\]
(respectively, $\set{a\in (\rvec{E}_\Gamma)^{{\mathcal
M}_c}_{g}}{G_g^{L_dc}(a,b)=0, \forall b\in ({\mathcal
D}_c)_{\beta(g)}}=\{0\}$).

\end{enumerate}
\end{proposition}
\begin{proof}
(i) $\Rightarrow$ (ii) Assume that $h\in {\mathcal M}_c$ and that
\begin{equation}\label{2.13}
(\rho^{{\mathcal T}^\Gamma\Gamma})^{-1}(T_h{\mathcal M}_c)\cap
F_h^\perp=\{0\}.
\end{equation}
 Let $U$ be an open
subset of $\Gamma$, with $h \in U$, and $\{\phi^{\gamma}\}_{\gamma
= 1, \dots , n-r}$ a set of independent real
$C^{\infty}$-functions on $U$ such that
\[
{\mathcal M}_{c} \cap U = \set{h' \in U}{\phi^{\gamma}(h') = 0,
\mbox{ for all } \gamma}.
\]
If $d$ is the differential of the Lie algebroid
$\widetilde{\tau}_\Gamma: {\mathcal T}^{\Gamma}\Gamma \to \Gamma$
then it is easy to prove that
\[
(\rho^{{\mathcal T}^{\Gamma}\Gamma})^{-1}(T_{h}{\mathcal M}_{c}) =
<\{d\phi^{\gamma}(h)\}>^{0}.
\]
Thus,
\begin{equation}\label{acotdim}
dim ((\rho^{{\mathcal T}^{\Gamma}\Gamma})^{-1}(T_{h}{\mathcal
M}_{c})) \geq n + r.
\end{equation}

On the other hand, $dim F_{h}^{\perp} = n-r$. Therefore, from
(\ref{2.13}) and (\ref{acotdim}), we obtain that
\[
dim ((\rho^{{\mathcal T}^{\Gamma}\Gamma})^{-1}(T_{h}{\mathcal
M}_{c})) = n + r
\]
and
\[
{\mathcal T}_{h}^{\Gamma}\Gamma = (\rho^{{\mathcal
T}^{\Gamma}\Gamma})^{-1}(T_{h}{\mathcal M}_{c}) \oplus
F_{h}^{\perp}.
\]
Consequently, using Lemma \ref{l2}, we deduce that
\begin{equation}\label{Descom}
 F_{h} = {\mathcal H}_{h} \oplus F_{h}^{\perp}.
\end{equation}

This implies that $dim {\mathcal H}_{h} = 2r$. Moreover, from
(\ref{Descom}), we also get that
\[
{\mathcal H}_{h} \cap {\mathcal H}_{h}^{\perp} \subseteq {\mathcal
H}_{h} \cap F_{h}^{\perp}
\]
and, since ${\mathcal H}_{h} \cap F_{h}^{\perp} = (\rho^{{\mathcal
T}^{\Gamma}\Gamma})^{-1}(T_{h}{\mathcal M}_{c}) \cap
F_{h}^{\perp}$ (see Lemma \ref{l2}), it follows that ${\mathcal
H}_{h} \cap {\mathcal H}_{h}^{\perp} = \{ 0 \}$.

Thus, we have proved that ${\mathcal H}_{h}$ is a symplectic
subspace of the symplectic vector space $({\mathcal
T}^{\Gamma}_{h}\Gamma, \Omega_{L_{d}}(h))$.

If $(\rho^{{\mathcal T}^{\Gamma}\Gamma})^{-1}(T_{g}{\mathcal
M}_{c}) \cap \bar{F}_{g}^{\perp} = \{ 0 \},$ for all $g \in
{\mathcal M}_{c}$ then, proceeding as above, we obtain that
$\bar{\mathcal H}_{g}$ is a symplectic subspace of the symplectic
vector space $({\mathcal T}_{g}^{\Gamma}\Gamma,
\Omega_{L_{d}}(g))$, for all $g \in {\mathcal M}_{c}$.

(ii) $\Rightarrow$ (i) Suppose that $h \in {\mathcal M}_{c}$ and
that ${\mathcal H}_{h}$ is a symplectic subspace of the symplectic
vector space $({\mathcal T}_{h}^{\Gamma}\Gamma,
\Omega_{L_{d}}(h))$.

If $(X_{h}, Y_{h}) \in (\rho^{{\mathcal
T}^{\Gamma}\Gamma})^{-1}(T_{h}{\mathcal M}_{c}) \cap
F_{h}^{\perp}$ then, using Lemma \ref{l2}, we deduce that $(X_{h},
Y_{h}) \in {\mathcal H}_{h}$.

Now, if $(X'_{h}, Y'_{h}) \in {\mathcal H}_{h}$ then, since
$(X_{h}, Y_{h}) \in F_{h}^{\perp}$, we conclude that
\[
\Omega_{L_{d}}(h)((X_{h}, Y_{h}), (X'_{h}, Y'_{h})) = 0.
\]
This implies that
\[
(X_{h}, Y_{h}) \in {\mathcal H}_{h} \cap {\mathcal H}_{h}^{\perp}
= \{0\}.
\]
Therefore, we have proved that $(\rho^{{\mathcal
T}^{\Gamma}\Gamma})^{-1}(T_{h}{\mathcal M}_{c}) \cap F_{h}^{\perp}
= \{0\}$.

If $\bar{\mathcal H}_{g} \cap \bar{\mathcal H}_{g}^{\perp} =
\{0\}$, for all $g \in {\mathcal M}_{c}$ then, proceeding as
above, we obtain that $(\rho^{{\mathcal
T}^{\Gamma}\Gamma})^{-1}(T_{g}{\mathcal M}_{c}) \cap
\bar{F}_{g}^{\perp} = \{0\}$, for all $g \in {\mathcal M}_{c}$.

(i) $\Rightarrow$ (iii) Assume that
\[
(\rho^{{\mathcal T}^\Gamma\Gamma})^{-1}(T_h{\mathcal M}_c)\cap
F_h^\perp =\{0\}
\]
and that $b\in (\lvec{E}_\Gamma)_{h}^{{\mathcal M}_c}$ satisfies
the following condition
\[
G_h^{L_d}(a,b)=0,\;\;\; \forall a\in ({\mathcal D}_c)_{\alpha(h)}.
\]
Then, $Y_h=(T_{\epsilon(\beta(h))}l_h)(b)\in T_h{\mathcal M}_c\cap
V_h\alpha$ and $(0,Y_h)\in (\rho^{{\mathcal
T}^\Gamma\Gamma})^{-1}(T_h{\mathcal M}_c)$.

Moreover, if $(X'_{h},Y'_h)\in F_h$, we have that
\[
X_h'=-(T_{\epsilon(\alpha(h))}(r_h\circ i))(a),\mbox{ with } a\in
({\mathcal D}_c)_{\alpha(h)}.
\]
Thus, using (\ref{OmegaLd}) and (\ref{1.21'}), we deduce that
\[
\Omega_{L_d}(h)((X_{h}',Y_h'),(0,Y_h))=\Omega_{L_d}(h)((X_h',0),(0,Y_h))=G_h^{L_d}(a,b)=0.\]
Therefore,
\[
(0,Y_h)\in (\rho^{{\mathcal T}^\Gamma\Gamma})^{-1}(T_h{\mathcal
M}_c)\cap F_h^{\perp}=\{0\},
\]
which implies that $b=0.$

If $(\rho^{{\mathcal T}^\Gamma\Gamma})^{-1}(T_g{\mathcal M}_c)\cap
\bar{F}_g^\perp=\{0\}$, for all $g\in {\mathcal M}_c$ then,
proceeding as above, we obtain that \[\set{a\in
(\rvec{E}_\Gamma)_g^{{\mathcal M}_c}}{G_g^{L_d c}(a,b)=0, \hbox{
for all } b\in ({\mathcal D}_c)_{\beta(g)}}=\{0\}.\]

(iii) $\Rightarrow$ (i) Suppose that $h\in {\mathcal M}_c$, that
\[ \set{b\in (\lvec{E}_\Gamma)^{{\mathcal M}_c}_{h}}{G_h^{L_d}(a,b)=0,\;\forall a\in ({\mathcal
D}_c)_{\alpha(h)}}=\{0\}\] and let $(X_h,Y_h)$ be an element of
the set $(\rho^{{\mathcal T}^\Gamma\Gamma})^{-1}(T_h{\mathcal
M}_c)\cap F_h^\perp.$

Then (see the proof of Lemma \ref{l2}), $X_h=0$ and $Y_h\in
T_h{\mathcal M}_c\cap V_h\alpha.$ Consequently,
\[
Y_h=(T_{\epsilon(\beta(h)}l_h)(b),\;\;\; \mbox{ with } b\in
(\lvec{E}_\Gamma)^{{\mathcal M}_c}_{h}.
\]

Now, if $a\in ({\mathcal D}_c)_{\alpha(h)},$ we have that
\[
X_h'=(T_{\epsilon(\alpha(h))}(r_h\circ i))(a)\in V_h\beta\mbox{
and } (X_h',0)\in F_h.\] Thus, from (\ref{1.21'}) and since
$(0,Y_h)\in F_h^{\perp},$ it follows that
$$G_h^{L_d}(a,b)=\Omega_{L_d}(h)((X_h',0)(0,Y_h))=0 .$$ Therefore,
$b=0$ which implies that $Y_h=0$.

If $\set{a\in (\rvec{E}_\Gamma)_g^{{\mathcal
M}_c}}{G_g^{L_dc}(a,b)=0,\forall b\in ({\mathcal
D}_c)_{\beta(g)}}=\{0\},$ for all $g\in {\mathcal M}_c$, then
proceeding as above we obtain that $(\rho^{{\mathcal
T}^\Gamma\Gamma})^{-1}(T_g{\mathcal M}_c)\cap
\bar{F}_g^\perp=\{0\},$ for all $g\in {\mathcal M}_c.$
\end{proof}

Using Theorem \ref{t4} and Proposition \ref{proposition2.7}, we
conclude
\begin{theorem}\label{t8}
Let $(L_d,{\mathcal M}_c,{\mathcal D}_c)$ be a discrete
nonholonomic  Lagrangian system. Then, the following conditions
are equivalent:
\begin{enumerate}
\item The discrete nonholonomic Legendre transformation
$\F^{-}(L_d,{\mathcal M}_c,{\mathcal D}_c)$ (respectively,
$\F^+(L_d,{\mathcal M}_c, {\mathcal D}_c)$) is a local
diffeomorphism. \item For every $h\in {\mathcal M}_c$
(respectively, $g\in {\mathcal M}_c)$
\[
(\rho^{{\mathcal T}^\Gamma\Gamma})^{-1}(T_h{\mathcal M}_c)\cap
F_h^\perp=\{0\}
\]
(respectively, $(\rho^{{\mathcal
T}^\Gamma\Gamma})^{-1}(T_g{\mathcal M}_c)\cap
\bar{F}_g^\perp=\{0\}).$

\item For every $h \in {\mathcal M}_{c}$ (respectively, $g \in
{\mathcal M}_{c}$) the dimension of the vector subspace ${\mathcal
H}_{h}$ (respectively, $\bar{\mathcal H}_{g}$) is $2r$ and the
restriction to the vector subbundle ${\mathcal H}$ (respectively,
$\bar{\mathcal H}$) of the Poincar\'{e}-Cartan $2$-section
$\Omega_{L_{d}}$ is nondegenerate.

\item For every $h\in {\mathcal M}_c$ (respectively, $g\in
{\mathcal M}_c$)
\[
\set{b\in (\lvec{E}_\Gamma)^{{\mathcal M}_c}_{h}}{G_h^{L_d
c}(a,b)=0, \forall a\in ({\mathcal D}_c)_{\alpha(h)}}=\{0\}
\]
(respectively, $\set{a\in (\rvec{E}_\Gamma)^{{\mathcal M}_c}_{g}}{
G_g^{L_dc}(a,b)=0, \forall b\in ({\mathcal
D}_c)_{\beta(g)}}=\{0\}$).
\end{enumerate}
\end{theorem}

\subsection{Nonholonomic evolution operators and regular discrete
nonholonomic Lagrangian systems}

First of all, we will introduce the definition of a nonholonomic
evolution operator.

\begin{definition} Let $(L_d,{\mathcal M}_c,{\mathcal D}_c)$ be a discrete
nonholonomic Lagrangian system and $\Upsilon_{nh}:{\mathcal
M}_c\to {\mathcal M}_c$ be a differentiable map. $\Upsilon_{nh}$
is said to be a discrete nonholonomic evolution operator for
$(L_d,{\mathcal M}_c,{\mathcal D}_c)$ if:
\begin{enumerate}
\item $\mbox{graph}(\Upsilon_{nh})\subseteq \Gamma_2$, that is,
$(g,\Upsilon_{nh}(g))\in \Gamma_2,$ for all $g\in {\mathcal M}_c$
and

\item $(g,\Upsilon_{nh}(g))$ is a solution of the discrete
nonholonomic equations, for all $g\in {\mathcal M}_c$, that is,
\[
d^o(L_d\circ l_g + L_d\circ r_{\Upsilon_{nh}(g)}\circ
i)(\epsilon(\beta(g)))_{|{\mathcal D}_c(\beta(g))}=0,\mbox{ for
all } g\in {\mathcal M}_c.
\]
\end{enumerate}
\end{definition}
\begin{remark}\label{r4.2-2}{\rm If $\Upsilon_{nh}:{\mathcal M}_c\to {\mathcal
M}_c$ is a differentiable map then, from (\ref{F-LMD}),
(\ref{F+LMD}) and (\ref{CDLT}), we deduce that $\Upsilon_{nh}$ is
a discrete nonholonomic evolution operator for $(L_d,{\mathcal
M}_c,{\mathcal D}_c)$ if and only if
\[
\F^-(L_d, {\mathcal M}_c,{\mathcal D}_c)\circ
\Upsilon_{nh}=\F^+(L_d,{\mathcal M}_c,{\mathcal D}_c).
\]}\end{remark}
Now, we will introduce the notion of a regular discrete
nonholonomic Lagrangian system.

\begin{definition}\label{d2.11}
A discrete nonholonomic Lagrangian system $(L_d,{\mathcal M}_c,
{\mathcal D}_c)$ is said to be regular if the discrete
nonholonomic Legendre transformations $\F^-(L_d,{\mathcal
M}_c,{\mathcal D}_c)$ and $\F^+(L_d,{\mathcal M}_c,{\mathcal
D}_c)$ are local diffeomorphims.
\end{definition}
From Theorem \ref{t8}, we deduce
\begin{corollary}\label{c3.12}
Let $(L_d,{\mathcal M}_c,{\mathcal D}_c)$ be a discrete
nonholonomic Lagrangian system. Then, the following conditions are
equivalent:
\begin{enumerate}
\item The system $(L_d,{\mathcal M}_c,{\mathcal D}_c)$ is regular.
\item The following relations hold
\[
(\rho^{{\mathcal T}^\Gamma\Gamma})^{-1}(T_h{\mathcal M}_c)\cap
F_h^\perp=\{0\},\mbox{ for all } h\in {\mathcal M}_c,
\]
\[
(\rho^{{\mathcal T}^\Gamma\Gamma})^{-1}(T_g{\mathcal M}_c)\cap
\bar{F}_g^\perp=\{0\},\mbox{ for all } g\in {\mathcal M}_c.
\]
\item ${\mathcal H}$ and $\bar{\mathcal H}$ are symplectic
subbundles of rank $2r$ of the symplectic vector bundle
$({\mathcal T}_{{\mathcal M}_c}^\Gamma\Gamma,\Omega_{L_d}).$

\item If $g$ and $h$ are points of ${\mathcal M}_c$ then the
$\R$-bilinear maps $G_{h}^{L_dc}$ and $\bar{G}_{g}^{L_dc}$ are
right and left nondegenerate, respectively.
\end{enumerate}
\end{corollary}
The map $G_{h}^{L_dc}$ (respectively, $\bar{G}_{g}^{L_dc}$) is right
nondegenerate (respectively, left nondegenerate) if
\[
G_h^{L_dc}(a,b)=0,\forall a\in ({\mathcal
D}_c)_{\alpha(h)}\Rightarrow b=0
\]
(respectively, $\bar{G}_g^{L_dc}(a,b)=0,\forall b\in ({\mathcal
D}_c)_{\beta(g)}\Rightarrow a=0$).

Every solution of the discrete nonholonomic equations for a
regular discrete nonholonomic Lagrangian system determines a
unique local discrete nonholonomic evolution operator. More
precisely, we may prove the following result:

\begin{theorem}
Let $(L_d,{\mathcal M}_c,{\mathcal D}_c)$ be a regular discrete
nonholonomic Lagrangian system and $(g_0,h_0)\in {\mathcal
M}_c\times {\mathcal M}_c$ be a solution of the discrete
nonholonomic equations for $(L_d,{\mathcal M}_c,{\mathcal D}_c).$
Then, there exist two open subsets $U_0$ and $V_0$ of $\Gamma$,
with $g_0\in U_0$ and $h_0\in V_0$, and there exists a local
discrete nonholonomic evolution operator
$\Upsilon_{nh}^{(L_d,{\mathcal M}_c,{\mathcal D}_c)}:U_0\cap
{\mathcal M}_c\to V_0\cap {\mathcal M}_c$ such that:
\begin{enumerate}
\item $\Upsilon_{nh}^{(L_d,{\mathcal M}_c,{\mathcal
D}_c)}(g_0)=h_0;$

\item $\Upsilon_{nh}^{(L_d,{\mathcal M}_c,{\mathcal D}_c)}$ is a
diffeomorphism and

\item $\Upsilon_{nh}^{(L_d,{\mathcal M}_c,{\mathcal D}_c)}$ is
unique, that is, if $U_0'$ is an open subset of $\Gamma$, with
$g_0\in U_0'$, and $\Upsilon_{nh}:U_0'\cap{\mathcal M}_c\to
{\mathcal M}_c$ is a (local) discrete nonholonomic evolution
operator then
\end{enumerate}

\vspace{-10pt}

\[
(\Upsilon_{nh}^{(L_d,{\mathcal M}_c,{\mathcal D}_c)})_{|U_0\cap
U_0'\cap {\mathcal M}_c}=(\Upsilon_{nh})_{|U_0\cap U_0'\cap
{\mathcal M}_c}.
\]
\end{theorem}
\begin{proof}
From remark \ref{1'}, we deduce that
\[
\F^+(L_d,{\mathcal M}_c,{\mathcal D}_c)(g_0)=\F^-(L_d,{\mathcal
M}_c,{\mathcal D}_c)(h_0)=\mu_0\in {\mathcal D}_c^*.
\]
Thus, we can choose two open subsets $U_0$ and $V_0$ of $\Gamma,$
with $g_0\in U_0$ and $h_0\in V_0,$ and an open subset $W_0$ of
$E_{\Gamma}^*$ such that $\mu_0\in W_0$ and
\[
\F^+(L_d,{\mathcal M}_c,{\mathcal D}_c):U_0\cap {\mathcal M}_c\to
W_0\cap {\mathcal D}_c^*,\;\;\; \F^-(L_d,{\mathcal M}_c,{\mathcal
D}_c):V_0\cap {\mathcal M}_c\to W_0\cap {\mathcal D}^*_c
\]
are diffeomorphisms. Therefore, from Remark \ref{r4.2-2}, we
deduce that
\[
\Upsilon_{nh}^{(L_d,{\mathcal M}_c,{\mathcal
D}_c)}=(\F^-(L_d,{\mathcal M}_c,{\mathcal D}_c)^{-1}\circ
\F^{+}(L_d,{\mathcal M}_c,{\mathcal D}_c))_{|U_0\cap {\mathcal
M}_c}:U_0\cap {\mathcal M}_c\to V_0\cap {\mathcal M}_c
\]
is a (local) discrete nonholonomic evolution operator. Moreover,
it is clear that $\Upsilon_{nh}^{(L_d,{\mathcal M}_c,{\mathcal
D}_c)}(g_0)=h_0$ and it follows that
$\Upsilon_{nh}^{(L_d,{\mathcal M}_c,{\mathcal D}_c)}$ is a
diffeomorphism.

Finally, if $U_0'$ is an open subset of $\Gamma$, with $g_0\in
U_0',$ and $\Upsilon_{nh}:U_0'\cap {\mathcal M}_c\to {\mathcal
M}_c$ is another (local) discrete nonholonomic evolution operator
then $(\Upsilon_{nh})_{|U_0\cap U_0'\cap {\mathcal M}_c}$ is also
a (local) discrete nonholonomic evolution operator. Consequently,
from Remark \ref{r4.2-2}, we conclude that
\[
\begin{array}{rcl}
(\Upsilon_{nh})_{|U_0\cap U_0'\cap {\mathcal
M}_c}&=&[\F^-(L_d,{\mathcal M}_c,{\mathcal D}_c)^{-1}\circ
\F^+(L_d,{\mathcal M}_c,{\mathcal D}_c)]_{|U_0\cap U_0'\cap
{\mathcal M}_c}\\&=&(\Upsilon_{nh}^{(L_d,{\mathcal M}_c,{\mathcal
D}_c)})_{|U_{0}\cap U_0'\cap {\mathcal M}_c}.
\end{array}
\]
\end{proof}
\subsection{Reversible discrete nonholonomic Lagrangian systems}
Let $(L_d,{\mathcal M}_c,$ ${\mathcal D}_c)$ be a discrete
nonholonomic Lagrangian system on a Lie groupoid
$\Gamma\rightrightarrows M$.

Following the terminology used in \cite{McLPerl} for the
particular case when $\Gamma$ is the pair groupoid $M\times M$, we
will introduce the following definition

\begin{definition}\label{d2.13-1}The discrete nonholonomic
Lagrangian system $(L_d,{\mathcal M}_c,{\mathcal D}_c)$ is said to
be reversible if
\[
L_d\circ i=L_d,\;\;\; i({\mathcal M}_c)={\mathcal M}_c,
\]
$i:\Gamma\to \Gamma$ being the inversion of the Lie groupoid
$\Gamma.$
\end{definition}
 For a reversible discrete nonholonomic Lagrangian system we have
the following result:
\begin{proposition}\label{p2.13-2} Let $(L_d,{\mathcal
M}_c,{\mathcal D}_c)$ be a reversible nonholonomic Lagrangian
system on a Lie groupoid $\Gamma$. Then, the following conditions
are equivalent:
\begin{enumerate}
\item The discrete nonholonomic Legendre transformation
$\F^-(L_d,{\mathcal M}_c,{\mathcal D}_c)$ is a local
diffeomorphism. \item The discrete nonholonomic Legendre
transformation $\F^+(L_d,{\mathcal M}_c,{\mathcal D}_c)$ is a
local diffeomorphism.
\end{enumerate}
\end{proposition}
\begin{proof}
If $h\in {\mathcal M}_c$ then, using (\ref{F-LMD}) and the fact
that $L_d\circ i=L_d,$ it follows that
\[
\F^-(L_d,{\mathcal M}_c,{\mathcal
D}_c)(h)(v_{\epsilon(\alpha(h))})=-v_{\epsilon(\alpha(h))}(L_d\circ
l_h^{-1})
\]
for $v_{\epsilon(\alpha(h))}\in ({\mathcal D}_c)_{\alpha(h)}.$
Thus, from (\ref{F+LMD}), we obtain that $$\F^-(L_d,{\mathcal
M}_c,{\mathcal
D}_c)(h)(v_{\epsilon(\alpha(h))})=-\F^+(L_d,{\mathcal
M}_c,{\mathcal D}_c)(h^{-1})(v_{\epsilon(\beta(h^{-1}))}).
$$

 This implies that
\[
\F^+(L_d,{\mathcal M}_c, {\mathcal D}_c)=-\F^-(L_d,{\mathcal
M}_c,{\mathcal D}_c)\circ i.
\]
Therefore, since the inversion is a diffeomorphism (in fact, we
have that $i^2={\rm id}$), we deduce the result
\end{proof}

Using Theorem \ref{t8}, Definition \ref{d2.11} and Proposition
\ref{p2.13-2}, we prove the following corollaries.

\begin{corollary}\label{c2.13-3}
Let $(L_d,{\mathcal M}_c,{\mathcal D}_c)$ be a reversible
nonholonomic Lagrangian system on a Lie groupoid $\Gamma$. Then,
the following conditions are equivalent:
\begin{enumerate}
\item The system $(L_d,{\mathcal M}_c,{\mathcal D}_c)$ is regular.
\item For all $h\in {\mathcal M}_c,$
\[(\rho^{{\mathcal T}^\Gamma\Gamma})^{-1}(T_h{\mathcal M}_c)\cap
F_h^\perp=\{0 \}.\] \item ${\mathcal H}=(\rho^{{\mathcal
T}^\Gamma\Gamma})^{-1}(T{\mathcal M}_c)\cap F$ is a symplectic
subbundle of the symplectic vector bundle $({\mathcal
T}^\Gamma_{{\mathcal M}_c}\Gamma,\Omega_{L_d}).$ \item The
$\mathbb{R}$-bilinear map
$G_h^{L_dc}:(\lvec{E}_\Gamma)_h^{{\mathcal M}_c}\times ({\mathcal
D}_c)_{\alpha(h)}\to \mathbb{R}$ is right nondegenerate, for all
$h\in {\mathcal M}_c.$ \end{enumerate}\end{corollary}

\begin{corollary}\label{2.13-4}
Let $(L_d,{\mathcal M}_c,{\mathcal D}_c)$ be a reversible
nonholonomic Lagrangian system on a Lie groupoid $\Gamma$. Then,
the following conditions are equivalent:
\begin{enumerate}
\item The system $(L_d,{\mathcal M}_c,{\mathcal D}_c)$ is regular.
\item For
all $g\in {\mathcal M}_c$,
\[
(\rho^{{\mathcal T}^\Gamma\Gamma})^{-1}(T_g{\mathcal M}_c)\cap
\bar{F}^\perp_g=\{0\}.\] \item $\bar{\mathcal H}=(\rho^{{\mathcal
T}^\Gamma\Gamma})^{-1}(T{\mathcal M}_{c})\cap \bar{F}$ is a
symplectic subbundle of the symplectic vector bundle $({\mathcal
T}_{{\mathcal M}_{c}}^\Gamma \Gamma,\Omega_{L_d}).$
 \item The $\mathbb{R}$-bilinear map
$\bar{G}_g^{L_dc}:({\mathcal D}_c)_{\beta(g)}\times
(\rvec{E}_\Gamma)_g^{{\mathcal M}_c}\to \mathbb{R}$ is left
nondegenerate, for all $g\in {\mathcal M}_c.$

\end{enumerate}
\end{corollary}

Next, we will prove that a reversible nonholonomic Lagrangian
system is dynamically reversible.

\begin{proposition}\label{2.13-5} Let $(L_d,{\mathcal
M}_c,{\mathcal D}_c)$ be a reversible nonholonomic Lagrangian
system on a Lie groupoid $\Gamma$ and $(g,h)$ be a solution of the
discrete nonholonomic Euler-Lagrange equations for $(L_d,{\mathcal
M}_c,{\mathcal D}_c).$ Then, $(h^{-1},g^{-1})$ is also a solution
of these equations. In particular, if the system $(L_d,{\mathcal
M}_c,{\mathcal D}_c)$ is regular and
$\Upsilon_{nh}^{(L_d,{\mathcal M}_c,{\mathcal D}_c)}$ is the
(local) discrete nonholonomic evolution operator for
$(L_d,{\mathcal M}_c,{\mathcal D}_c)$ then
$\Upsilon_{nh}^{(L_d,{\mathcal M}_c,{\mathcal D}_c)}$ is
reversible, that is,
\[
\Upsilon_{nh}^{(L_d,{\mathcal M}_c,{\mathcal D}_c)}\circ i\circ
\Upsilon_{nh}^{(L_d,{\mathcal M}_c,{\mathcal D}_c)}=i.
\]
\end{proposition}
\begin{proof}
Using that $i({\mathcal M}_c)={\mathcal M}_c,$ we deduce that
\[
(h^{-1},g^{-1})\in \Gamma_2\cap ({\mathcal M}_c\times {\mathcal
M}_c).
\]
Now, suppose that $\beta(g)=\alpha(h)=x$ and that $v\in ({\mathcal
D}_c)_x.$ Then, since $L_d\circ i=L_d,$ it follows that
\[
\begin{array}{rcl}
d^o[L_d\circ l_{h^{-1}} + L_d\circ r_{g^{-1}}\circ
i](\varepsilon(x))(v)&=&v(L_d\circ i\circ r_h\circ i) + v(L_d\circ
i\circ l_g)\\&=&v(L_d\circ l_g)+ v(L_d\circ r_h\circ i)=0.
\end{array}\]

Thus, we conclude that $(h^{-1},g^{-1})$ is a solution of the
discrete nonholonomic Euler-Lagrange equations for $(L_d,{\mathcal
M}_c,{\mathcal D}_c).$

If the system $(L_d,{\mathcal M}_c,{\mathcal D}_c)$ is regular and
$g\in {\mathcal M}_c,$ we have that
$(g,\Upsilon_{nh}^{(L_d,{\mathcal M},{\mathcal D}_c)}(g))$ is a
solution of the discrete nonholonomic Euler-Lagrange equations for
$(L_d,{\mathcal M}_c,{\mathcal D}_c).$ Therefore,
$(i(\Upsilon_{nh}^{(L_d,{\mathcal M},{\mathcal D}_c)}(g)),i(g))$
is also a solution of the dynamical equations which implies that
\[
\Upsilon_{nh}^{(L_d,{\mathcal M},{\mathcal
D}_c)}(i(\Upsilon_{nh}^{(L_d,{\mathcal M},{\mathcal
D}_c)}(g)))=i(g).
\]
\end{proof}

\begin{remark}{\rm Proposition \ref{2.13-5} was proved in
\cite{McLPerl} for the particular case when $\Gamma$ is the pair
groupoid. }
\end{remark}

\subsection{Lie groupoid morphisms and reduction}

Let  $(\Phi, \Phi_{0})$ be a  Lie groupoid morphism between the
Lie groupoids ${\Gamma} \rightrightarrows M$ and ${\Gamma}'
\rightrightarrows M'$.

Denote by $(E(\Phi),\Phi_0)$ the corresponding morphism between
the Lie algebroids $E_\Gamma$ and $E_{\Gamma'}$ of $\Gamma$ and
$\Gamma'$, respectively (see Section \ref{secLg}).

If $L_d:\Gamma\to \mathbb{R}$ and $L_d':\Gamma'\to \mathbb{R}$ are
discrete Lagrangians on $\Gamma$ and $\Gamma'$ such that
\[
L_d=L_d'\circ \Phi
\]
then, using Theorem 4.6 in \cite{groupoid}, we have that
\[
(D_{DEL}L_d)(g,h)(v)=(D_{DEL}L_d')(\Phi(g),\Phi(h))(E_x(\Phi)(v))
\]
for $(g, h)\in \Gamma_2$ and $v\in (E_\Gamma)_x,$ where
$x=\beta(g)=\alpha(h)\in M.$

Using this fact, we deduce the following result:

\begin{corollary}
Let $(\Phi,\Phi_0)$ be a Lie groupoid morphism between the Lie
groupoids $\Gamma\rightrightarrows M$ and
$\Gamma'\rightrightarrows M'$. Suppose that $L_d':\Gamma'\to
\mathbb{R}$ is a discrete Lagrangian on $\Gamma'$, that
$(L_d=L_d'\circ \Phi,{\mathcal M}_c,{\mathcal D}_c)$ is a discrete
nonholonomic Lagrangian system on $\Gamma$ and that $(g,h)\in
\Gamma_2\cap ({\mathcal M}_c\times{\mathcal M}_c)$. Then:

\begin{enumerate}
\item The pair $(g,h)$ is a solution of the discrete nonholonomic
problem $(L_d,{\mathcal M}_c, \linebreak {\mathcal D}_c)$ if and
only if $(D_{DEL}L_d')(\Phi(g),\Phi(h))$  vanishes over the set
\linebreak $(E_{\beta(g)}\Phi)(({\mathcal D}_c)_{\beta(g)}).$

\item If $(L_d',{\mathcal M}_c',{\mathcal D}'_c)$ is a discrete
nonholonomic Lagrangian system on $\Gamma'$ such that
$(\Phi(g),\Phi(h))\in {\mathcal M}_c'\times {\mathcal M}_c'$ and
$(E_{\beta(g)}(\Phi))(({\mathcal D}_{c})_{\beta(g)}) = ({\mathcal
D}_c')_{\Phi_0(\beta(g))}$ then $(g,h)$ is a solution for the
discrete nonholonomic problem $(L_d,{\mathcal M}_c,{\mathcal
D}_c)$ if and only if $(\Phi(g),\Phi(h))$ is a solution for the
discrete nonholonomic problem $(L_d',{\mathcal M}_c',{\mathcal
D}'_c)$.
\end{enumerate}
\end{corollary}

\subsection{ Discrete nonholonomic Hamiltonian evolution operator}

Let $(L_d, {\mathcal M}_c,$ ${\mathcal D}_c)$ a regular discrete
nonholonomic system. Assume,  without the loss of generality, that
the discrete nonholonomic Legendre transformations $\F^-(L_d,
{\mathcal M}_c, {\mathcal D}_c): {\mathcal M}_c\longrightarrow
{\mathcal D}_c^*$ and $\F^+(L_d, {\mathcal M}_c, {\mathcal D}_c):
{\mathcal M}_c\longrightarrow {\mathcal D}_c^*$ are global
diffeomorphisms. Then, $\gamma_{nh}^{(L_d, {\mathcal
M}_c, {\mathcal D}_c)}=\F^-(L_d, {\mathcal
M}_c, {\mathcal D}_c)^{-1}\circ\F^+(L_d, {\mathcal
M}_c, {\mathcal D}_c)$ is the discrete nonholonomic
 evolution operator and one may define the
\emph{discrete nonholonomic Hamiltonian evolution operator},
$\tilde{\gamma}_{nh}: {\mathcal D}_c^*\to {\mathcal D}_c^*$, by
\begin{equation}\label{dheo}
\tilde{\gamma}_{nh}=\F^+(L_d, {\mathcal M}_c, {\mathcal D}_c)\circ
\gamma_{nh}^{(L_d, {\mathcal
M}_c, {\mathcal D}_c)}\circ \F^+(L_d, {\mathcal M}_c, {\mathcal
D}_c)^{-1}\; .
\end{equation}

From Remark \ref{r4.2-2}, we have  the following alternative
definitions
\begin{eqnarray*}
 \tilde{\gamma}_{nh}&=&\F^-(L_d,
{\mathcal M}_c, {\mathcal D}_c)\circ \gamma_{nh}^{(L_d, {\mathcal
M}_c, {\mathcal D}_c)}\circ
\F^-(L_d, {\mathcal M}_c, {\mathcal D}_c)^{-1},\\
\tilde{\gamma}_{nh}&=&\F^+(L_d, {\mathcal M}_c, {\mathcal
D}_c)\circ \F^-(L_d, {\mathcal M}_c, {\mathcal
D}_c)^{-1} \end{eqnarray*} of the discrete Hamiltonian
evolution operator. The following commutative diagram illustrates
the situation

\unitlength=0.9mm \linethickness{0.4pt}
\begin{picture}(120.00,46.00)(10,25)
\put(55.00,60.00){\makebox(0,0)[cc]{${\mathcal M}_c$}}
\put(95.00,60.00){\makebox(0,0)[cc]{${\mathcal M}_c$}}
\put(58.67,60.00){\vector(1,0){30.33}}
\put(71.67,61.00){\makebox(0,0)[cb]{$\gamma_{nh}^{(L_d, {\mathcal
M}_c, {\mathcal D}_c)}$}}
\put(31.00,30.00){\makebox(0,0)[cc]{${\mathcal D}^*_c$}}
\put(75,30.00){\makebox(0,0)[cc]{${\mathcal D}^*_c$}}
\put(119.00,30.00){\makebox(0,0)[cc]{${\mathcal D}^*_c$}}
\put(52.00,56.67){\vector(-2,-3){16.00}}
\put(56.67,56.67){\vector(2,-3){16.00}}
\put(93.00,56.67){\vector(-2,-3){16.00}}
\put(98.33,56.67){\vector(2,-3){16.00}}
\put(37.00,29.67){\vector(1,0){31.33}}
\put(80.33,29.67){\vector(1,0){31.33}}
\put(40.33,40.00){\makebox(0,0)[rb]{\footnotesize $\F^-(L_d,
{\mathcal M}_c, {\mathcal D}_c)$}}
\put(63.33,47.00){\makebox(0,0)[lb]{\footnotesize$\F^+(L_d,
{\mathcal M}_c, {\mathcal D}_c)$}}
\put(84.00,40.00){\makebox(0,0)[lb]{\footnotesize$\F^-(L_d,
{\mathcal M}_c, {\mathcal D}_c)$}}
\put(106.33,47.00){\makebox(0,0)[lb]{\footnotesize$\F^+(L_d,
{\mathcal M}_c, {\mathcal D}_c)$}}
\put(50.00,33.13){\makebox(0,0)[cb]{$\tilde{\gamma}_{nh}$}}
\put(98.00,33.13){\makebox(0,0)[cb]{$\tilde{\gamma}_{nh}$}}
\end{picture}

\begin{remark}
{\rm The discrete nonholonomic evolution operator is an
application from ${\mathcal D}^*_c$ to itself. It is remarkable
that ${\mathcal D}_c^*$  is also the appropriate nonholonomic
momentum space for a continuous nonholonomic system defined by a
Lagrangian $L:E_{\Gamma}\to \R$ and the constraint distribution
${\mathcal D}_c$. Therefore, in the regular case, the solution of
the continuous nonholonomic Lagrangian system also determines a
flow from ${\mathcal D}^*_c$ to itself. We consider that this
would be a good starting point to compare the discrete and
continuous dynamics and eventually to establish a backward error
analysis for nonholonomic systems. } \end{remark}

\subsection{The discrete nonholonomic momentum map}

Let $(L_d,{\mathcal M}_c,{\mathcal D}_c)$ be a regular discrete
nonholonomic Lagrangian system on a Lie groupoid
$\Gamma\rightrightarrows M$ and $\tau:E_\Gamma\to M$ be the Lie
algebroid of $\Gamma.$

Suppose that ${\frak g}$ is a Lie algebra and that $\Psi:{\frak
g}\to \Sec{\tau}$ is a $\mathbb{R}$-linear map. Then, for each
$x\in M,$ we consider the vector subspace ${\frak g}^x$ of ${\frak
g}$ given by
\[
{\frak g}^x=\set{\xi\in {\frak g}}{\Psi(\xi)(x)\in ({\mathcal
D}_c)_x}
\]
and the disjoint union of these vector spaces
\[
{\frak g}^{{\mathcal D}_c}=\bigcup_{x\in M}{\frak g}^x.\] We will
denote by $({\frak g}^{{\mathcal D}_c})^*$  the disjoint union of
the dual spaces, that is,
\[
({\frak g}^{{\mathcal D}_c})^*=\bigcup_{x\in M}({\frak g}^x)^*.\]

Next, we define the \emph{discrete nonholonomic momentum  map}
$J^{nh}:\Gamma\to ({\frak g}^{{\mathcal D}_c})^*$ as follows: $
J^{nh}(g)\in ({\frak g}^{\beta(g)})^*$ and
\[
J^{nh}(g)(\xi)=\Theta_{L_d}^+(\Psi(\xi)^{(1,1)})(g)=\lvec{\Psi(\xi)}(g)(L_d),
\mbox{ for } g\in \Gamma \mbox{ and } \xi\in {\frak g}^{\beta(g)}.
\]

If $\widetilde{\xi}:M\to {\frak g}$ is a smooth map such that
$\widetilde{\xi}(x)\in {\frak g}^x,$ for all $x\in M,$ then
we may consider the smooth function $J_{\widetilde{\xi}}^{nh}:\Gamma\to
\mathbb{R}$ defined by
\[
J^{nh}_{\widetilde{\xi}}(g)=J^{nh}(g)(\widetilde{\xi}(\beta(g))),\;\;\forall
g\in \Gamma.
\]
\begin{definition} The Lagrangian $L_d$ is said to be ${\frak
g}$-invariant with respect $\Psi$ if
\[
\Psi(\xi)^{(1,1)}(L_d)=\lvec{\Psi(\xi)}(L_d)-\rvec{\Psi(\xi)}(L_d)=0,\;\;\;
\forall \xi\in {\frak g}.
\]
\end{definition}
Now, we will prove the following result
\begin{theorem}\label{2.14-2}
Let $\Upsilon_{nh}^{(L_d,{\mathcal M}_c,{\mathcal D}_c)}:{\mathcal
M}_c\to {\mathcal M}_c$ be the local discrete nonholonomic
evolution operator for the regular system $(L_d,{\mathcal M}_c,{\mathcal
D}_c).$ If $L_d$ is ${\frak g}$-invariant with respect to
$\Psi:{\frak g}\to \Sec{\tau}$ and $\widetilde{\xi}:M\to {\frak g}$
is a smooth map such that $\widetilde{\xi}(x)\in {\frak g}^x$,
for all $x\in M,$ then
\[
\begin{array}{rcl}
J_{\widetilde{\xi}}^{nh}(\Upsilon_{nh}^{(L_d,{\mathcal
M}_c,{\mathcal D}_c)}(g))-J_{\widetilde{\xi}}^{nh}(g)&=&\\[8pt] &\kern-85pt=&\kern-50pt
\lvec{\Psi(\widetilde{\xi}(\beta(\Upsilon_{nh}^{(L_d,{\mathcal
M}_c,{\mathcal D}_c)}(g)))
-\widetilde{\xi}(\beta(g)))}(\Upsilon_{nh}^{(L_d,{\mathcal
M}_c,{\mathcal D}_c)}(g))(L_d)\end{array}
\]
for $g\in {\mathcal M}_c.$
\end{theorem}
\begin{proof}Using that the Lagrangian $L_d$ is ${\frak g}$-invariant with respect to $\Psi$,
we have that
 \begin{equation}\label{2.17-1}
 \begin{array}{rcl}
 \rvec{\Psi(\tilde{\xi}(\alpha (\Upsilon_{nh}^{(L_d,{\mathcal M}_c,{\mathcal D}_c)}(g))))}(\Upsilon_{nh}^{(L_d,{\mathcal
 M}_c,{\mathcal  D}_c)}(g))(L_d)&=&\\[8pt]&\kern-250pt=&\kern-120pt\lvec{\Psi(\tilde{\xi}(\alpha
(\Upsilon_{nh}^{(L_d,{\mathcal M}_c,{\mathcal
D}_c)}(g))))}(\Upsilon_{nh}^{(L_d,{\mathcal M}_c,{\mathcal
D}_c)}(g))(L_d).\end{array}
\end{equation}
Also, since $(g,\Upsilon_{nh}^{(L_d,{\mathcal M}_c,{\mathcal
D}_c)}(g))$ is a solution of the discrete nonholonomic equations:
\begin{eqnarray}\label{2.17-2}
\lvec{\Psi(\tilde{\xi}(\beta(g)))}(g)(L_d)&=&
\rvec{\Psi(\tilde{\xi}(\alpha(\Upsilon_{nh}^{(L_d,{\mathcal
M}_c,{\mathcal D}_c)}(g))))}(\Upsilon_{nh}^{(L_d,{\mathcal
M}_c,{\mathcal D}_c)}(g))(L_d).
\end{eqnarray}
Thus, from (\ref{2.17-1}) and  (\ref{2.17-2}), we find that
\[
\lvec{\Psi(\tilde{\xi}(\beta(g))}(g)(L_d)=\lvec{\Psi(\tilde{\xi}(\beta(g)))
}(\Upsilon_{nh}^{(L_d,{\mathcal M}_c,{\mathcal D}_c)}(g))(L_d).
\]
Therefore, \begin{eqnarray*}
J^{nh}_{\tilde{\xi}}(\Upsilon_{nh}^{(L_d,{\mathcal M}_c,{\mathcal
D}_c)}(g))- J^{nh}_{\tilde{\xi}}(g)&=&
\lvec{\Psi\left(\tilde{\xi}(\beta(\Upsilon_{nh}^{(L_d,{\mathcal
M}_c,{\mathcal D}_c)}(g)))\right)}(\Upsilon_{nh}^{(L_d,{\mathcal
M}_c,{\mathcal
D}_c)}(g))(L_d)\\&&-\lvec{\Psi\left(\tilde{\xi}(\beta(g))\right)}(g)(L_d)\\
&=&\lvec{\Psi\left(\tilde{\xi}(\beta(\Upsilon_{nh}^{(L_d,{\mathcal
M}_c,{\mathcal D}_c)}(g)))\right)}(\Upsilon_{nh}^{(L_d,{\mathcal
M}_c,{\mathcal
D}_c)}(g))(L_d)\\&&-\lvec{\Psi(\tilde{\xi}(\beta(g))
}(\Upsilon_{nh}^{(L_d,{\mathcal M}_c,{\mathcal D}_c)}(g))(L_d)\\
&\kern-90pt=&\kern-50pt\lvec{\Psi\left(\tilde{\xi}(\beta(\Upsilon_{nh}^{(L_d,{\mathcal
M}_c,{\mathcal D}_c)}(g)))-\tilde{\xi}(\beta(g))\right)
}(\Upsilon_{nh}^{(L_d,{\mathcal M}_c,{\mathcal D}_c)} (g))(L_d).
\end{eqnarray*}
\end{proof}

Theorem \ref{2.14-2} suggests us to introduce the following
definition

\begin{definition}\label{2.14-3} An element $\xi\in {\frak g}$ is
said to be a \emph{ horizontal symmetry} for the discrete
nonholonomic system $(L_d,{\mathcal M}_c,{\mathcal D}_c)$ and the
map $\Psi:{\frak g}\to \Sec{\tau}$ if
\[
\Psi(\xi)(x)\in ({\frak D}_c)_x,\;\; \mbox{ for all } x\in M.
\]
\end{definition}
Now, from Theorem \ref{2.14-2}, we conclude that

\begin{corollary}
If $L_d$ is ${\frak g}$-invariant with respect to $\Psi$ and
$\xi\in {\frak g}$ is a horizontal  symmetry for $(L_d,{\mathcal
M}_c,{\mathcal D}_c)$ and $\Psi: {\frak g} \to \Sec{\tau}$ then
$J^{nh}_{\tilde{\xi}}:\Gamma\to \mathbb{R}$ is a constant of the
motion for $\Upsilon_{nh}^{(L_d,{\mathcal M}_c,{\mathcal D}_c)}$,
that is,
\[
J_{\tilde\xi}^{nh}\circ \Upsilon_{nh}^{(L_d,{\mathcal M}_c,{\mathcal
D}_c)}=J_{\tilde\xi}^{nh}.\]
\end{corollary}

\section{Examples}
\subsection{Discrete holonomic Lagrangian systems on a Lie
groupoid}

Let us examine the case when the system is subjected to holonomic
constraints.

Let $L_d:\Gamma \to \R$ be a discrete Lagrangian on a Lie groupoid
$\Gamma \rightrightarrows M$. Suppose that ${\mathcal M}_c\subseteq \Gamma$
is a Lie subgroupoid of $\Gamma$ over $M'\subseteq M$, that is,
${\mathcal M}_c$ is a Lie groupoid over $M'$ with structural maps
\[
\alpha_{|{\mathcal M}_c}:{\mathcal M}_c\to M',\;
\beta_{|{\mathcal M}_c}:{\mathcal M}_c\to M',\;
\epsilon _{|M'}:M'\to {\mathcal M}_c,\;
i_{|{\mathcal M}_c}:{\mathcal M}_c\to {\mathcal M}_c,
\]
the canonical inclusions $i_{{\mathcal M}_c}: {\mathcal
M}_c\longrightarrow \Gamma$ and $i_{M'}: M'\longrightarrow M$
are injective immersions and the pair $(i_{{\mathcal M}_c},i_{M'})$
is a Lie groupoid morphism. We may assume, without the loss of
generality, that $M'=M$ (in other case, we will replace the
Lie groupoid $\Gamma$ by the Lie subgroupoid $\Gamma '$ over
$M'$ defined by $\Gamma '=\alpha ^{-1}(M')\cap \beta ^{-1}(M')$).

Then, if $L_{{\mathcal M}_c}=L_d\circ i_{{\mathcal M}_c}$ and
$\tau _{{\mathcal M}_c}:E_{{\mathcal M}_c}\to M$ is the Lie algebroid of ${\mathcal M}_c$,
we have that the discrete (unconstrained) Euler-Lagrange equations for the
Lagrangian function $L_{{\mathcal M}_c}$ are:
\begin{equation}\label{saq}
\lvec{X}(g)(L_{{\mathcal M}_c})-\rvec{X}(h)(L_{{\mathcal
M}_c})=0,\quad (g, h) \in ({\mathcal M}_c)_2,
\end{equation}
for $X\in \Sec{\tau_{{\mathcal M}_c}}$.

We are interested in writing these equations in terms of the
Lagrangian $L_d$ defined on the Lie groupoid $\Gamma$. From
Corollary 4.7 (iii) in \cite{groupoid}, we deduce that $(g, h)\in
({\mathcal M}_c)_2$ is a solution of Equations \ref{saq} if and
only if  $D_{DEL} L_d (g, h)$ vanishes over ${\rm Im}
(E_{\beta(g)}(i_{{\mathcal M}_c}))$. Here, $E(i_{{\mathcal M}_c}):
E_{{\mathcal M}_c}\to E_\Gamma$ is the Lie algebroid morphism
induced between $E_{{\mathcal M}_c}$ and  $E_\Gamma$ by the Lie
groupoid morphism $(i_{{\mathcal M}_c}, Id)$. Therefore, we may
consider the discrete holonomic system as the discrete
nonholonomic system $(L_d, {\mathcal M}_c, {\mathcal D}_c)$, where
${\mathcal D}_c= (E (i_{{\mathcal M}_c}))( E_{{\mathcal
M}_c})\cong E_{{\mathcal M}_c}$.

In the particular case, when the  subgroupoid ${\mathcal M}_c$ is
determined by the vanishing set of $n-r$ independent real $C^{\infty}$-functions
$\phi ^\gamma:\Gamma \to \R$:
\[
{\mathcal M}_{c}  = \set{g \in \Gamma}{\phi^{\gamma}(g) = 0,
\mbox{ for all } \gamma },
\]
then the discrete holonomic equations are equivalent to:
\begin{eqnarray*}
\lvec{Y}(g)(L_d)-\rvec{Y}(h)(L_d)&=&\lambda_{\gamma}d^o\phi^{\gamma}(\epsilon(\beta(g)))(Y(\beta(g)),\\
\phi^{\gamma}(g)=\phi^{\gamma}(h)&=&0,
\end{eqnarray*}
for all $Y\in \Sec{\tau}$, where $d^o$ is the standard
differential on $\Gamma$. This algorithm is a generalization of
the Shake algorithm for holonomic systems (see
\cite{CoSMa,LeRe,mawest,McLPerl} for similar results on the pair
groupoid $Q\times Q$).

\subsection{Discrete nonholonomic Lagrangian systems on the pair
groupoid} Let $(L_{d}, {\mathcal M}_{c}, {\mathcal D}_{c})$ be a
discrete nonholonomic Lagrangian system on the pair group\-oid $Q
\times Q \rightrightarrows Q$ and suppose that $(q_{0}, q_{1})$ is
a point of ${\mathcal M}_{c}$. Then, using the results of Section
\ref{DMLg}, we deduce that $((q_{0}, q_{1}), (q_{1}, q_{2})) \in
(Q \times Q)_{2}$ is a solution of the discrete nonholonomic
Euler-Lagrange equations for $(L_{d}, {\mathcal M}_{c}, {\mathcal
D}_{c})$ if and only if
\[
\begin{array}{l}(D_{2}L_{d}(q_{0}, q_{1}) + D_{1}L_{d}(q_{1},
q_{2}))_{|{\mathcal D}_{c}(q_{1})} = 0,\\
(q_{1}, q_{2}) \in {\mathcal M}_{c},
\end{array}
\]
or, equivalently,
\[
\begin{array}{l}D_{2}L_{d}(q_{0}, q_{1}) + D_{1}L_{d}(q_{1}, q_{2}) = \displaystyle \sum_{j=1}^{n-r}
\lambda_{j}A^{j}(q_{1}),\\
(q_{1}, q_{2}) \in {\mathcal M}_{c},
\end{array}
\]
where $\lambda_{j}$ are the Lagrange multipliers and $\{A^{j}\}$
is a local basis of the annihilator ${\mathcal D}_{c}^{0}$. These
equations were considered in \cite{CoSMa} and \cite{McLPerl}.

Note that if $(q_{1}, q_{2}) \in \Gamma = Q \times Q$ then, in
this particular case, $G^{L_{d}}_{(q_{1}, q_{2})}: T_{q_{1}}Q
\times T_{q_{2}}Q \to \R$ is just the $\R$-bilinear map
$(D_{2}D_{1}L_{d})(q_{1}, q_{2})$.

On the other hand, if $(q_{1}, q_{2}) \in {\mathcal M}_{c}$ we
have that
\[
\lvec{(TQ)}^{{\mathcal M}_{c}}_{(q_{1}, q_{2})} = \set{v_{q_{2}}
\in T_{q_{2}}Q}{(0, v_{q_{2}}) \in T_{(q_{1}, q_{2})}{\mathcal
M}_{c}},
\]
\[
\rvec{(TQ)}^{{\mathcal M}_{c}}_{(q_{1}, q_{2})} = \set{v_{q_{1}}
\in T_{q_{1}}Q}{(v_{q_{1}}, 0) \in T_{(q_{1}, q_{2})}{\mathcal
M}_{c}}.
\]

Thus, the system $(L_{d}, {\mathcal M}_{c}, {\mathcal D}_{c})$ is
regular if and only if for every $(q_{1}, q_{2}) \in {\mathcal
M}_{c}$ the following conditions hold:

$$\left.\begin{array}{l}\mbox{If }v_{q_{1}} \in \rvec{(TQ)}^{{\mathcal
M}_{c}}_{(q_{1}, q_{2})}\mbox{ and}\\[8pt] \langle D_{2}D_{1}L_{d}(q_{1},
q_{2})v_{q_{1}}, v_{q_{2}}
\rangle = 0,\;\;\forall v_{q_{2}} \in {\mathcal
D}_{c}(q_{2})\end{array}\right\} \Longrightarrow v_{q_{1}} = 0,$$

and

$$\left.\begin{array}{l}\mbox{ If } v_{q_{2}} \in \lvec{(TQ)}^{{\mathcal
M}_{c}}_{(q_{1}, q_{2})} \mbox{ and } \\[8pt] \langle D_{2}D_{1}L_{d}(q_{1},
q_{2})v_{q_{1}}, v_{q_{2}} \rangle = 0, \;\;\; \forall v_{q_{1}} \in
{\mathcal D}_{c}(q_{1}) \end{array}\right\}\Longrightarrow
v_{q_{2}} = 0.$$

The first condition was obtained in \cite{McLPerl} in order to
guarantee the existence of a unique local nonholonomic evolution
operator $\Upsilon_{nh}^{(L_{d}, {\mathcal M}_{c}, {\mathcal
D}_{c})}$ for the system $(L_{d}, {\mathcal M}_{c}, {\mathcal
D}_{c})$. However, in order to assure that $\Upsilon_{nh}^{(L_{d},
{\mathcal M}_{c}, {\mathcal D}_{c})}$ is a (local) diffeomorphism
one must assume that the second condition also holds.

\begin{example}[\emph{Discrete Nonholonomically Constrained  particle}]
{\rm
 Consider the discrete nonholonomic system determined by:
\begin{enumerate}
\item[a)] A discrete Lagrangian $L_d:\R^3\times \R^3\to \R$:
\[
L_d(x_0, y_0, z_0, x_1, y_1, z_1)=\frac{1}{2}\left[
\left(\frac{x_1-x_0}{h}\right)^2+ \left(\frac{y_1-y_0}{h}\right)^2
+ \left(\frac{z_1-z_0}{h}\right)^2\right] .
\]
\item[b)] A constraint distribution of $Q=\R^3$,
\[
 {\mathcal D}_c=\hbox{span} \left\{ X_1=\frac{\partial}{\partial
x}+y\frac{\partial}{\partial z}, X_2=\frac{\partial}{\partial
y}\right\}.
\]
\item[c)] A discrete constraint submanifold ${\mathcal M}_c$ of
$\R^3\times \R^3$ determined by the constraint
\[
\displaystyle{\phi(x_0, y_0, z_0, x_1, y_1,
z_1)=\frac{z_1-z_0}{h}-
\left(\frac{y_1+y_0}{2}\right)\left(\frac{x_1-x_0}{h}\right). }
\]
\end{enumerate}
The system $(L_d, {\mathcal M}_c,{\mathcal D}_c)$ is a discretization of
a classical continuous nonholonomic system: the nonholonomic free particle
(for a discussion on this continuous system see, for instance,
\cite{BlKrMaMu,Cort}). Note that if $E_{(\R^3\times \R^3)}\cong T\R^3$
is the Lie algebroid of the pair groupoid $\R^3\times \R^3 \rightrightarrows
\R^3$ then
\[
T_{(x_1,y_1,z_1,x_1,y_1,z_1)}{\mathcal M}_c\cap E_{(\R^3\times \R^3)}(x_1,y_1,z_1)
= {\mathcal D}_c (x_1,y_1,z_1).
\]
Since
\[\lvec{X_1}=\frac{\partial}{\partial
x_1}+y_1\frac{\partial}{\partial z_1},\quad
\rvec{X_1}=-\frac{\partial}{\partial
x_0}-y_0\frac{\partial}{\partial z_0},\quad
\lvec{X_2}=\frac{\partial}{\partial y_1},\quad
\rvec{X_2}=-\frac{\partial}{\partial y_0},
\]
then, the discrete nonholonomic equations are:
\begin{eqnarray}
\left( \frac{x_2-2x_1+x_0}{h^2} \right) + y_1\left(
\frac{z_2-2z_1+z_0}{h^2} \right)&=&0,\label{eq1}\\
\frac{y_2-2y_1+y_0}{h^2}&=&0,\label{eq2}
\end{eqnarray}
which together with the constraint equation determine a well posed system of
difference equations.

We have that
\[
\begin{array}{l}
D_2D_1L_d=-\frac{1}{h}\{ dx_0\wedge dx_1+dy_0\wedge dy_1+dz_0\wedge dz_1\}\\[8pt]
\begin{array}{rcl}
(\rvec{T\R^3})^{{\mathcal M}_c}_{(x_0,y_0,z_0,x_1,y_1,z_1)}&=&\{
a_0\frac{\partial}{\partial x_0}+b_0\frac{\partial}{\partial y_0}+
\kern2pt c_0\frac{\partial}{\partial z_0}\in T_{(x_0,y_0,z_0)}\R^3 \,/\, \\[5pt]
&&
c_0=\frac{1}{2} ( a_0 (y_1+y_0) - b_0 (x_1-x_0)) \}.
\end{array}\\[8pt]
\begin{array}{rcl}
(\lvec{T\R^3})^{{\mathcal M}_c}_{(x_0,y_0,z_0,x_1,y_1,z_1)}&=&\{
a_1\frac{\partial}{\partial x_1}+b_1\frac{\partial}{\partial y_1}+
\kern2pt c_1\frac{\partial}{\partial z_1}\in T_{(x_1,y_1,z_1)}\R^3 \,/\, \\[5pt]
&& c_1=\frac{1}{2} ( a_1 (y_1+y_0) + b_1 (x_1-x_0)) \}.
\end{array}
\end{array}
\]
Thus, if we consider the open subset of ${\mathcal M}_c$ defined
by
\[
\set{(x_0,y_0,z_0,x_1,y_1,z_1)\in {\mathcal M}_c}{
2+y_1^2+y_1y_0\neq 0, 2+y_0^2+y_0y_1\neq 0}
\]
then in this subset the discrete nonholonomic system is regular.

Let $\Psi: {\frak g}=\R^2\longrightarrow {\frak X}(\R^3)$ given by
$\Psi(a, b)=a\frac{\partial}{\partial x}+b\frac{\partial}{\partial
z} $. Then ${\frak g}^{{\mathcal D}_c}=\hbox{span} \{
\Psi(\tilde{\xi})=X_1\}$, where $\tilde{\xi}: \R^3\to \R^2$ is
defined by $\tilde{\xi}(x, y, z)=(1, y)$. Moreover, the Lagrangian
$L_d$ is ${\frak g}$-invariant with respect to $\Psi$. Therefore,
\[
\begin{array}{l}
\kern-55pt J^{nh}_{\tilde{\xi}}(x_1, y_1, z_1, x_2, y_2, z_2)-
J^{nh}_{\tilde{\xi}}(x_0, y_0, z_0, x_1, y_1,
z_1)\\ \kern70pt=\lvec{\Psi( 0, y_2-y_1)}(x_1, y_1, z_1, x_2, y_2, z_2)(L_d),
\end{array}
\]
that is,
\[
\left(\frac{x_2-x_1}{h^2}+y_2
\frac{z_2-z_1}{h^2}\right)-\left(\frac{x_1-x_0}{h^2}+y_1
\frac{z_1-z_0}{h^2}\right)=(y_2-y_1)\left(\frac{z_2-z_1}{h^2}\right).
\]
This equation is precisely Equation (\ref{eq1}). }

\end{example}

\subsection{Discrete nonholonomic Lagrangian systems on a Lie
group} Let $G$ be a Lie group. $G$ is a Lie groupoid over a single
point and the Lie algebra ${\frak g}$ of $G$ is just the Lie
algebroid associated with $G$.

If $g, h \in G$, $v_{h} \in T_{h}G$ and $\alpha_{h} \in
T_{h}^{*}G$ we will use the following notation:
\[\begin{array}{rclrcl}
g v_{h} &= &(T_{h}l_{g})(v_{h}) \in T_{gh}G,& v_{h} g &=&
(T_{h}r_{g})(v_{h}) \in T_{hg}G,\\
g \alpha_{h} &=& (T_{gh}^{*}l_{g^{-1}})(\alpha_{h}) \in
T_{gh}^{*}G, & \alpha_{h} g &=& (T_{hg}^{*}r_{g^{-1}})(\alpha_{h})
\in T_{hg}^{*}G.
\end{array}
\]

Now, let $(L_{d}, {\mathcal M}_{c}, {\mathcal D}_{c})$ be a
discrete nonholonomic Lagrangian system on the Lie group $G$, that
is, $L_{d}: G \to \R$ is a discrete Lagrangian, ${\mathcal M}_{c}$
is a submanifold of $G$ and ${\mathcal D}_{c}$ is a vector
subspace of ${\frak g}$.

If $g_{1} \in {\mathcal M}_{c}$ then $(g_{1}, g_{2}) \in G \times
G$ is a solution of the discrete nonholonomic Euler-Lagrange
equations for $(L_{d}, {\mathcal M}_{c}, {\mathcal D}_{c})$ if and
only if
\begin{equation}\label{qer}
\begin{array}{l} g_{1}^{-1} dL_{d}(g_{1}) - dL_{d}(g_{2})g_{2}^{-1} =
\displaystyle \sum_{j=1}^{n-r} \lambda_{j}\mu^{j},\\
g_{k} \in {\mathcal M}_{c},\ k=1,2\end{array}
\end{equation} where $\lambda_{j}$ are the Lagrange multipliers and
$\{\mu^{j}\}$ is a basis of the annihilator ${\mathcal D}_{c}^{0}$
of ${\mathcal D}_{c}$. These equations were obtained in
\cite{McLPerl} (see Theorem 3 in \cite{McLPerl}).

Taking $p_k=dL_{d}(g_k)g_k^{-1}$, $k=1,2$ then
\begin{equation}\label{qera}
\begin{array}{l}
p_{2}-Ad^*_{g_{1}}p_1=
\displaystyle -\sum_{j=1}^{n-r} \lambda^{j}\mu_{j},\\
g_{k} \in {\mathcal M}_{c}, k=1, 2
\end{array}
\end{equation}
where $\map{Ad}{G\times\mathfrak{g}}{\mathfrak{g}}$ is the adjoint
action of $G$ on $\mathfrak{g}$. These equations were obtained in
\cite{FZ} and called \emph{discrete Euler-Poincar\'e-Suslov
equations}.

On the other hand, from (\ref{OmegaLd}), we have that
\[
\Omega _{L_d} ((\rvec{\eta},\lvec{\mu}),(\rvec{\eta}',
\lvec{\mu}'))=\rvec{\eta}'(\lvec{\mu} (L_d))-
\rvec{\eta }(\lvec{\mu}' (L_d)).
\]
Thus, if $g \in G$ then, using (\ref{1.21'}), it follows that the
$\R$-bilinear map $G_{g}^{L_{d}}: {\frak g} \times {\frak g}
\to \R$ is given by
\[
G_{g}^{L_{d}}(\xi, \eta) = -\lvec{\eta}(g)(\rvec{\xi}(L_{d})).
\]
Therefore, the system $(L_{d}, {\mathcal M}_{c}, {\mathcal D}_{c})$ is
regular if and only if for every $g \in {\mathcal M}_{c}$ the
following conditions hold:

$$\eta \in {\frak g} / \lvec{\eta}(g) \in T_{g}{\mathcal M}_{c}
\mbox{ and } \lvec{\eta}(g)(\rvec{\xi}(L_{d})) = 0, \forall \xi
\in {\mathcal D}_{c} \Longrightarrow \eta = 0,$$

$$
\xi \in {\frak g} / \rvec{\xi}(g) \in T_{g}{\mathcal M}_{c} \mbox{
and } \lvec{\eta}(g)(\rvec{\xi}(L_{d})) = 0, \forall \eta \in
{\mathcal D}_{c} \Longrightarrow \xi = 0.$$

We illustrate this situation with two simple examples previously
considered in \cite{FZ}.

\subsubsection{The discrete Suslov system} (See \cite{FZ})
The Suslov system studies the motion of a rigid body suspended at
its centre of mass under the action of the following nonholonomic
constraint: the body angular velocity is orthogonal to some fixed
direction.

The configuration space is $G=SO(3)$ and the elements of the Lie
algebra $\frak{so} (3)$ may be identified with $\R^3$ and
represented by coordinates $(\omega_x, \omega_y,\omega_z)$.
Without loss of generality, let us choose as fixed direction the
third vector of the body frame $\bar{e}_1, \bar{e}_2, \bar{e}_3$.
Then, the nonholonomic constraint is $\omega_z=0$.

The discretization of this system is modelled by considering the
discrete Lagrangian $L_d: SO(3)\longrightarrow \R$ defined by
$L_d(\Omega)=\frac{1}{2} \hbox{Tr } (\Omega J)$, where $J$
represents the mass matrix (a symmetric positive-definite matrix
with components $(J_{ij})_{1\leq i, j\leq 3}$).

The constraint submanifold ${\mathcal M}_c$ is determined by the
constraint $\hbox{Tr }(\Omega E_3)=0$ (see \cite{FZ}) where
\[ E_1=\left(
\begin{array}{lll}
0&0&0\\
0&0&-1\\
0&1&0
\end{array}
\right), \quad E_2=\left(
\begin{array}{lll}
0&0&1\\
0&0&0\\
-1&0&0
\end{array}
\right), \quad E_3=\left(
\begin{array}{lll}
0&-1&0\\
1&0&0\\
0&0&0
\end{array}
\right).
\]
is the standard basis of $\frak{so}(3)$, the Lie algebra of
$SO(3)$.

 The vector subspace ${\mathcal D}_c=\hbox{span} \{ E_1,E_2\}$. Therefore, ${{\mathcal
D}^0_c}=\hbox{span} \{ E^3\}$. Moreover, the exponential map of
$SO(3)$ is a diffeomorphism from an open subset of ${\mathcal
D}_c$ (which contains the zero vector) to an open subset of
${\mathcal M}_c$ (which contains the identity element $I$). In
particular, $T_I {\mathcal M}_c={\mathcal D}_c$.

On the other hand, the \emph{discrete Euler-Poincar\'e-Suslov
equations} are gi\-ven by
\begin{eqnarray*}
\lvec{E_i}(\Omega_1)(L_d)-\rvec{E_i}(\Omega_2)(L_d)=0, \; \; \;
\hbox{Tr }(\Omega_{i} E_3)=0, \; \; \; i \in \{1,2\}.
\end{eqnarray*}
After some straightforward operations, we deduce that the above
equations are equivalent to:
\begin{eqnarray*}
\hbox{Tr }\left((E_i\Omega_2-\Omega_1E_i) J\right)=0, \; \; \;
\hbox{Tr }(\Omega_{i} E_3)=0, \; \; \;  i \in \{1,2\}
\end{eqnarray*}
or, considering the components $\Omega_k=(\Omega^{(k)}_{ij})$ of
the elements of $SO(3)$, we have that:
\begin{eqnarray*}
\left(\begin{array}{l}
J_{23}\Omega^{(1)}_{33}-J_{33}\Omega^{(1)}_{32}+J_{22}\Omega^{(1)}_{23}\\
-J_{23}\Omega^{(1)}_{22}
+J_{12}\Omega^{(1)}_{13}-J_{13}\Omega^{(1)}_{12}
\end{array}\right)&=&\left(
\begin{array}{l}-J_{23}\Omega^{(2)}_{33}-J_{22}\Omega^{(2)}_{32}-J_{12}\Omega^{(2)}_{31}\\
+J_{33}\Omega^{(2)}_{23}
+J_{23}\Omega^{(2)}_{22}+J_{13}\Omega^{(2)}_{21}
\end{array}\right)\\
\left(\begin{array}{l}-J_{13}\Omega^{(1)}_{33}+J_{33}\Omega^{(1)}_{31}-J_{12}\Omega^{(1)}_{23}\\+J_{23}\Omega^{(1)}_{21}
-J_{11}\Omega^{(1)}_{13}+J_{13}\Omega^{(1)}_{11}\end{array}\right)&=&\left(\begin{array}{l}
J_{13}\Omega^{(2)}_{33}+J_{12}\Omega^{(2)}_{32}+J_{11}\Omega^{(2)}_{31}\\-J_{33}\Omega^{(2)}_{13}
-J_{23}\Omega^{(2)}_{12}-J_{13}\Omega^{(2)}_{11}\end{array}\right)\\
\Omega^{(1)}_{12}=\Omega^{(1)}_{21}, &\, &
\Omega^{(2)}_{12}=\Omega^{(2)}_{21}.
\end{eqnarray*}

Moreover, since the discrete Lagrangian verifies that
\[
L_d(\Omega)=\frac{1}{2} \hbox{Tr } (\Omega J)=\frac{1}{2} \hbox{Tr
} (\Omega^t J)=L_d(\Omega^{-1})
\]
and also the constraint satisfies $\hbox{Tr }(\Omega E_3)=-\hbox{Tr
}(\Omega^{-1} E_3)$, then this discretization of the Suslov system
is reversible. The regularity condition in  $\Omega \in SO(3)$ is
in this particular case:
\[
\eta\in {\frak so}(3)\; / \hbox{Tr }(E_1\Omega \eta J)=0, \quad
\hbox{Tr }(E_2\Omega \eta J)=0 \hbox{ and } \hbox{Tr }(\Omega\eta
E_3)=0\Longrightarrow \eta=0
\]
It is easy to show that the system is regular in a neighborhood of
the identity $I$.

\subsubsection{The discrete Chaplygin sleigh} (See \cite{F, FZ})\label{s4.3.2}
The Chaplygin sleigh system describes the motion of a rigid body
sliding on a horizontal plane. The body is supported at three
points, two of which slide freely without friction while the third
is a knife edge, a constraint that allows no motion orthogonal to
this edge (see \cite{NF}).

The configuration space of this system is the group $SE(2)$ of Euclidean
motions of $\R^2$. An element $\Omega \in SE(2)$ is represented by a
matrix
\[
\Omega = \left(
\begin{array}{ccc}
\cos\theta &-\sin\theta & x\\
\sin\theta &\cos\theta & y\\
0&0&1
\end{array}
\right) \qquad\mbox{ with }\theta ,x,y\in \R .
\]
Thus, $(\theta ,x,y)$ are local coordinates on $SE(2)$.

A basis of the Lie algebra $\mathfrak{se}(2)\cong \R ^3$ of $SE(2)$
is given by
\[
e=\left(
\begin{array}{ccc}
0&-1&0\\
1&0&0\\
0&0&0
\end{array}
\right),\qquad
e_1=\left(
\begin{array}{ccc}
0&0&1\\
0&0&0\\
0&0&0
\end{array}
\right) ,\qquad e_2=\left(
\begin{array}{ccc}
0&0&0\\
0&0&1\\
0&0&0
\end{array}
\right)
\]
and we have that
\[
[e,e_1]=e_2,\quad [e,e_2]=-e_1,\quad [e_1,e_2]=0.
\]
An element $\xi \in \mathfrak{se}(2)$ is of the form
\[
\xi =\omega \, e+v_1\, e_1+v_2\, e_2
\]
and the exponential map $\exp: \mathfrak{se}(2)\cong \R^3\to
SE(2)$ of $SE(2)$ is given by
\[
\exp (\omega, v_1,v_2)=(\omega ,v_1\frac{\sin \omega}{\omega}+v_2 ( \frac{\cos \omega -1}{\omega } ),
-v_1 ( \frac{\cos \omega -1}{\omega} )+v_2  \frac{\sin \omega}{\omega }  ),\mbox{ if }\omega \neq 0,
\]
and
\[
\exp (0, v_1,v_2)=(0,v_1,v_2).
\]
Note that the restriction of this map to the open subset $U=]-\pi ,\pi [\times \R^2\subseteq \R ^3
\cong \mathfrak{se}(2)$ is a diffeomorphism onto the open subset $\exp (U)$ of $SE(2)$.

A discretization of the Chaplygin sleigh may be constructed as follows:
\begin{itemize}
\item[-] The discrete Lagrangian $L_d: SE(2)\longrightarrow \R$ is given by
\[
L_d(\Omega)=\frac{1}{2} \hbox{Tr } ( \Omega \J \Omega^T)-\hbox{Tr
} ( \Omega \J),
\]
 where $\J$ is the matrix:
\[
\J=\left(
\begin{array}{ccc}
(J/2)+ma^2&mab&ma\\
mab&(J/2)+mb^2&mb\\
ma&mb&m
\end{array}
\right)
\]
(see \cite{FZ}).
\item[-] The vector subspace ${\mathcal D}_c$ of $\mathfrak{se}(2)$ is
\[
{\mathcal D}_c=\hbox{ span } \{ e, e_1\}= \set{(\omega, v_1,
v_2)\in \mathfrak{se}(2)}{v_2=0}.
\]
\item[-] The constraint submanifold ${\mathcal M}_c$ of
$SE(2)$ is
\begin{equation}\label{defM}
{\mathcal M}_c=\exp (U\cap {\mathcal D}_c).
\end{equation}
Thus, we have that
\begin{eqnarray*}
{\mathcal M}_c&=&\set{(\theta , x,y )\in SE(2)}{-\pi <\theta
<\pi,\, \theta \neq 0, (1-\cos \theta )x-y\sin \theta =0}\\ &&
\cup \set{(0,x,0)\in SE(2)}{x\in\R}.
\end{eqnarray*}
From (\ref{defM}) it follows that $I\in {\mathcal M}_c$ and
$T_I{\mathcal M}_c= {\mathcal D}_c$. In fact, one may prove that
\[
T_{(0,x,0)}{\mathcal M}_c=\hbox{ span }\{ \frac{\partial}{\partial \theta}_{|(0,x,0)}
+\frac{x}{2}\frac{\partial}{\partial y}_{|(0,x,0)},\frac{\partial}{\partial x}_{|(0,x,0)}\},
\]
for $x\in \R$.
\end{itemize}
\begin{figure}
   \includegraphics[width=6cm]{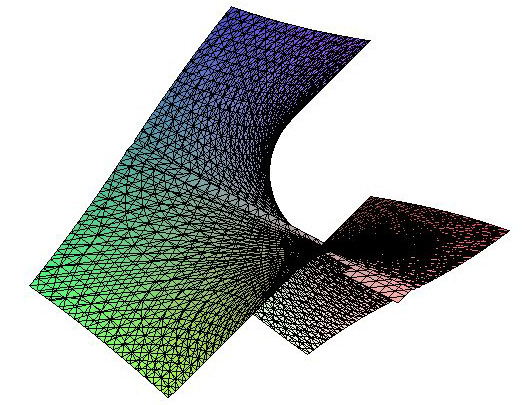}
  \caption{Submanifold ${\mathcal M}_c$ }\label{Figure2}
\end{figure}
Now, the discrete Euler-Poincar\'e-Suslov equations are:
\begin{eqnarray*}
\lvec{e}(\theta _1,x_1, y_1)(L_d)-\rvec{e}(\theta _2,x_2,y_2)(L_d)=0,\\
\quad \lvec{e_1}(\theta _1,x_1, y_1)(L_d)-\rvec{e_1}(\theta _2,x_2,y_2)(L_d)=0,
\end{eqnarray*}
and the condition $(\theta _k,x_k,y_k)\in {\mathcal M}_c$, with $k\in \{1,2\}$.
We rewrite these equations as the following system of difference
equations:
\begin{eqnarray*}
\left(\begin{array}{l}
-am\cos\theta _1-bm\sin\theta _1+am\\
+mx_1\cos\theta _1+my_1\sin\theta _1
\end{array}\right)
&\kern-5pt=&\kern-5pt
\left(\begin{array}{l}
mx_{2}+am\cos\theta_{2}\\
-bm\sin\theta_{2}-am \end{array}\right)
\\
\left(\begin{array}{l}
 am\, y_{1}\cos\theta _1-am\,x_1\sin \theta _1
 -bm\,x_1\cos\theta _1\\
-bm\,y_1\sin\theta _1+(a^2m+b^2m+J)\sin \theta _1
\end{array}\right)
&\kern-5pt=&\kern-5pt
\left(\begin{array}{l} am\,y_{2}-bm\,x_2\\
+(a^2m+b^2m+J)\sin \theta _2
\end{array}\right)
\end{eqnarray*}
together with the condition
\[
(\theta _k,x_k,y_k)\in {\mathcal M}_c,\,\, k\in\{1,2\}.
\]
On the other hand, one may prove that the discrete nonholonomic
Lagrangian system $(L_{d}, {\mathcal M}_{c}, {\mathcal D})$ is
reversible.

 Finally, consider a point $(0,x,0)\in {\mathcal M}_c$ and an
 element $\eta \equiv (\omega ,v_1,v_2)\in \mathfrak{se}(2)$ such that
\[
\lvec{\eta}(0,x,0)\in T_{(0,x,0)}{\mathcal M}_c,\,\,
\lvec{\eta}(0,x,0)(\rvec{e}(L_d))=0, \, \,
\lvec{\eta}(0,x,0)(\rvec{e_1}(L_d))=0.
\]
Then, if we assume that $a^2m+J+am\frac{x}{2}\neq 0$ it follows
that $\eta =0$.

Thus, the discrete nonholonomic Lagrangian system $(L_d,
{\mathcal M}_c,{\mathcal D}_c)$ is regular in a neighborhood of
the identity $I$.

\subsection{Discrete nonholonomic Lagrangian systems on an action
Lie group\-oid} \label{TaLg}
Let $H$ be a Lie group with identity element ${\mathfrak e}$ and
$\cdot:M\times H\to M$, $(x,h)\in M\times H\mapsto xh,$ a right
action of $H$  on $M$. Thus, we may consider the action Lie
groupoid $\Gamma =M\times H$ over $M$ with structural maps given by
\begin{equation}\label{*tilde}
\begin{array}{l}
\tilde{\alpha}(x,h)=x,\;\;\; \tilde{\beta}(x,h)=xh,\;\;\;
\tilde{\epsilon }(x)=(x,{\mathfrak e}),\\
\tilde{m}((x,h),(xh,h'))=(x,hh'),\;\;\; \tilde{i}(x,h)=(xh,
h^{-1}).
\end{array}
\end{equation}
Now, let ${\mathfrak h}=T_{{\mathfrak e}}H$ be the Lie algebra of
$H$ and $\Phi:{\mathfrak h}\to {\mathfrak X}(M)$ the map given by
\[
\Phi(\eta)=\eta_M,\;\;\;\; \mbox{for } \eta\in {\mathfrak h},
\]
where $\eta_M$ is the infinitesimal generator of the action
$\cdot:M\times H\to M$ corresponding to $\eta$. Then, $\Phi$ is a
Lie algebra morphism and the corresponding action Lie algebroid
$pr_1:M\times {\mathfrak{ h}}\to M$ is just the Lie algebroid of
$\Gamma=M\times H$.

We have that $\Sec{pr_1}\cong \set{\tilde{\eta}:M\to {\mathfrak
h}}{\tilde{\eta} \mbox{ is  smooth }}$ and that the Lie algebroid
structure $(\lcf\cdot,\cdot\rcf_\Phi,\rho_\Phi)$ on $pr_1:M\times
H\to M$ is defined by
\[
\lcf
\tilde{\eta},\tilde{\mu}\rcf_{\Phi}(x)=[\tilde{\eta}(x),\tilde{\mu}(x)]
+
(\tilde{\eta}(x))_M(x)(\tilde{\mu})-(\tilde{\mu}(x))_M(x)(\tilde{\eta}),\;\;\;
\rho_\Phi(\tilde{\eta})(x)=(\tilde{\eta}(x))_M(x),
\]
for $\tilde{\eta}, \tilde{\mu}\in \Sec{pr_1}$ and $x\in M.$
Here, $[\cdot,\cdot]$ denotes the Lie bracket of ${\mathfrak h}$.

If $(x,h)\in \Gamma=M\times H$ then the left-translation
$l_{(x,h)}:\tilde{\alpha}^{-1}(xh)\to \tilde{\alpha}^{-1}(x)$ and
the right-translation $r_{(x,h)}:\tilde{\beta}^{-1}(x)\to
\tilde{\beta}^{-1}(xh)$ are given
\begin{equation}\label{tilde+}
l_{(x,h)}(xh,h')=(x,hh'),\;\;\;
r_{(x,h)}(x(h')^{-1},h')=(x(h')^{-1},h'h).
\end{equation}

Now, if $\eta\in {\mathfrak h}$ then $\eta$ defines a constant
section $C_\eta:M\to {\mathfrak h}$ of $pr_1:M\times {\mathfrak
h}\to M$ and, using (\ref{linv}), (\ref{rinv}), (\ref{*tilde}) and
(\ref{tilde+}), we have that the left-invariant and the
right-invariant vector fields $\lvec{C}_\eta$ and $\rvec{C}_\eta$,
respectively, on $M\times H$ are defined by
\begin{equation}\label{*+}
\rvec{C}_\eta(x,h)=(-\eta_M(x),\rvec{\eta}(h)),\;\;\;\;
\lvec{C}_\eta(x,h)=(0_x,\lvec{\eta}(h)),
\end{equation}
for $(x,h)\in \Gamma=M\times H.$

Note that if $\{\eta_i\}$ is a basis of ${\mathfrak h}$ then
$\{C_{\eta_i}\}$ is a global basis of $\Sec{pr_1}.$

On the other hand, we will denote by $\exp _\Gamma :E_\Gamma =
M\times \frak h\to \Gamma =M\times H$ the map given by
\[
\exp _\Gamma (x,\eta )=(x,\exp _H (\eta )),\quad \mbox{for }
(x,\eta )\in E_\Gamma =M\times \frak h,
\]
where $\exp _H:\frak h\to H$ is the exponential map of the Lie
group $H$. Note that if $\Phi _{(x,e)}:\R\to \Gamma =M\times H$ is
the integral curve of the left-invariant vector field
$\lvec{C}_\eta$ on $\Gamma =M\times H$ such that
$\Phi _{(x,e)}(0)=(x,e)$ then (see (\ref{*+}))
\[
\exp _\Gamma (x,\eta ) =\Phi _{(x,e)}(1).
\]
Next, suppose that $L_d: \Gamma=M\times H\to \R$ is a Lagrangian
function, ${\mathcal D}_c$ is a constraint distribution such that
$\{ X^\alpha \}$ is a local basis of sections of the annihilator
${\mathcal D}^0_c$, and ${\mathcal M}_c\subseteq \Gamma$ is the
discrete constraint submanifold.

For every $h\in H$ (resp., $x\in M$) we will denote by $L_h$
(resp., $L_x$) the real function on $M$ (resp., on $H$) given by
$L_h(y)=L_d(y,h)$ (resp., $L_x(h')=L_d(x,h'))$. A composable
pair $((x,h_k),(xh_k,h_{k+1}))\in \Gamma_2\cap ({\mathcal
M}_c\times {\mathcal M}_c)$ is a solution of the discrete
nonholonomic Euler-Lagrange equations for the system
$(L_d,{\mathcal M}_c, {\mathcal D}_c)$ if
\[
\lvec{C}_\eta(x,h_k)(L_d)-\rvec{C}_\eta(xh_k,h_{k+1})(L_d)=\lambda
_\alpha X^\alpha (xh_k) (\eta) , \mbox{ for all } \eta\in {\mathfrak
h},
\]

or, in other terms (see (\ref{*+}))
\[
\{(T_{{\mathfrak e}} l_{h_k})(\eta)\}(L_x)-\{(T_{{\mathfrak e}}
r_{h_{k+1}})(\eta)\}(L_{xh_k})+\eta_M(xh_k)(L_{h_{k+1}})=\lambda
_\alpha X^\alpha (xh_k) (\eta),
\]
for all $\eta\in {\mathfrak h.}$

\subsubsection{The discrete Veselova system}
As a concrete example of a nonholonomic system on a transformation
Lie groupoid we consider a discretization of the Veselova system
(see \cite{VeVe}). In the continuous theory~\cite{CoLeMaMa}, the
configuration manifold is the transformation Lie algebroid
$pr_1:S^2\times\mathfrak{so}(3)\to S^2$ with Lagrangian
\[
L_c(\gamma,\omega)=\frac{1}{2}\omega\cdot
I\omega-mgl\gamma\cdot\e,
\]
where $S^2$ is the unit sphere in $\R^3$, $\omega\in\R^3\simeq\mathfrak{so}(3)$ is the angular
velocity, $\gamma$ is the direction opposite to the gravity and
$\e$ is a unit vector in the direction from the fixed point to the
center of mass, all them expressed in a frame fixed to the body.
The constants $m$, $g$ and $l$ are respectively the mass of the
body, the strength of the gravitational acceleration and the
distance from the fixed point to the center of mass. The matrix
$I$ is the inertia tensor of the body. Moreover, the constraint
subbundle ${\mathcal D}_c \to S^2$ is
given by
\[
\gamma \in S^2 \mapsto {\mathcal D}_c{}(\gamma )=\set{\omega \in
\R^3\simeq \mathfrak{so}(3)}{\gamma \cdot \omega =0}.
\]
Note that the section $\phi :S^2 \to S^2\times
\mathfrak{so}(3)^*$, $(x,y,z)\mapsto ((x,y,z),xe^1+ye^2+ze^3)$,
where $\{e_1,e_2,e_3\}$ is the canonical basis of $\R ^3$ and
$\{e^1,e^2,e^3\}$ is the dual basis, is a global basis for
${\mathcal D}_c^0$.

If $\omega \in \mathfrak{so}(3)$ and $\widehat{\omega}$
is the skew-symmetric matrix of order 3 such that
$\widehat{\omega}v=\omega \times v$ then the Lagrangian
function $L_c$ may be expressed as follows
\[
L_c(\gamma,\omega )=\frac{1}{2}\tr(\widehat{\omega} \I
\widehat{\omega} ^T)-m\,g\,l\,\gamma\cdot\e ,
\]
where $\I = \frac{1}{2} \tr (I)I_{3 \times 3}-I$. Here, $I_{3
\times 3}$ is the identity matrix. Thus, we may define a discrete
Lagrangian $L_d: \Gamma = S^2\times SO(3)\to\R$ for the system by
(see \cite{groupoid})
\[
L_d(\gamma,\Omega )=-\frac{1}{h}\tr(\I
\Omega)-h\,m\,g\,l\,\gamma\cdot\e.
\]
On the other hand, we consider the open subset of
$SO(3)$
\[
V=\set{\Omega \in SO(3)}{\tr \Omega \neq \pm 1}
\]
and the real function $\psi :S^2\times V\to \R$ given
by
\[
\psi (\gamma, \Omega )=\gamma\cdot (\widehat{\Omega -\Omega ^T}).
\]
One may check that the critical points of $\psi$ are
\[
C_\psi =\set{(\gamma ,\Omega )\in S^2\times
V}{\Omega\gamma-\gamma=0}.
\]
Thus, the subset ${\mathcal M}_c$ of $\Gamma =S^2\times SO(3)$
defined by
\[
{\mathcal M}_c= \set{(\gamma ,\Omega )\in (S^2\times
V)-C_\psi}{\gamma \cdot (\widehat{\Omega -\Omega ^T})=0},
\]
is a submanifold of $\Gamma$ of codimension one. ${\mathcal M}_c$ is
the discrete constraint submanifold.

We have that the map $\exp _\Gamma :S^2\times \mathfrak{so}(3)\to
S^2\times SO(3)$ is a diffeomorphism from an open subset of ${\mathcal D}_c$,
which contains the zero section, to an open subset of ${\mathcal M}_c$, which
contains the subset of $\Gamma$ given by
\[
\tilde{\epsilon} (S^2)=\{(\gamma ,e)\in S^2\times SO(3)\}.
\]
So, it follows that
\[
({\mathcal D}_c)(\gamma )=T_{(\gamma ,e)}{\mathcal M}_c\cap E_\Gamma (\gamma),
\qquad \mbox{ for }\gamma \in S^2.
\]
Following the computations of \cite{groupoid} we get the
nonholonomic discrete Euler-Lagrange equations, for $((\gamma
_k,\Omega _k),(\gamma _{k}\Omega _k,\Omega _{k+1}))\in \Gamma _2$
\[
\begin{array}{l}
M_{k+1}-\Omega_k^T M_k\Omega_k
+mglh^2(\widehat{\gamma_{k+1}\times\e})=\lambda \widehat{\gamma
_{k+1}},\\
\gamma _k (\widehat{\Omega _k-\Omega _k^T})=0,\,\,
\gamma _{k+1} (\widehat{\Omega _{k+1}-\Omega _{k+1}^T})=0,
\end{array}
\]
where $M=\Omega \I-\I \Omega^T$. Therefore, in terms of the axial
vector $\Pi$ in $\R^3$ defined by $\hat{\Pi}=M$, we can write the
equations in the form
\[
\begin{array}{l}
\Pi_{k+1}=\Omega _k^T\Pi_k-mglh^2\gamma_{k+1}\times\e +\lambda
\gamma_{k+1},\\
\gamma _k (\widehat{\Omega _k-\Omega _k^T})=0,\,\,
\gamma _{k+1} (\widehat{\Omega _{k+1}-\Omega _{k+1}^T})=0.
\end{array}
\]
Note that, using the expression of an arbitrary element of $SO(3)$ in
terms of the Euler angles (see Chapter 15 of \cite{MaRa}), we deduce
that the discrete constraint submanifold ${\mathcal M}_c$ is reversible,
that is, $i({\mathcal M}_c)={\mathcal M}_c$. However, the discrete
nonholonomic Lagrangian system $(L_d,{\mathcal M}_c,{\mathcal D}_c)$
is not reversible. In fact, it is easy to prove that $L_d\circ i\neq L_d$.

On the other hand, if $\gamma \in S^2$ and $\xi, \eta \in \R^3\cong \mathfrak{so}(3)$
then it follows that
\[
\rvec{C}_\xi (\gamma ,I_3) (\lvec{C}_\eta (L_d))=-\xi\cdot I\eta.
\]
Consequently, the nonholonomic system $(L_d,{\mathcal
M}_c,{\mathcal D}_c)$ is regular in a neighborhood (in ${\mathcal
M}_c$) of the submanifold $\tilde{\epsilon} (S^2)$.

\subsection{Discrete nonholonomic Lagrangian systems on an Atiyah
Lie group\-oid}

Let $p: Q \to M=Q/G$ be a principal $G$-bundle and  choose a local
trivialization $G \times U$, where $U$ is an open subset of $M$.
Then, one may identify the open subset $(p^{-1}(U) \times
p^{-1}(U))/G \simeq ((G \times U)\times (G \times U))/G$ of the
Atiyah groupoid $(Q \times Q)/G$ with the product manifold $(U
\times U) \times G$. Indeed, it is easy to prove that the map
\[
((G \times U)\times (G \times U))/G \to (U \times U)\times G,
\]
\[
[((g, x), (g', y))] \to   ((x, y), g^{-1}g')),
\]
is bijective. Thus, the restriction to $((G \times U)\times (G
\times U))/G$ of the Lie groupoid structure on $(Q \times Q)/G$
induces a Lie groupoid structure in $(U \times U) \times G$ with
source, target and identity section given by
\[
\begin{array}{lr}
\alpha: (U \times U) \times G \to U; & ((x, y), g) \to x,
\\
\beta: (U \times U)\times G \to U; & ((x, y), g) \to y,
\\
\epsilon: U \to (U \times U) \times G; & x \to ((x, x), {\frak
e}),
\end{array}
\]
and with multiplication $m: ((U \times U)\times G)_{2} \to (U
\times U)\times G$ and inversion $i: (U \times U) \times G \to (U
\times U) \times G$ defined by
\begin{equation}
\label{mi}
\begin{array}{rcl}
 m(((x, y), g), ((y, z), h))& = &((x, z), g h),
\\
i((x, y), g) & = & ((y, x), g^{-1}).
\end{array}
\end{equation}
The Lie algebroid $A((U \times U)\times G)$ may be
identified with the vector bundle $TU \times
{\mathfrak g}\to U$. Thus, the fibre over the point
$x \in U$ is the vector space
$T_{x}U \times {\mathfrak g}$. Therefore, a section of $A((U \times
U)\times G)$ is a pair $(X, \tilde{\xi})$, where $X$ is a vector
field on $U$ and $\tilde{\xi}$ is a map from $U$ on ${\mathfrak
g}$. The space $\Sec{A((U \times U)\times G)}$ is generated by
sections of the form $(X, 0)$ and $(0, C_{\xi})$, with $X \in
{\mathfrak X}(U)$, $\xi \in {\mathfrak g}$ and $C_{\xi}: U \to
{\mathfrak g}$ being the constant map $C_{\xi}(x) = \xi$, for all
$x \in U$ (see \cite{groupoid} for more details).

Now, suppose that $L_d: (U \times U)\times G \to \R$ is a
Lagrangian function, ${\mathcal D}_c$ a vector subbundle of
$TU\times {\mathfrak g}$ and ${\mathcal M}_c$ a constraint
submanifold on $(U \times U)\times G$. Take a basis of sections
$\{Y^{\alpha}\}$ of the annihilator ${\mathcal D}_c^o$. Then, the
discrete nonholonomic equations are
\[
\lvec{(X_\alpha, \tilde{\eta}_\alpha)}((x, y), g_{k})(L_d) -
\rvec{(X_\alpha, \tilde{\eta}_\alpha)}((y, z),
g_{k+1})(L_d)  = 0,
\]
with $(X_\alpha,\tilde{\eta}_\alpha ):U\to TU\times \mathfrak{g}$
a basis of the space $\Sec{\tau _{{\mathcal D}_c}}$ and
$( ((x,y),g_k),((y,z),$ $g_{k+1}))\in ({\mathcal M}_c\times{\mathcal M}_c)
\cap ((U\times U)\times G)_2$. The above equations may be also written
as
\[
\begin{array}{lcr}
\lvec{(X, 0)}((x, y), g_{k})(L_d) - \rvec{(X, 0)}((y, z),
g_{k+1})(L_d) & =& \lambda_{\alpha} Y^{\alpha}(y)(X(y)),
\\
\lvec{(0, C_{\xi})}((x, y), g_{k})(L_d) - \rvec{(0, C_{\xi})}((y,
z), g_{k+1})(L_d) & = & \lambda_{\alpha}
Y^{\alpha}(y)(C_{\xi}(y)),
\end{array}
\]
with $X \in {\mathfrak X}(U)$, $\xi \in {\mathfrak g}$ and $(((x,
y), g_{k}), ((y, z), g_{k+1})) \in ({\mathcal M}_c\times {\mathcal
M}_c)\cap ((U \times U)\times G)_{2}$. An equivalent expression
of these equations is
\begin{equation}\label{disLP1}
\begin{array}{l}
D_{2}L_d((x, y), g_{k}) + D_{1}L_d ((y, z), g_{k+1}) =
\lambda_{\alpha}\mu^{\alpha}(y),
\\
p_{k+1}(y, z) = Ad^*_{g_{k}} p_{k}(x, y)-\lambda_{\alpha}
\tilde{\eta}^{\alpha}(y),
\end{array}
\end{equation}
where $p_{k}(\bar{x}, \bar{y}) = d(r_{g_{k}}^*L_{(\bar{x},
\bar{y}, \;)})({\mathfrak e}) $ for $(\bar{x}, \bar{y}) \in U
\times U$ and we write $Y^{\alpha}\equiv (\mu^{\alpha},
\tilde{\eta}^{\alpha})$, $\mu^{\alpha}$ being a 1-form on $U$ and
$\tilde{\eta}^{\alpha}: U \to {\mathfrak g}^*$ a smooth map.
\subsubsection{A discretization of the equations of motion of a rolling ball
without sliding on a rotating table with constant angular
velocity}

A (homogeneous) sphere of radius $r>0$,  mass $m$ and inertia
about any axis $I$ rolls without sliding on a horizontal table
which rotates with constant angular velocity $\Omega$ about a
vertical axis through one of its points. Apart from the constant
gravitational force, no other external forces are assumed to act
on the sphere (see \cite{NF}).

The configuration  space for the continuous system is
$Q=\R^2\times SO(3)$ and we shall use the notation $(x, y; R)$ to
represent a typical point in $Q$. Then, the nonholonomic
constraints are
\begin{eqnarray*}
\dot{x}+\frac{r}{2}\tr (\dot{R}R^T E_2)&=&-\Omega y,\\
\dot{y}-\frac{r}{2}\tr (\dot{R}R^T E_1)&=&\Omega x,
\end{eqnarray*}
where $\{E_1, E_2, E_3\}$ is the standard basis of $\frak{so}(3)$.

The matrix $\dot{R}R^T$ is skew symmetric, therefore we may write
\[
\dot{R}R^T = \left(
\begin{array}{ccc}
0&-w_3&w_2\\
w_3&0&-w_1\\
-w_2&w_1&0
\end{array}
\right)
\]
where $(w_1, w_2, w_3)$ represents the angular velocity vector of
the sphere measured with respect to the inertial frame. Then, we
may rewrite the constraints in the usual form:
\begin{eqnarray*}
\dot{x}-rw_2&=&-\Omega y,\\
\dot{y}+rw_1&=&\Omega x.
\end{eqnarray*}

The Lagrangian for the rolling ball is:
\begin{eqnarray*}
L_c(x, y; R, \dot{x},\dot{y};\dot{R})&=&\frac{1}{2}m(\dot{x}^2 +
\dot{y}^2) +\frac{1}{4}I \tr (\dot{R}R^T(\dot{R}R^T)^T)\\
&=&\frac{1}{2}m(\dot{x}^2 + \dot{y}^2) + \frac{1}{2}I (\omega_1^2
+ \omega_2^2 + \omega_3^2).
\end{eqnarray*}
Moreover, it is clear that $Q=\R^2\times SO(3)$ is the total space
of a trivial principal $SO(3)$-bundle over $\R^2$ and the bundle
projection $\phi:Q\to M=\R^2$ is just the canonical projection on
the first factor. Therefore, we may consider the corresponding
Atiyah algebroid $E'=TQ/SO(3)$ over $M=\R^2$. We will identify the tangent
bundle to $SO(3)$ with ${\frak {so}}(3)\times SO(3)$ by using
right translation.

Under this identification between $T(SO(3))$ and ${\frak {so}}(3)\times SO(3)$
the tangent action of $SO(3)$ on $T(SO(3))\cong
{\frak {so}}(3)\times SO(3)$ is the trivial action
\begin{equation}\label{Action}
({\frak {so}}(3)\times SO(3))\times SO(3)\to
{\frak {so}}(3)\times SO(3),\;\;\; ((\omega ,R),S)\mapsto (\omega,RS).
\end{equation}

Thus, the Atiyah algebroid $TQ/SO(3)$ is isomorphic to the product
manifold $T\R^2\times {\frak {so}}(3)$ and the vector bundle
projection is $\tau_{\R^2}\circ pr_1,$ where $pr_1:T\R^2\times
{\frak {so}}(3)\to T\R^2$ and $\tau_{\R^2}:T\R^2\to \R^2$ are the
canonical projections.

A section of $E'=TQ/SO(3)\cong T\R^2\times {\frak {so}}(3)\to
\R^2$ is a pair $(X,u)$, where $X$ is a vector field on $\R^2$ and
$u:\R^2\to {\frak {so}}(3)$ is a smooth map. Therefore, a global
basis of sections of $T\R^2\times {\frak {so}}(3)\to \R^2$ is
\[\begin{array}{rrr}
s_1'=(\displaystyle\frac{\partial}{\partial x},0),&
s_2'=(\displaystyle\frac{\partial}{\partial y},0),&\\[5pt]
s_3'=(0,E_1),& s_4'=(0,E_2),& s_5'=(0,E_3).
\end{array}
\]

The anchor map $\rho':E'=TQ/SO(3)\cong T\R^2\times {\frak
{so}}(3)\to T\R^2$ is the projection over the first factor and if
$\lcf\cdot, \cdot \rcf'$ is the Lie bracket on the space
$\Sec{E'=TQ/SO(3)}$ then the only non-zero fundamental Lie brackets
are
\[
\lcf s_3',s_4'\rcf'=s_5',\;\;\;\lcf s_4',s_5'\rcf'=s_3',\;\;\;
\lcf s_5',s_3'\rcf'=s_4'.
\]

Moreover, the Lagrangian function $L_c=T$ and the constraint
functions are $SO(3)$-invariant. Consequently, $L_c$ induces a
Lagrangian function $L'_c$ on $E'=TQ/SO(3)$
\begin{eqnarray*}
L_c'(x,y,\dot{x},\dot{y};\omega)&=&\frac{1}{2}m(\dot{x}^2 +
\dot{y}^2) + \frac{1}{4}I \tr (\omega \omega^T),\\
&=& \frac{1}{2}m(\dot{x}^2 +
\dot{y}^2) - \frac{1}{4}I \tr (\omega ^2),
\end{eqnarray*}
where $(x,y,\dot{x},\dot{y})$ are the standard coordinates on
$T\R^2$ and $\omega\in {\frak {so}}(3)$. The constraint functions
defined on $E'=TQ/SO(3)$ are:
\begin{equation}\label{ConstraintE'}
\begin{array}{rcl}
\dot{x}+\frac{r}{2}\tr(\omega E_2)&=&-\Omega y,\\
\dot{y}-\frac{r}{2}\tr(\omega E_1)&=&\Omega x.
\end{array}
\end{equation}
We have a nonholonomic system on  the Atiyah algebroid
$E'=TQ/SO(3)\cong T\R^2\times {\frak {so}}(3).$ This kind of
systems was recently analyzed by J. Cort\'es {\it et al}
\cite{CoLeMaMa} (in particular, this example was carefully
studied).

Eqs. (\ref{ConstraintE'}) define an affine subbundle of the vector bundle
$E'\cong T\R ^2\times \mathfrak{so}(3)\to \R^2$ which is modelled over
the vector subbundle ${\mathcal D}'_c$ generated by the sections
\[
{\mathcal D}'_c=\{s_5', rs_1'+s_4', rs_2'-s_3'\}.
\]
Our objective is to discretize this example directly on the Atiyah
algebroid. The Atiyah groupoid is now identified to $\R^2\times
\R^2\times SO(3)\rightrightarrows \R^2$. We may construct  the
discrete Lagrangian by
\[
L_d'(x_0, y_0, x_1, y_1; W_1)= L_c'(x_0,
y_0, \frac{x_1-x_0}{h},  \frac{y_1-y_0}{h}; (\log W_1)/h)
\]
where $\log: SO(3)\longrightarrow {\frak {so}}(3)$ is the
(local)-inverse of the exponential map $\hbox{exp}: {\frak
{so}}(3)\longrightarrow SO(3)$. For simplicity instead of this
procedure we use the following approximation:
\[
\log W_1/h\approx \frac{W_1-I_{3\times 3}}{h}
\]
where $I_{3\times 3}$ is the identity matrix.

Then
\begin{eqnarray*}
L'_d(x_0, y_0, x_1, y_1; W_1)&\kern-4pt=& \kern-4pt L_c'(x_0, y_0,
\frac{x_1-x_0}{h}, \frac{y_1-y_0}{h}; \frac{W_1-I_{3\times 3}}{h})\\
&\kern-4pt=&\kern-4pt\frac{1}{2}m\left[\left(\frac{x_1-x_0}{h}\right)^2 \kern-2pt+
\left(\frac{y_1-y_0}{h}\right)^2\right] \kern-2pt+\frac{I}{(2h)^2}
\tr (I_{3\times 3}-W_1)
\end{eqnarray*}
Eliminating constants, we may consider as discrete Lagrangian
\[
L'_d=\frac{1}{2}m\left[\left(\frac{x_1-x_0}{h}\right)^2 +
\left(\frac{y_1-y_0}{h}\right)^2\right] -\frac{I}{2h^2}
\tr (W_1)
\]

The {\bf discrete constraint submanifold} ${\mathcal M}'_c$ of
$\R^2\times \R^2\times SO(3)$ is determined by the constraints:
\begin{eqnarray*}
\frac{x_1-x_0}{h}+\frac{r}{2h}\tr (W_1E_2)&=&-\Omega \frac{y_1+y_0}{2},\\
\frac{y_1-y_0}{h}-\frac{r}{2h}\tr (W_1E_1)&=&\Omega
\frac{x_1+x_0}{2},
\end{eqnarray*}
We have that the system $(L'_d, {\mathcal M}'_c,{\mathcal D}'_c)$
is not reversible. Note that the Lagrangian function $L'_d$ is reversible.
However, the constraint submanifold ${\mathcal M}_c'$ is not reversible.

The discrete nonholonomic Euler-Lagrange equations for the system
$(L'_d, , {\mathcal M}'_c,$ $ {\mathcal D}'_c)$ are:
\begin{eqnarray*}
\lvec{s_5'}(x_0, y_0, x_1, y_1; W_1)(L'_d)- \rvec{s_5'}(x_1, y_1,
x_2, y_2; W_2)(L'_d)&=&0\\
\lvec{(rs_1'+s_4')}(x_0, y_0, x_1, y_1; W_1)(L'_d)-
\rvec{(rs_1'+s_4')}(x_1, y_1,
x_2, y_2; W_2)(L'_d)&=&0\\
\lvec{(rs_2'-s_3')}(x_0, y_0, x_1, y_1; W_1)(L'_d)-
\rvec{(rs_2'-s_3')}(x_1, y_1, x_2, y_2; W_2)(L'_d)&=&0
\end{eqnarray*}
with the constraints defining ${\mathcal M}_c$.

On the other hand, the vector fields $\lvec{s}'_5$,
$\rvec{s}_5'$, $\lvec{rs'_1+s_4'}$, $\rvec{rs'_1+s_4'}$,
$\lvec{rs'_2-s_3'}$ and $\rvec{rs'_2-s_3'}$ on
$(\R ^2\times \R ^2)\times SO(3)$ are given by
\[
\begin{array}{l}
\lvec{s}'_5=((0,0),\lvec{E}_3),\;  \rvec{s}_5'=((0,0),\rvec{E}_3),\\
\lvec{rs'_1+s_4'}=((0,r\frac{\partial}{\partial x}),\lvec{E}_2),\;
\rvec{rs'_1+s_4'}=((-r\frac{\partial}{\partial x},0),\rvec{E}_2),\\
\lvec{rs'_2-s_3'}=((0,r\frac{\partial}{\partial
y}),-\lvec{E}_1),\; \rvec{rs'_2-s_3'}=((0,
-r\frac{\partial}{\partial y}),-\lvec{E}_1),
\end{array}
\]
where $\lvec{E}_i$ (respectively, $\rvec{E}_i$) is
the left-invariant (respectively, right-invariant)
vector field on $SO(3)$ induced by $E_i\in \mathfrak{so}(3)$,
for $i\in \{1,2,3\}$. Thus, we deduce the following
system of equations:

\begin{eqnarray*}
\tr\, ((W_1-W_2)E_3)&=&0,\\
rm\left(\frac{x_2-2x_1+x_0}{h^2}\right)+\frac{I}{2h^2}\tr\,
((W_1-W_2)E_2)&=&0,\\
rm\left(\frac{y_2-2y_1+y_0}{h^2}\right)-\frac{I}{2h^2}\tr\,
((W_1-W_2)E_1)&=&0,\\
\frac{x_2-x_1}{h}+\frac{r}{2h}\tr (W_2E_2)+\Omega \frac{y_2+y_1}{2}&=&0,\\
\frac{y_2-y_1}{h}-\frac{r}{2h}\tr (W_2E_1)-\Omega
\frac{x_2+x_1}{2}&=&0
\end{eqnarray*}
where $(x_0, x_1, y_0, y_1; W_1)$ are known. Simplifying we obtain
the following system of equations:
\begin{eqnarray}
\frac{x_2-2x_1+x_0}{h^2}+\frac{I\Omega}{I+mr^2}\label{aqe1}
\frac{y_2-y_0}{2h}&=&0\\
\frac{y_2-2y_1+y_0}{h^2}-\frac{I\Omega}{I+mr^2}
\frac{x_2-x_0}{2h}&=&0\label{aqe2}\\
\tr\, ((W_1-W_2)E_3)&=&0\\
\frac{x_2-x_1}{h}+\frac{r}{2h}\tr (W_2E_2)+\Omega \frac{y_2+y_1}{2}&=&0,\\
\frac{y_2-y_1}{h}-\frac{r}{2h}\tr (W_2E_1)-\Omega
\frac{x_2+x_1}{2}&=&0.
\end{eqnarray}
Now, consider the open subset $U$ of $\R^2\times \R ^2\times SO(3)$
\[
U=(\R^2\times \R^2)\times \set{W\in SO(3)}{W-\tr (W)I_{3\times 3}
\mbox{ is regular}}.
\]
Then, using Corollary \ref{c3.12} (iv), we deduce that the
discrete nonholonomic Lagrangian system $(L'_d,{\mathcal
M}'_c,{\mathcal D}'_c)$ is regular in the open subset $U'$ of
${\mathcal M}'_c$ given by $U'=U\cap {\mathcal M}'_c$.

\begin{figure}[h]
\centering
  \includegraphics[width=5.2cm]{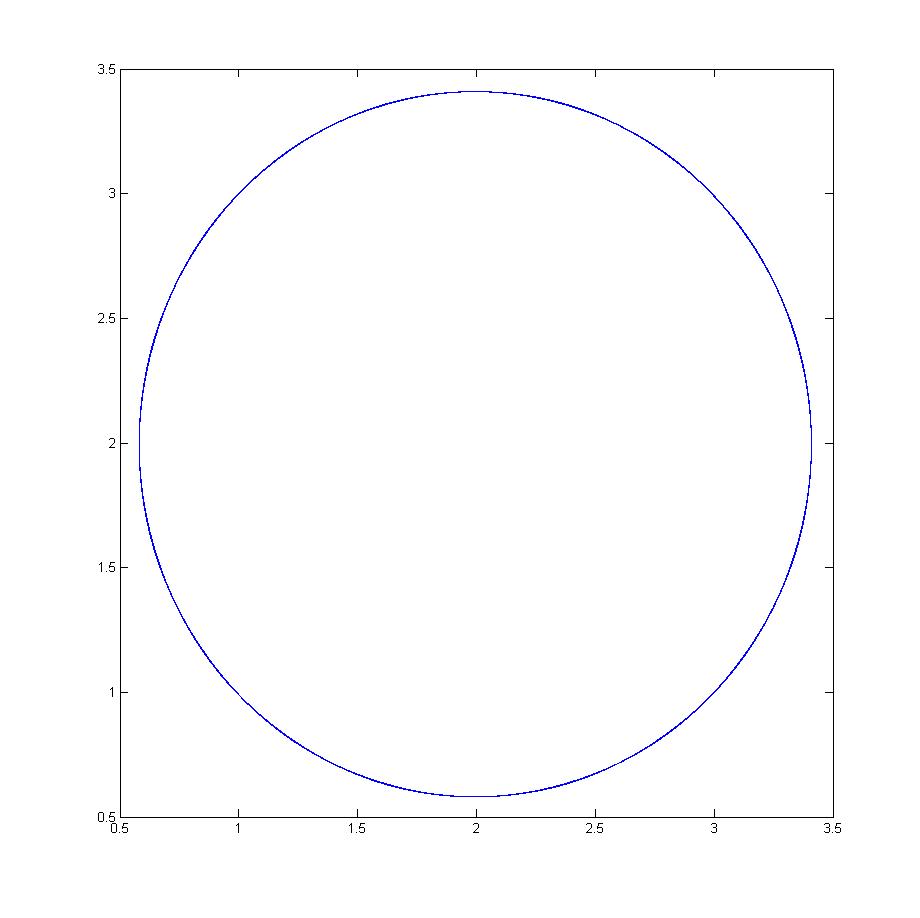}\qquad
  \includegraphics[width=5cm]{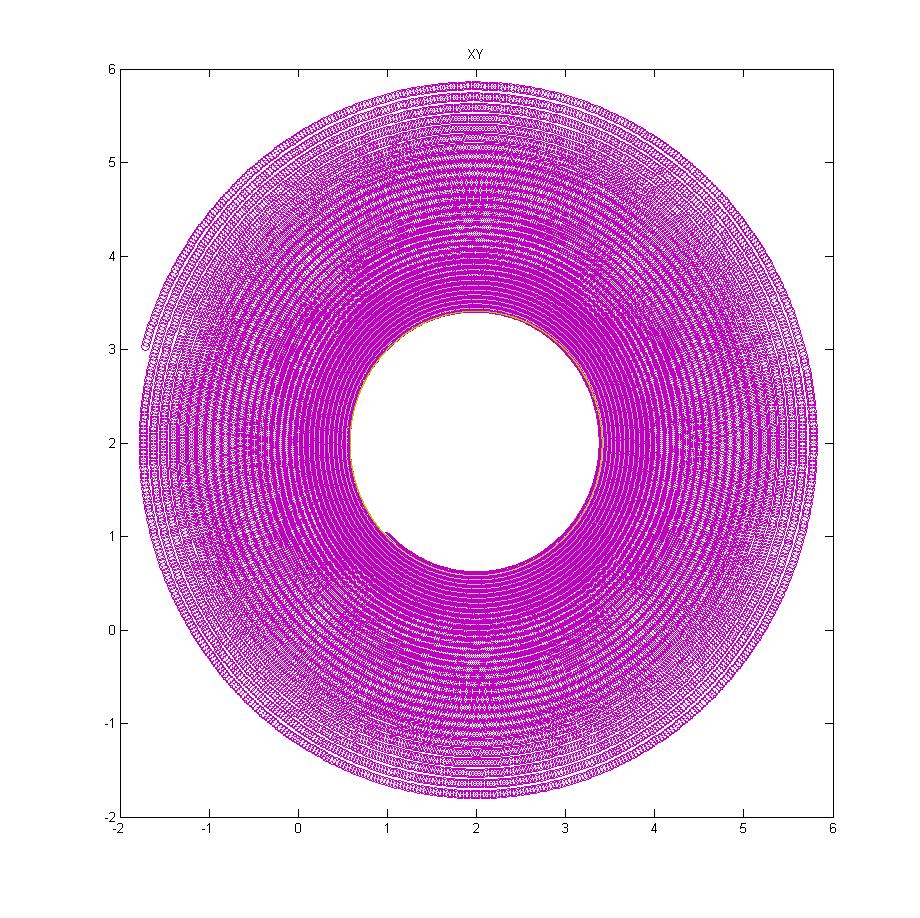}
  \caption{Orbits for the discrete nonholonomic equations of motion
  (left) and a standard numerical method (right)
(initial conditions $x(0)=0.99$, $y(0)=1$, $x(1)=1$, $y(1)=0.99$
and $h=0.01$ after $20000$ steps).}\label{figure1}
\end{figure}

If we denote by $u_k=(x_{k+1}-x_k)/h$ and  $v_k=(y_{k+1}-y_k)/h$,
$k\in\N$ then from  Equations (\ref{aqe1}) and (\ref{aqe2}) we
deduce that
\[
\left(
\begin{array}{c}
u_{k+1}\\
v_{k+1}
\end{array}
\right) =A\left(
\begin{array}{c}
u_{k}\\
v_{k}
\end{array}
\right)=\frac{1}{4+\alpha^2h^2}\left(
\begin{array}{cc}
4-\alpha^2 h^2&-4\alpha h\\
4\alpha h&4-\alpha^2 h^2
\end{array}
\right)\left(
\begin{array}{c}
u_{k}\\
v_{k}
\end{array}
\right)\] or in other terms
\begin{eqnarray*}
x(k+2)&=&\frac{8x(k+1)+(\alpha^2h^2-4)x(k)-4\alpha
h(y(k+1)-y(k))}{\alpha^2h^2+4}\\
y(k+2)&=&\frac{8y(k+1)+(\alpha^2h^2-4)y(k)+4\alpha(x(k+1)-x(k))}{\alpha^2h^2+4};
\end{eqnarray*}
where $\alpha=\frac{I\Omega}{I+mr^2}$. Since $A\in SO(2)$, the
discrete  nonholonomic model predicts that the point of contact of
the ball will sweep out a circle on the table in agreement with
the  continuous model. Figure \ref{figure1}  shows the excellent
behaviour of the proposed numerical method

\subsection{Discrete Chaplygin systems}

Now, we present the theory for a particular (but typical) example
of discrete nonholonomic systems: \emph{discrete Chaplygin
systems}. This kind of systems was considered in the case of the
pair groupoid in \cite{CoSMa}.

For any groupoid $\Gamma\rightrightarrows M$, the map $\chi:
\Gamma \rightarrow M\times M$, $g\mapsto (\alpha(g), \beta(g))$ is
a morphism over $M$  from $\Gamma$ to the pair groupoid $M\times
M$ (usually called the \emph{anchor} of $\Gamma$). The induced
morphism of Lie algebroids is precisely  the anchor $\rho:
E_{\Gamma}\rightarrow TM$ of $E_{\Gamma}$ (the Lie algebroid of
$\Gamma$).

\begin{definition}
A \emph{discrete Chaplygin system}  on the groupoid $\Gamma$ is a
discrete nonholonomic problem $(L_d,{\mathcal
M}_c,{\mathcal D}_c)$ such that
\begin{itemize}
\item[-] $(L_d,{\mathcal M}_c, {\mathcal D}_c)$
is a regular discrete nonholonomic Lagrangian system;
\item[-] $\chi_{{\mathcal M}_c}=\chi\circ i_{{\mathcal M}_c}: {\mathcal M}_c\longrightarrow
M\times M$ is a diffeomorphism;
\item[-] $\rho\circ i_{{\mathcal D}_c}: {\mathcal
D}_c\longrightarrow TM$ is an isomorphism of vector bundles.
\end{itemize}
\end{definition}

Denote by $\tilde{L}_d: M\times M\longrightarrow \R$ the discrete
Lagrangian defined by $\tilde{L}_d=L_d\circ i_{{\mathcal
M}_c}\circ (\chi_{{\mathcal M}_c})^{-1}$.

In the following, we want to express the dynamics on $M\times M$,
by finding relations between de dynamics defined by the
nonholonomic system on $\Gamma$ and $M\times M$.

From our hypothesis, for any vector field $Y\in {\mathfrak X}(M)$
there exists a unique section $X\in \Sec{\tau_{{\mathcal D}_c}}$
such that $\rho\circ i_{{\mathcal D}_c}\circ X=Y$.

 Now, using (\ref{linv}), (\ref{rinv}) and (\ref{LA}), it follows that
\[
T_g \alpha (\rvec{X}(g))=-Y(\alpha(g))\ \hbox{ and }\ T_g \beta
(\lvec{X}(g))=Y(\beta(g))
\]
with some abuse of notation. In other words,
\[
{\mathcal T}_g \chi (X^{(1,0)}(g))=Y^{(1,0)}(\alpha(g), \beta(g))\
\hbox{ and }\ {\mathcal T}_g \chi
(X^{(0,1)}(g))=Y^{(0,1)}(\alpha(g), \beta(g))
\]
for $g\in{\mathcal M}_c$, where ${\mathcal T}\chi :
{\mathcal T}^\Gamma\Gamma \cong V\beta\oplus _\Gamma V\alpha
\to {\mathcal T}^{M\times M}(M\times M)\cong T(M\times M)$
is the prolongation of the morphism $\chi$ given by
\[
({\mathcal T}_g\chi )(X_g,Y_g)=((T_g\alpha )(X_g),(T_g\beta )(Y_g)),
\]
for $g\in \Gamma$ and $(X_g,Y_g)\in {\mathcal T}^\Gamma
_g\Gamma \cong V_g\beta\oplus  V_g\alpha$.

Since $\chi_{{\mathcal M}_c}$ is a diffeomorphism, there exists
a unique $X'_g\in T_g{\mathcal M}_c$ (respectively, $\bar{X}'_g
\in T_g{\mathcal M}_c$) such that
\[
({\mathcal T}_g\chi _{{\mathcal M}_c})(X'_g)=Y^{(1,0)}(\alpha (g),\beta (g))=
(-Y(\alpha (g)),0_{\beta (g)})
\]
(respectively, $({\mathcal T}_g\chi _{{\mathcal M}_c})(\bar{X}'_g)
=Y^{(0,1)}(\alpha (g),\beta (g))=
(0_{\alpha (g)},Y(\beta (g)))$) for all $g\in {\mathcal M}_c$.

Thus,
\[
\begin{array}{l}
X'_g\in T_g{\mathcal M}_c\cap V_g\beta , \quad \rvec{X}(g)-X'_g=Z'_g \in
V_g\alpha \cap V_g\beta, \\
\bar{X}'_g\in T_g{\mathcal M}_c\cap V_g\alpha ,\quad  \lvec{X}(g)-\bar{X}'_g=\bar{Z}'_g \in
V_g\alpha \cap V_g\beta,
\end{array}
\]
for all $g\in {\mathcal M}_c$.

Now, if $(g, h)\in \Gamma_2\cap ({\mathcal M}_c\times {\mathcal
 M}_c)$ then
\begin{eqnarray*}
\lvec{X}(g)(L_d)-\rvec{X}(h)(L_d)&=&
\bar{X}'_g(L_d)+\bar{Z}'_g(L_d) -{X}'_h(L_d)-Z'_h(L_d)\\
&=&\lvec{Y}(\alpha(g),\beta(g))(\tilde{L}_d)
-\rvec{Y}(\alpha(h),\beta(h))(\tilde{L}_d)\\
&&+\bar{Z}'_g(L_d)-Z'_h(L_d).
\end{eqnarray*}
Therefore, if we use the following notation
\[
\begin{array}{l}
(\alpha(g),\beta(g))=(x,y),\qquad (\alpha(h),\beta(h))=(y,z) \\
F_Y^+(x,y)=-\bar{Z}'_{\chi^{-1}_{{\mathcal M}_c}(x, y)}(L_d),\quad
F_Y^-(y,z)=Z'_{\chi^{-1}_{{\mathcal M}_c}(y, z)}(L_d),
\end{array}
\]
then
\begin{eqnarray*}
\lvec{X}(g)(L_d)-\rvec{X}(h)(L_d)&=&
\lvec{Y}(x,y)(\tilde{L}_d)-\rvec{Y}(y,z)(\tilde{L}_d)\\
& & -F_Y^+(x,y)+F_Y^-(y,z).
\end{eqnarray*}
In conclusion, we have proved that $(g,h)$ is a solution of the
discrete nonholonomic Euler-Lagrange equations for the system
$(L_d,{\mathcal M}_c,{\mathcal D}_c)$ if and only if $((x,y),(y,z))$
is a solution of the reduced equations
\[
\lvec{Y}(x,y)(\tilde{L}_d)-\rvec{Y}(y,z)(\tilde{L}_d) =
F_Y^+(x,y)-F_Y^-(y,z),\qquad Y\in \mathfrak X(M).
\]
Note that the above equations are the standard forced discrete
Euler-Lagrange equations (see \cite{mawest}).

\subsubsection{The discrete two wheeled planar mobile
robot}

We now consider a discrete version of the two-wheeled planar
mobile robot \cite{Cort,CoLeMaMa}. The position and orientation of
the robot is determined, with respect a fixed cartesian reference,
by an element $\Omega=(\theta,x,y)\in SE(2)$, that is, a matrix
\[
\Omega = \left(%
\begin{array}{ccc}
  \cos \theta  & -\sin\theta & x \\
  \sin \theta & \cos\theta & y \\
  0 & 0 & 1 \\
\end{array}%
\right).
\]
Moreover, the different
positions of the two wheels are described by elements $(\phi,
\psi)\in \T^2$. Therefore, the configuration space is $SE(2)\times
\T^2$. The system is subjected to three nonholonomic constraints:
one  constraint induced by the condition of no lateral sliding of
the robot and the other two by the rolling conditions of both
wheels.

It is well known that this system is $SE(2)$-invariant and then
the system may be described as a nonholonomic system on the Lie
algebroid $\mathfrak{se}(2)\times T\T^2\to \T ^2$ (see
\cite{CoLeMaMa}). In this case, the Lagrangian is
\begin{eqnarray*}
L&=&\frac{1}{2}\left(J\omega^2+m(v^1)^2+m(v^2)^2+2m_0l\omega v^2+J_2\dot{\phi}^2+J_2\dot{\psi}^2\right)\\
 &=& \frac{1}{2}\tr (\xi \J\xi ^T)+\frac{J_2}{2}\dot{\phi}^2+\frac{J_2}{2}\dot{\psi}^2
 \end{eqnarray*}
 where
 \[
 \xi = \omega\, e+v^1\,e_1+v^2\, e_2=\left(
\begin{array}{ccc}
0&-\omega&v^1\\
\omega&0&v^2\\
0&0& 0
\end{array}
\right)\qquad \hbox{and} \qquad \J=\left(
\begin{array}{ccc}
J/2&0&m_0l\\
0&J/2&0\\
m_0l&0& m
\end{array}
\right)
\]
Here, $m=m_0+2m_1$, where $m_0$ is the mass of the robot without
the two wheels, $m_1$ the mass of each wheel, $J$ its the moment
of inertia with respect to the vertical axis, $J_2$ the axial
moments of inertia of the wheels and $l$ the distance between the
center of mass of the robot and the intersection point of the
horizontal symmetry axis of the robot and the horizontal line
connecting the centers of the two wheels.

The nonholonomic constraints are
\begin{equation}\label{nonhcons}
\begin{array}{rcl}
v^1+\frac{R}{2}\dot{\phi}+\frac{R}{2}\dot{\psi}&=&0,\\
v^2&=&0,\\
\omega+\frac{R}{2c}\dot{\phi}-\frac{R}{2c}\dot{\psi}&=&0,
\end{array}
\end{equation}
determining a submanifold ${\mathcal M}$ of ${\frak se}(2)\times
T\T^2$, where $R$ is the radius of the two wheels and $2c$ the
lateral length of the robot.

In order to discretize the above nonholonomic system, we consider
the Atiyah groupoid $\Gamma =SE(2)\times (\T^2\times \T ^2)
\rightrightarrows \T^2$. The Lie algebroid of $SE(2)\times
(\T^2\times \T ^2) \rightrightarrows \T^2$ is $T\T ^2\times
\mathfrak{se}(2)\to \T ^2$. Then:
\begin{itemize}
\item[-] The discrete Lagrangian $L_d: SE(2)\times (\T^2\times
\T^2)\rightarrow \R$ is given by:
\[
\begin{array}{rcl}
L_d(\Omega_k, \phi_k,\psi_k, \phi_{k+1},
\psi_{k+1})&=&\frac{1}{2h^2}\hbox{Tr } ( (\Omega_k-I_{3\times 3})
\J (\Omega_k-I_{3\times 3})^T)\\
&&+\frac{J_1}{2}\frac{(\Delta\phi_k)^2}{h^{2}}+\frac{J_1}{2}\frac{(\Delta
\psi_k)^2}{h^{2}},
\end{array}
\]
where $I_{3\times 3}$ is the identity matrix,
$\Delta\phi_k=\phi_{k+1}-\phi_k$, $\Delta\psi_k=\psi_{k+1}-\psi_k$
and
\[
\Omega_k=\left(
\begin{array}{lll}
\cos\theta_k&-\sin\theta_k& x_k\\
\sin\theta_k&\cos\theta_k& y_k\\
0&0&1
\end{array}
\right).
\]
We obtain that
\begin{eqnarray*}
L_d&=& \frac{1}{2h^2}\left( m x^2_k+my_k^2
-2lm_0x_k(1-\cos\theta_k)\right.\\&&\left.+2J(1-\cos\theta_k)+2lm_0y_k\sin\theta_k\right)
+\frac{1}{2}{J_1}\frac{(\Delta\phi_k)^2}{h^{2}}+\frac{1}{2}J_1\frac{(\Delta
\psi_k)^2}{h^{2}}.
\end{eqnarray*}

\item[-] The constraint vector subbundle of $\mathfrak{se}(2)\times
T\T^2$ is generated by the sections:
\[
\left\{ s_1=\frac{R}{2}e_1+\frac{R}{2c}e-\frac{\partial}{\partial
\phi}, s_2=\frac{R}{2}e_1-\frac{R}{2c}e-\frac{\partial}{\partial
\psi}\right\}.
\]
\item[-] The continuous constraints of the two-wheeled planar robot are
written in matrix form (see \ref{nonhcons}):
\[
\xi=\left(
\begin{array}{ccc}
0&-\omega&v^1\\
\omega&0& v^2\\
0&0&0
\end{array}
\right)= \left(
\begin{array}{ccc}
0&\frac{R}{2c}\dot{\phi}-\frac{R}{2c}\dot{\psi}&-\frac{R}{2}\dot{\phi}-\frac{R}{2}\dot{\psi}\\
-\frac{R}{2c}\dot{\phi}+\frac{R}{2c}\dot{\psi}&0&0\\
0&0&0
\end{array}
\right)
\]

We discretize the previous constraints using the exponential on
$SE(2)$ (see Section \ref{s4.3.2}) and discretizing the velocities on the right hand side
\[
\Omega_k=\left(
\begin{smallmatrix}
\cos \left(\frac{R}{2c}\Delta{\phi}_k-\frac{R}{2c}\Delta{\psi}_k\right)&
\sin\left(\frac{R}{2c}\Delta{\phi}_k-\frac{R}{2c}\Delta{\psi}_k\right)
&-c\frac{\Delta {\phi}_k+\Delta{\psi}_k}{\Delta
{\phi}_k-\Delta{\psi}_k} \sin \left(\frac{R}{2c}\Delta
{\phi}_k-\frac{R}{2c}\Delta{\psi}_k  \right)\\
-\sin\left(\frac{R}{2c}\Delta{\phi}_k-\frac{R}{2c}\Delta{\psi}_k\right)&\cos
\left(\frac{R}{2c}\Delta{\phi}_k-\frac{R}{2c}\Delta{\psi}_k\right)&c\frac{\Delta
{\phi}_k+\Delta{\psi}_k}{\Delta {\phi}_k-\Delta{\psi}_k}\left( 1-
\cos \left(\frac{R}{2c}\Delta
{\phi}_k-\frac{R}{2c}\Delta{\psi}_k  \right)\right)\\
0&0&1
\end{smallmatrix}
\right)
\]
if $\Delta{\phi}_k\not= \Delta{\psi}_k$ and
\[
\Omega_k=\left(
\begin{array}{ccc}
1&0 &-R \Delta{\phi}_k\\
0&1&0\\
0&0&1
\end{array}
\right)
\]
if $\Delta{\phi}_k= \Delta{\psi}_k$.

 Therefore, the constraint submanifold ${\mathcal M}_c$ is defined as
\begin{eqnarray}
\theta_k&=&-\frac{R}{2c}\Delta{\phi}_k+\frac{R}{2c}\Delta{\psi}_k\label{co1}\\
x_k&=&-c\frac{\Delta {\phi}_k+\Delta{\psi}_k}{\Delta
{\phi}_k-\Delta{\psi}_k} \sin \left(\frac{R}{2c}\Delta
{\phi}_k-\frac{R}{2c}\Delta{\psi}_k  \right)\label{co2}\\
y_k&=&c\frac{\Delta {\phi}_k+\Delta{\psi}_k}{\Delta
{\phi}_k-\Delta{\psi}_k}\left( 1- \cos \left(\frac{R}{2c}\Delta
{\phi}_k-\frac{R}{2c}\Delta{\psi}_k  \right)\right)\label{co3}
\end{eqnarray}
if $\Delta{\phi}_k\not= \Delta{\psi}_k$ and $\theta_k=0$,
$x_k=-R\Delta{\phi}_k$ and $y_k=0$ if $\Delta{\phi}_k=
\Delta{\psi}_k$.
\end{itemize}
We have that the discrete nonholonomic system $(L_d,{\mathcal
M}_c, {\mathcal D}_c)$ is reversible. Moreover, if $\epsilon
_\Gamma : \T^2\to SE(2) \times (\T ^2\times \T ^2)$ is the
identity section of the Lie groupoid $\Gamma =SE(2)\times
(\T^2\times \T^2)$ then it is clear that
\[
\epsilon _\Gamma (\T^2)=\{ I_{3\times 3} \}\times \Delta_{\T
^2\times \T ^2}\subseteq {\mathcal M}_c.
\]
Here, $\Delta_{\T ^2\times \T ^2}$ is the diagonal in $\T ^2\times
\T ^2$. In addition, the system $(L_d, {\mathcal M}_c, {\mathcal
D}_c)$
 is regular in a neighborhood $U$ of the submanifold
 $\epsilon _\Gamma (\T^2)=\{ I_{3\times 3} \}\times
\Delta_{\T ^2\times \T ^2}$ in ${\mathcal M}_c$. Note that
\[
T_{(I_{3\times 3},\phi _1,\psi _1,\phi _1,\psi _1)}
{\mathcal M}_c\cap E_\Gamma (\phi _1,\psi _1)=
{\mathcal D}_c (\phi _1,\psi _1),
\]
for $(\phi _1,\psi _1)\in \T ^2$, where
$E_\Gamma =\mathfrak{se}(2)\times T\T ^2$ is
the Lie algebroid of the Lie groupoid $\Gamma
= SE(2)\times (\T ^2\times \T ^2)$.

On the other hand, it is easy to show that the
system $(L_d,U,{\mathcal D}_c)$ is a discrete
Chaplygin system.

The reduced Lagrangian on $\T^2\times \T^2$ is
\[
\tilde{L}_d
=\left\{
\begin{array}{l}
\displaystyle \frac{1}{h^{2}} (mc^2  ( \frac{\Delta
{\phi}_k+\Delta{\psi}_k}{\Delta {\phi}_k-\Delta{\psi}_k}  )^2
( 1- \cos ( \frac{R}{2c}\Delta {\phi}_k-\frac{R}{2c}\Delta{\psi}_k  )  ) \\
\displaystyle +J  (  1-\cos   ( \frac{R}{2c}\Delta
{\phi}_k-\frac{R}{2c}\Delta{\psi}_k ))) +\frac{1}{2}J_1
\frac{(\Delta \phi _k)^2}{h^{2}}+\frac{1}{2} J_{1}\frac{(\Delta
\psi _k)^2}{h^{2}} \quad \hbox{if}\quad \Delta{\phi}_k\not=
\Delta{\psi}_k
\\[15pt]
\displaystyle (J_1+\frac{mR^2}{2})\frac{(\Delta\phi_k)^2}{h^{2}},
\quad \hbox{if}\quad \Delta{\phi}_k= \Delta{\psi}_k
\end{array}
\right.
\]
The discrete nonholonomic equations are:
\begin{eqnarray*}
\lvec{s_1}\big|_{(\Omega_1,\phi_1, \psi_1, \phi_2, \psi_2)}(L_d)-
\rvec{s_1}\big|_{(\Omega_2,\phi_2, \psi_2, \phi_3, \psi_3)}(L_d)&=&0\\
\lvec{s_2}\big|_{(\Omega_1,\phi_1, \psi_1, \phi_2, \psi_2)}(L_d)-
\rvec{s_2}\big|_{(\Omega_2\phi_2, \psi_2, \phi_3, \psi_3)}(L_d)&=&0
\end{eqnarray*}
These equations in coordinates are:
\begin{eqnarray}
&&2J_1(\phi_3-2\phi_2+\phi_1)=lRm_0(\cos\theta_2+\cos\theta_1)+\frac{JR}{c}(\sin\theta_2-\sin\theta_1)\nonumber\\
&&-\frac{R\cos\theta_1}{c}(lm_0y_1+cmx_1)
+\frac{R\sin\theta_1}{c}(lm_0x_1-cmy_1)\nonumber\\&&
+\frac{R}{c}(cmx_2+lm_0(y_2-2c))\label{equa1}\\
&&\hspace{-0.5cm}2J_1(\psi_3-2\psi_2+\psi_1)=lRm_0(\cos\theta_2+\cos\theta_1)-\frac{JR}{c}(\sin\theta_2-\sin\theta_1)\nonumber\\
&&+\frac{R\cos\theta_1}{c}(lm_0y_1-cmx_1)
-\frac{R\sin\theta_1}{c}(lm_0x_1+cmy_1)\nonumber\\&&
+\frac{R}{c}(cmx_2-lm_0(y_2+2c))\label{equa2}
\end{eqnarray}
Substituting constraints (\ref{co1}), (\ref{co2}) and (\ref{co3})
in Equations (\ref{equa1}) and (\ref{equa2}) we obtain a set of
equations of the type $0=f_1(\phi_1, \phi_2, \phi_3, \psi_1,
\psi_2, \psi_3)$ and $0=f_1(\phi_1, \phi_2, \phi_3, \psi_1,
\psi_2, \psi_3)$ which are the reduced equations of the Chaplygin
system.

\section{Conclusions and Future Work}
In this paper we have elucidated the geometrical framework for
nonholonomic discrete Mechanics on Lie groupoids. We have proposed
 discrete nonholonomic equations that are general enough to
produce practical integrators for continuous nonholonomic systems
(reduced or not). The geometric properties related with these
equations have been completely studied and the applicability of
these developments has been stated in several interesting
examples.

Of course,  much work remains to be done to clarify the nature of
discrete nonholonomic mechanics. Many of this future work was
stated in \cite{McLPerl} and, in particular, we emphasize:
\begin{itemize}
\item[-] a complete backward error analysis which explain the very
good energy behavior showed in examples or the preservation of a
discrete energy (see \cite{FZ});
\item[-] related with the previous question, the construction of a
discrete exact model for a continuous nonholonomic system (see
\cite{IMMM,mawest,McLPerl});
\item[-] to study discrete nonholonomic systems which preserve a volume form on the constraint
surface mimicking the continuous case (see, for instance,
\cite{FeJo,ZeBl} for this last case);
\item[-] to analyze the discrete hamiltonian framework and the
construction of  integrators depending on different
discretizations;
\item[-] and the construction of a discrete nonholonomic connection in the
case of Atiyah groupoids (see \cite{TL,groupoid}).
\end{itemize}

Related with some of the previous questions, in the conclusions of
the paper of R. McLachlan and M. Perlmutter \cite{McLPerl}, the
authors raise the question of the possibility of the definition of
generalized constraint forces dependent on all the points
$q_{k-1}$, $q_k$ and $q_{k+1}$ (instead of just $q_k$) for the
case of the pair groupoid. We think that the discrete nonholonomic
Euler-Lagrange equations can be generalized to consider this case
of general constraint forces that, moreover, are closest to the
continuous model (see \cite{LeMaSa,McLPerl}).

\end{document}